\documentclass[11pt]{amsart}
\usepackage{hyperref}
\usepackage{latexsym}
\usepackage{amsmath}
\usepackage{amsfonts}
\usepackage{amssymb}
\usepackage{amsthm}

\usepackage{enumerate}

\usepackage[a4paper, twoside=false, vmargin={2cm,3cm}, includehead]{geometry}

  \newcommand{\LE}{\mathcal{LE}}
   \renewcommand{\H}{\mathcal{H}}
  
 \newcommand{\cB}{\mathcal{X}}
\newcommand{\cX}{\mathcal{X}}

\renewcommand{\P}{\mathcal{P}}
\newcommand{\G}{\mathcal{G}}
\newcommand{\X}{\mathcal{X}}
\newcommand{\cI}{\mathcal{I}}

\newcommand{\cZ}{\mathcal{Z}}
\newcommand{\gG}{\Gamma}
\newcommand{\C}{\mathbb{C}}
\newcommand{\T}{\mathbb{T}}
\newcommand{\Q}{\mathbb{Q}}
\newcommand{\R}{\mathbb{R}}
\newcommand{\E}{\mathbb{E}}
\newcommand{\N}{\mathbb{N}}
\newcommand{\bbZ}{\mathbb{Z}}
\newcommand{\Z}{\mathbb{Z}}

\newcommand{\norm}[1]{\left\Vert #1\right\Vert}
\newcommand{\nnorm}[1]{\lvert\!|\!| #1|\!|\!\rvert}
\theoremstyle{plain}
\newtheorem{theorem}{Theorem}[section]

\newtheorem*{theoremA'}{Theorem A'}
\newtheorem*{theoremB'}{Theorem B'}
\newtheorem*{theoremC'}{Theorem C'}

\newtheorem{problem}{Problem}
\newtheorem*{problem*}{Problem}
\newtheorem*{theorem*}{Theorem}

\newtheorem*{problem1*}{Special Case of Problem~\ref{Conj:StructMultiGeneral}}
\newtheorem*{problem2*}{Special Case of Problem~\ref{Conj:CondConvGeneral}}
\newtheorem*{problem2**}{Variant of Special Case of Problem~\ref{Conj:CondConvGeneral}}
\newtheorem*{problem10*}{Special Case of Problem~\ref{Conj:ComplexityBound}}
\newtheorem*{problem15*}{Special Case of Problem~\ref{Conj:CommCharPairIndep}}
\newtheorem*{problem19*}{Special Case of Problem~\ref{Conj:ConvNonComm}}
\newtheorem*{problem20*}{Special Case of Problem~\ref{Conj:RecNonComm}}
\newtheorem*{problem21*}{Special Case of Problem~\ref{Conj:RecNonCommEasy}}
\newtheorem*{problem25*}{Special Case of Problem~\ref{C:ConjChar}}
\newtheorem*{problem31*}{Special Case of Problem~\ref{Conj:ConvNonComm2}}

\newtheorem*{Correspondence1}{Furstenberg Correspondence Principle}

\theoremstyle{definition}
\newtheorem{definition}[theorem]{Definition}
\newtheorem{example}{Example}

\theoremstyle{remark}

\begin{document}
\title[\tiny{Some open problems on multiple ergodic averages}]{Some open problems on multiple ergodic averages}

\author{Nikos Frantzikinakis}
\address[Nikos Frantzikinakis]{University of Crete, Department of Mathematics and Applied Mathematics, Voutes University Campus, Heraklion 71003, Greece} \email{frantzikinakis@gmail.com}

\begin{abstract}
We survey some recent developments and give a list of open problems regarding
multiple recurrence and convergence phenomena of $\Z^d$ actions in ergodic theory
and related applications in combinatorics and number theory.
\end{abstract}



\maketitle


\section{Introduction}
\subsection{Objective}
Ergodic theory has a long history of interaction with other
mathematical fields and in particular with combinatorics and number theory.
The seminal work of  H.~Furstenberg \cite{Fu77}, where an ergodic proof of the theorem of Szemer\'edi \cite{Sz75} on arithmetic progressions was given,
linked problems in  ergodic theory, combinatorics,  and number theory,
and  provided an ideal ground for cross-fertilization. In combinatorics, it has produced several far reaching extensions of the  theorem of Szemer\'edi
 some of  which  still have no proof  that avoids ergodic theory.
 In number theory, it provided some key ideas and tools in the proof of Green and Tao  \cite{GT08} that the primes contain arbitrarily long arithmetic progressions and subsequent extensions  of this result.
 On the other direction, the field of ergodic
theory has tremendously benefited as well, since the problems of
combinatorial and number-theoretic nature gave a boost to
the in depth study of several interesting recurrence and convergence phenomena that would have otherwise been ignored.

The connecting link between combinatorics and ergodic theory is that
regularity properties of sets of integers with positive density
correspond to multiple recurrence properties of measure preserving
systems.
One   establishes these  recurrence
properties  by  analyzing  the limiting behavior of some closely related
multiple ergodic averages. The study of these ergodic averages has developed into a central part of ergodic theory and has generated  tools and ideas that have found applications in areas outside ergodic theory.
The purpose of these notes is to survey some of these developments and give
 a list of  open problems  on three closely related topics:
\medskip

\begin{enumerate}
\item  The limiting behavior of multiple ergodic averages.
\medskip

\item
Multiple recurrence properties of measure preserving systems.

\medskip
\item Problems of arithmetic nature regarding, for example,  the existence of  patterns on sets of integers
with positive   upper density\footnote{The {\em
upper  density} $\bar{d}(E)$ of a set $E\subset\Z^d$ is defined by
$\bar{d}(E) \mathrel{\mathop:}= \limsup_{N \to\infty}\frac{|E\cap
[-N,N]^d|}{|[-N,N]^d |}$.}.
\end{enumerate}
\medskip





 The list  of problems is greatly influenced by my personal interests  and is by no means meant to be a comprehensive list
 of open problems in the area widely known as ergodic Ramsey theory. Almost exclusively,
  problems related to $\Z^d$ actions
   are considered and even within this confined class there are a few important topics not touched upon.
 For material and a list of problems that goes  beyond the scope of this set of notes   we refer the reader to the
 survey articles
   \cite{Be96, Be06a, Be06b} 
   and the references therein.

 A first version of these notes was posted in 2011. Since then, several of the recorded problems have been partially or completely solved and new techniques and problems have surfaced.
 This updated version is meant to address these developments.


 \subsection{The general framework}
 We are given a
 measure space  $(X,\X,\mu)$ with $\mu(X)=1$,
 and invertible
measure preserving transformations $T_1,\ldots,T_\ell\colon X\to X$
that commute, that is, satisfy $T_iT_j=T_jT_i$ for all $i,j\in \{1,\ldots, \ell\}$.
In all applications we are interested in we can further assume that the
 measure space  is {\em Lebesgue}, that is,
  $X$  can be given the structure of a Polish space  (i.e. metrizable, separable, complete) such that $\X$  is its Borel $\sigma$-algebra.
We call $(X,\X, \mu, T_1,\ldots, T_\ell)$ a {\em measure preserving system}, or simply,
a {\em system}.  We are also given
bounded measurable functions $f_1,\ldots, f_\ell\colon X \to \C$,
and sequences  $a_1,\ldots,a_\ell \colon \N\to \Z$.

The first family of problems we are interested in concerns the study of the limiting behavior of the so called multiple ergodic averages
\begin{equation}\label{E:GeneralMulties}
\frac{1}{N}\sum_{n=1}^N T_1^{a_1(n)}f_1\cdots
T_\ell^{a_\ell(n)}f_\ell
\end{equation}
where  $Tf\mathrel{\mathop:}=f\circ T$ and
$T^k\mathrel{\mathop:}=T\circ \cdots \circ T$. One would like to
know whether these averages converge as $N\to \infty$ (in $L^2(\mu)$
or pointwise), find some structured factors that control their
limiting behavior (called characteristic factors), and if possible,
find a formula, or a usable way to extract information, for the
limiting function. When $\ell=1$, such problems have been studied
extensively and in several cases solved even for pointwise convergence (see the survey paper
\cite{RW95}
  for a variety of
 related  results). Our main concern here is to study the averages \eqref{E:GeneralMulties} when $\ell\geq 2$. To get manageable problems,
 one typically  restricts
the class of eligible sequences and usually assumes that they are polynomial sequences,
sequences arising from smooth functions,
 sequences related to the prime numbers,
 or random sequences of integers.
  We typically  also assume that the
transformations commute, or,  to get started, that they are all
equal. On the other hand,  because of the nature of the implications
in combinatorics that we are interested in,
 it is not desirable to impose assumptions  on the  structure of each individual measure  preserving transformation.
 In several cases, the steps taken in order to  attack such problems include  $(i)$
elementary uniformity  estimates, $(ii)$ ergodic  structural results,
and $(iii)$ analysis of systems with special structure (often algebraic).
We discuss  these steps in more detail in subsequent sections.


The  second family of problems concerns the study of expressions of the form
\begin{equation}\label{E:GeneralRecurrence}
\mu(A\cap T_1^{a_1(n)}A\cap \cdots \cap T_\ell^{a_\ell(n)}A)
\end{equation}
where $A\in \X$ has positive measure.
One wants to know whether such expressions are positive  for some  $n\in \N$, or even better,
for lots of $n\in \N$ (for instance on the average), and if possible, get some explicit lower bound
that depends only  on the measure of the set $A$ and on $\ell$ (in some cases this is going to be of the form $(\mu(A))^{\ell+1}$). Such multiple recurrence results are often obtained by carrying out an in depth
analysis of the limiting behavior of the averages \eqref{E:GeneralMulties}. Usually they are not hard
 if an explicit formula of the limiting function is known, but they can be very tricky in the absence of such a
formula, even when we work with very special systems of algebraic nature.

Concerning the third family of problems,  and restricting ourselves
to subsets of $\Z$,  one is interested to know, for example, whether
every set of integers with positive upper density
contains patterns
of the form
$$
m, \ m+a_1(n),  \ldots, \ m+a_\ell(n)
$$
for some (or lots of) $m\in \Z$ and $ n\in \N$.
Using a correspondence principle of H.~Furstenberg (see Section~\ref{SS:correspondence}), one can
translate such statements to multiple recurrence statements in
ergodic theory; an equivalent problem is
  then to show that the expressions \eqref{E:GeneralRecurrence} are positive for some $n\in \N$ when all the
   measure preserving  transformations $T_1,\ldots, T_\ell$  are equal. Similar questions can be asked on higher dimensions and concern   patterns that can be
   found on  subsets of $\Z^d$ with positive upper density. Such questions
correspond to multiple recurrence statements when the  transformations $T_1,\ldots, T_\ell$  commute.
 This approach  was initiated  forty years ago by H.~Furstenberg  who gave  an
 alternative proof of Szemer\'edi's theorem using ergodic theory in his foundational article \cite{Fu77}. Subsequently,   H.~Furstenberg and Y.~Katznelson  gave the
 first proof
of the multidimensional Szemer\'edi theorem \cite{FuK79} and the
density Hales-Jewett theorem \cite{FuK91}, and  V.~Bergelson and
A.~Leibman  proved a polynomial extension of Szemer\'edi's theorem
\cite{BL96} (currently no proof that avoids ergodic theory is known
for this result). And the story does not end there, in the last two
decades new powerful tools in ergodic theory were developed and
used, and are currently being used, in order to prove several other deep
results in density Ramsey theory. The reader will find  a variety of  such
applications in subsequent sections and the extended bibliography
section. Several additional applications can
be found in   the survey articles  \cite{Be96, Be06a, Be06b} and the references therein.


An additional motivation for studying such problems  has to do
with potential implications
in number theory, in particular, there is a connection to  problems of finding patterns in the set of the prime numbers. Knowing that every
 set of integers with positive upper density contains patterns of a certain sort could be  an important  first  step towards proving an analogous result for the set of   primes. This idea originates from work of   B.~Green and T.~Tao   \cite{GT08}, where it was used to show
that the primes contain arbitrarily long arithmetic progressions. It was also subsequently  used by T.~Tao and T.~Ziegler \cite{TZ08} to show that the primes contain arbitrarily long polynomial progressions.

\bigskip
\noindent {\bf Acknowledgments.} I would like to thank
B.~Kra and E.~Lesigne   for helpful comments on a previous version of this article and
the referees for valuable comments.

\smallskip

\section{Some useful  tools and observations}

\subsection{Furstenberg correspondence principle}\label{SS:correspondence}
In order to reformulate
statements in combinatorics as multiple recurrence statements in
ergodic theory, we  use  the following correspondence principle of
Furstenberg \cite{Fu77, Fu81a} (the formulation given is from \cite{Be87b}):
\begin{Correspondence1} Let $d,\ell\in \N$, $E\subset \Z^d$ be a
  set of integers, and  ${\bf v}_1,\ldots, {\bf v}_\ell\in \Z^d$.  Then there exist a system $(X,\X,\mu, T_1,\ldots, T_\ell)$ and a set
  $A\in\X$, with $\mu(A)=\bar{d}(E)$, and such that
  \begin{equation}\label{E:correspondence}
    \bar{d}(E\cap (E-n_1{\bf v}_1)\cap\dots\cap
    (E-n_\ell{\bf v}_\ell))\geq \mu(A\cap T_1^{n_1}A\cap\dots \cap T^{n_\ell}A),
  \end{equation}
  for every  $n_1,\ldots,n_\ell\in\Z$. Furthermore, if ${\bf v}_1=\cdots={\bf v}_\ell$, one can take
  $T_1=\cdots=T_\ell$.
\end{Correspondence1}
Let $a_1,\ldots,a_\ell\colon \N\to \Z$ be a collection of sequences.  Using the previous principle, we deduce that in order to show that every set of integers with positive upper density contains patterns of the form $m, m+a_1(n), \ldots, m+a_\ell(n)$ for some $m,n\in\N$,
 it suffices
to show that for  every system  $(X,\X,\mu, T)$ and  set $A\in \X$ with $\mu(A)>0$ one has
$$
\mu(A\cap T^{a_1(n)}A\cap\dots \cap T_\ell^{a_\ell(n)}A)>0
$$
for some $n\in\N$.  It is often more convenient to  verify that  a functional variant of this multiple recurrence property holds on average, namely, that for every $f\in L^\infty(\mu)$  that is non-negative and not identically zero one has
\begin{equation}\label{E:positive}
\liminf_{N\to\infty}\frac{1}{N}\sum_{n=1}^N\int f\cdot T^{a_1(n)}f\cdot \ldots \cdot T^{a_\ell(n)}f\, d\mu >0.
\end{equation}
Thus, one is reduced to study the limiting behavior (in $L^2(\mu)$)
of the averages
$$
\frac{1}{N}\sum_{n=1}^N T^{a_1(n)}f\cdots T^{a_\ell(n)}f
$$
in a depth that is sufficient to verify the positiveness property \eqref{E:positive}.

\subsection{Characteristic factors}\label{SS:2.2} A notion  that underlies the
study of the limiting behavior of several multiple ergodic averages
is that of the \emph{characteristic factors}.   Implicit use of
this notion was already made by
H.~Furstenberg in \cite{Fu77}, but the term ``characteristic factor'' was coined  in a paper
of H.~Furstenberg and B.~Weiss \cite{FuW96}.

 Given a system $(X,\X,\mu, T_1,\ldots, T_\ell)$, we say that the
sub-$\sigma$-algebras $\X_1,\ldots,\X_\ell$ of $\X$ are \emph{characteristic factors} for the averages
\begin{equation}\label{E:GeneralMulties'}
A_N(f_1,\ldots,f_\ell)\mathrel{\mathop:}=\frac{1}{N}\sum_{n=1}^N
T_1^{a_1(n)}f_1\cdots T_\ell^{a_\ell(n)}f_\ell
\end{equation}
if the following two conditions hold:
 \begin{itemize}
\item  $\X_i$ \ is \ $T_i$-invariant \ for \ $i=1,\ldots, \ell$,
\item we have \
    $
    A_N(f_1,\ldots,f_\ell)-A_N(\tilde{f}_1,\ldots,\tilde{f}_\ell)\to^{L^2(\mu)} 0,
    $
 \  where \ $\tilde{f}_i\mathrel{\mathop:}=\E(f_i|\X_i)$ \ for \ $i=1,\ldots,
 \ell$, for all $f_1,\ldots, f_\ell \in L^\infty(\mu)$.
\end{itemize}
Equivalently, the second condition states that  if $\E(f_i|\X_i)=0$ for some $i\in \{1,\ldots,\ell\}$, then
  $A_N(f_1,\ldots,f_\ell)\to^{L^2(\mu)} 0$.
If in addition one has  $\X_1=\cdots=\X_\ell$, then we call this common sub-$\sigma$-algebra a   \emph{characteristic factor}
for the averages  \eqref{E:GeneralMulties'}.
 When analyzing the limiting behavior of multiple ergodic averages, an important intermediate goal is to
produce characteristic factors that are   as simple  as
possible, and typically simple for our purposes  means   that the
corresponding factor systems have  very special (often algebraic) structure. For instance, it can be shown (see \cite{Fu77})
that for ergodic systems $(X,\X,\mu,T)$,
the Kronecker factor $\mathcal{K}_T$, which is induced by the eigenfunctions of $T$, is a characteristic factor for the averages
$$
\frac{1}{N}\sum_{n=1}^N T^{an}f_1\cdot T^{bn}f_2,
$$
whenever $a,b$ are distinct integers. It is also  well known  that the system $(X,\mathcal{K}_T,\mu,T)$ is measure theoretically isomorphic to an ergodic rotation on a compact Abelian group with the Haar measure.



\subsection{Uniformity seminorms}
One way  to find simple characteristic factors
 for the averages \eqref{E:GeneralMulties'} is to  show that their $L^2(\mu)$-norm
 is controlled by the so called
\emph{Gowers-Host-Kra  uniformity seminorms}. Similar  seminorms were   first
introduced in combinatorics by T.~Gowers \cite{Gow01} and their
ergodic counterpart (which is more relevant for our study) was introduced
by B.~Host and B.~Kra \cite{HK05a}. For a given   ergodic system
$(X,\X,\mu,T)$ and function $f\in L^\infty(\mu)$, they are defined
as follows:
\begin{gather*}
\label{E:first}
 \nnorm{f}_{1}\mathrel{\mathop:}=\Big| \int f \ d\mu\Big|\ ;\\
\label{E:recur}
\nnorm f_{k+1}^{2^{k+1}} \mathrel{\mathop:}=\lim_{N\to\infty}\frac{1}{N}\sum_{n=1}^{N}
\nnorm{\bar{f}\cdot T^nf}_{k}^{2^{k}}.
\end{gather*}
It is shown in \cite{HK05a} that for every $k\in \N$  the  above limit
 exists, and $\nnorm{\cdot}_k$, thus defined, is a seminorm on
$L^\infty(\mu)$. For non-ergodic systems the seminorms  can be
similarly defined, the only difference is that $\nnorm{\cdot}_1$ is
defined by  $\nnorm{f}_{1}\mathrel{\mathop:}=\norm{\int f \
d\mu_x}_{L^2(\mu)}$, where $\mu=\int \mu_x \ d\mu(x)$ is the ergodic
decomposition of the measure $\mu$ with respect to $T$. If further
clarification is needed, we  write  $\nnorm{\cdot}_{k,\mu}$, or
$\nnorm{\cdot}_{k,T}$. We remark that if a measure preserving system
is \emph{weak mixing} (meaning,  the product system $(X\times X,
\mu\times \mu, T\times T)$ is ergodic),  then $\nnorm{f}_k=|\int f\ d\mu|$ for every $k\in \N$.

\subsection{Nilsystems and nilsequences}\label{SS:nil} The analytic part of the modern  theory of characteristic factors requires familiarity with  variants of the several  classical Fourier analysis results
that apply to functions defined on  general nilmanifolds. In this non-Abelian setup the role of rotations on the circle play the nilsystems and the role
of complex exponential sequences play the nilsequences.
\subsubsection{Nilsystems}  A {\em $k$-step nilmanifold} is a homogeneous space $X=G/\Gamma$, where  $G$ is a $k$-step nilpotent Lie group
and $\Gamma$ is a discrete cocompact subgroup of $G$.
 A {\it $k$-step  nilsystem} is
 a system of the form $(X, \G/\gG, m_X, T_a)$,  where $X=G/\Gamma$ is a $k$-step nilmanifold, $a\in G$,  $T_a\colon X\to X$ is defined by  $T_{a}(g\gG) \mathrel{\mathop:}= (ag) \gG$, $g\in G$,  $m_X$  is
the normalized Haar measure on $X$, and  $\G/\gG$  is the completion
of the Borel $\sigma$-algebra of $G/\gG$.

 Examples of nilsystems include all rotations on  compact Abelian Lie groups,
 and more generally, every  nilpotent affine transformation on a compact
 Abelian Lie group is measure theoretically isomorphic to
 a nilsystem (see Example 1).
But these examples do not cover all the possible nilsystems (see
Example 2).

\begin{example}
\label{ex:one} On the space $G=\bbZ\times\R^2$,
define multiplication as follows: \\
if $g_1=(m_1,x_1,y_1)$ and $g_2=(m_2,x_2,y_2)$, let
$$
g_1\cdot g_2=(m_1+m_2,x_1+x_2, y_1+y_2+m_1x_2).
$$ Then $G$ is a $2$-step
nilpotent group and  the discrete subgroup $\gG=\bbZ^3$ is
cocompact. If  $a=(m,\alpha ,\beta)$. It turns out that  the nilrotation $x\mapsto ax$ defined on $X=G/\Gamma$ is measure theoretically
isomorphic to the unipotent affine transformation $S \colon
\mathbb{T}^2\to\mathbb{T}^2$ (with the Haar measure $m_{\T^2}$) defined by
$$
S(x,y)=(x+\alpha,y+m x+\beta), \quad x,y\in \T.
$$
\end{example}

\begin{example}
\label{ex:two} On the space $G=\R^3$,
define multiplication as follows: \\
if $g_1=(x_1,y_1,z_1)$ and $g_2=(x_2,y_2,z_2)$, let
$$
g_1\cdot g_2=(x_1+x_2,y_1+y_2, z_1+z_2+x_1y_2).
$$ Then $G$ is a $2$-step
nilpotent group and  the discrete subgroup $\gG=\bbZ^3$ is
cocompact. Let $a=(\alpha,\beta,0)$, where $\alpha,\beta\in [0,1)$ are
linearly independent.  It turns out that  the nilrotation $x\mapsto ax$ defined on $X=G/\Gamma$  is  measure theoretically isomorphic to the
skew product transformation $S \colon \mathbb{T}^3\to\mathbb{T}^3$
that has  the form
$$
S(x,y)=(x+\alpha,y+\beta, z+f(x,y)), \quad x,y\in \T,
$$
where  $f\colon \mathbb{T}^2\to\mathbb{T}$ is defined by
$$
f(x,y)= (x+\alpha) [y +\beta] - x [y] -\alpha y, \quad x,y\in \T.
$$
It can be shown that  the system induced by $S$  is \emph{not} measure theoretically isomorphic to a
unipotent affine transformation on some finite dimensional torus.
\end{example}

\subsubsection{Nilsequences}
Following \cite{BHK05} we define a \emph{basic $k$-step
nilsequence} to be  a sequence of the form $(F(b^nx))$ where $X=G/\Gamma$ is a $k$-step nilmanifold,
 $b\in G$,
$x\in X$, and $F\in C(X)$. If in the previous definition we allow $F$ to be Riemann-integrable we call $(F(b^nx))$ a \emph{basic generalized  $k$-step
nilsequence}. A {\em   $k$-step nilsequence} is a uniform limit of {\em basic   $k$-step nilsequences} and similarly we define
{\em generalized   $k$-step nilsequences}.
 We remark that in \cite{BHK05} and subsequent work, the notion of a generalized nilsequence is not used, but for the purpose of formulating some problems later  in this article it seems safer to extend the class of eligible sequences.

More generally, if  $b_1,\ldots,b_\ell\in G$ commute,  $x\in X$, and
$F\in C(X)$ (or $F$ is  Riemann integrable on $X$),
 we call the sequence $(F(b_1^{n_1}\cdots b_\ell^{n_\ell} x))$ a
 \emph{basic (generalized) $k$-step nilsequence in $\ell$-variables}.
A \emph{$k$-step (generalized) nilsequence in $\ell$-variables}, is a uniform limit of basic $k$-step (generalized) nilsequences in $\ell$-variables.

As is easily verified,
the collection of $k$-step (generalized) nilsequences in $\ell$-variables, with the topology of uniform convergence, forms a closed algebra.

 One can verify that for  $\alpha,\beta\in \R$, the sequence    $(e^{i (n^2\alpha+n\beta)})$ is a basic $2$-step nilsequence defined on the nilmanifold given in Example $1$ above.
 Also, for  $\alpha,\beta\in \R$, the sequence  $(e^{i [n\alpha]n\beta})$ is a basic generalized $2$-step
 nilsequence defined on the nilmanifold given by the direct product of the nilmanifolds
 in Examples $1$ and $2$ above. It can be shown (see \cite{BL07}) that if $\alpha$ and $\beta$ are rationally independent, then the sequence $(e^{i [n\alpha]n\beta})$ is asymptotically orthogonal to all sequences of the form   $(e^{i (n^2\gamma+n\delta)})$, meaning, it satisfies
 $$
 \lim_{N\to\infty} \frac{1}{N}\sum_{n=1}^N e^{i [n\alpha]n\beta}\cdot e^{i (n^2\gamma+n\delta)}=0$$
 for every $\gamma,\delta\in\R$.

\subsection{The Host-Kra factors  and their structure} \label{SS:HK} It is shown in \cite{HK05a} that the seminorms
$\nnorm{\cdot}_{k}$  induce $T$-invariant sub-$\sigma$-algebras
 $\cZ_{k-1}$   that satisfy
\begin{equation}
\label{E:DefZ_l}
\text{\em for } \ f\in L^\infty(\mu),\ \ \E(f|\cZ_{k-1})=0\quad \text{ \em if and
  only if } \quad   \nnorm f_{k} = 0.
\end{equation}
As a consequence, if for some $k_1,\ldots, k_\ell \in \N$  we have
an estimate of the form
\begin{equation}\label{E:SemiEst}
\limsup_{N\to\infty} \norm{A_N(f_1,\ldots, f_\ell)}_{L^2(\mu)} \leq C
\min_{i=1,\ldots,\ell}{\nnorm{f_i}_{k_i,T_i}}
\end{equation}
for some constant $C=C_{f_1,\ldots, f_\ell}$, then  the factors $\cZ_{k_1-1,T_1},\ldots, \cZ_{k_\ell-1,T_\ell}$ are characteristic for mean convergence of the averages \eqref{E:GeneralMulties'}.  Under such circumstances, one gets  characteristic factors with the sought-after algebraic structure. This is a consequence of an
important result of   B.~Host and B.~Kra \cite{HK05a}  which states that for ergodic systems
the factor system $(X, \cZ_{k}, \mu, T)$ is an inverse  limit of  $k$-step
 nilsystems. Note that for ergodic systems the factor $\cZ_1$ coincides with the Kronecker factor.

 Depending on the problem, it may be more
useful to reinterpret  the  structure theorem of Host and Kra
 as a decomposition result:
For every ergodic system $(X,\X,\mu,T)$ and function $f\in L^\infty(\mu)$,  for every $k\in \N$ and  $\varepsilon>0$, there exist
measurable functions $f_s$ (the structured component), $f_u$ (the uniform component), $f_e$ (the error term),
with $L^\infty(\mu)$ norm at most $2\norm{f}_{L^\infty(\mu)}$, such that
\begin{itemize}
\item $f=f_s+f_u+f_e$;
\item $\nnorm{f_u}_{k+1}=0$; \  $\norm{f_e}_{L^1(\mu)}\leq \varepsilon$; \ and
\item  $(f_s (T^nx))_{n\in\N}$
is a basic $k$-step nilsequence for $\mu$-almost every $x\in X$.
\end{itemize}
Such a decomposition also holds for
  non-ergodic  systems  (see Proposition~$3.1$ in  \cite{CFH11}).

Combining the hypothetical seminorm estimates \eqref{E:SemiEst} with
the aforementioned  structure theorem (or the decomposition result),
the problem of analyzing the limiting behavior of the averages
\eqref{E:GeneralMulties'} is
 reduced to a new   problem which typically amounts to proving  certain
 equidistribution properties of sequences on nilmanifolds.
 A lot of tools for handling such equidistribution problems have been developed in recent years, thus making such a reduction very much worthwhile.
 Some
examples  of useful equidistribution results on nilmanifolds  can be found in
\cite{AGH63,Fr08,  Fr09, FrK05a, GT09c,GT14, GT09d, HK05a, Lei05a, Lei05b, Lei06, Lei07, Lei09, Les89, Les91, Zi05}.

\subsection{A general strategy} When one is aiming to prove a
multiple recurrence or mean convergence property in ergodic theory by analyzing the limiting behavior of the
multiple ergodic averages  \eqref{E:GeneralMulties'}, very often one goes through  the following three
steps:
 \begin{itemize}
\item  Produce seminorm estimates like those in \eqref{E:SemiEst}.

\item  Use a structure theorem or a decomposition result  to reduce matters to
the analysis of the averages  \eqref{E:GeneralMulties'} for  nilsystems.

\item Use qualitative or quantitative equidistribution results on nilmanifolds to complete the analysis.
\end{itemize}
This approach is also implicit in the
foundational paper of Furstenberg \cite{Fu77} and subsequent works \cite{BL96,  FuK79, FuKO82}.
An important difference
is that in these works, in place of nilsystems, one uses the much larger class of
 distal systems, but for  several recent applications (for example for the precise evaluation of limits) the class of distal systems is  too
complicated to be able to deal with  directly.

The reader can find several examples demonstrating  how this general strategy is used  in order to prove multiple recurrence and convergence results,
as well as related applications in combinatorics, in the following
articles:
     \cite{As10,  BHK05,   BLL08, BLM11,  Chu09,     CF11, CFH11, CL88a,  Fr04, Fr08,  Fr10, Fr15a,  FrJLW10,   FrK05a,  FrK05b, FrK06, FrLW09,   FrW09,  FuW96, Gr11,  HK04,    HK05a, HK05b, HK09, HK11a,
Joh11, Lei05c,  Les93a, Po11,  Ru95,  WZ11, Zi06,  Zi07}. Depending
on the problem, the difficulty of each step varies; typically the
first step is elementary and is carried out by successive uses of
the Cauchy-Schwarz inequality and an  estimate of van der
Corput
(or Hilbert space
variants of it), the second step involves the use of (a modification of) the structure theorem of B.~Host and B.~Kra,
and the third step is a combination of algebraic and analytic
techniques. In the next subsections we explain some techniques that help us execute the first and the third steps.


\subsection{The polynomial exhaustion technique}\label{SS:PET}
We briefly explain a technique that is often used in order to produce seminorm
estimates of the type \eqref{E:SemiEst}. It is based on an induction
scheme (often called PET induction) introduced by  V.~Bergelson in
\cite{Be87a}. Let $\mathcal{F}\mathrel{\mathop:}=\{a_1,\ldots, a_\ell\}$ be a
family of real valued sequences, and suppose that one wishes to
establish seminorm estimates of the form
\begin{equation}\label{E:SemiEsti}
\limsup_{N\to\infty} \norm{A_N(f_1,\ldots, f_\ell)}_{L^2(\mu)} \leq C
\min_{i=1,\ldots,\ell}{\nnorm{f_i}_{k_i}},
\end{equation}
for some constant $C=C_{f_1,\ldots, f_\ell}$, where
$$
A_N(f_1,\ldots,f_\ell)\mathrel{\mathop:}=\frac{1}{N}\sum_{n=1}^N
T^{[a_1(n)]}f_1\cdots T^{[a_\ell(n)]}f_\ell.
$$
Variants of this technique could also be used to get similar estimates for some multiple ergodic
averages involving commuting transformations, but this is usually a much more difficult task.

The main idea is to
use the following Hilbert space variant  of van der Corput's fundamental estimate
\begin{theorem*}
 Let $v_1,\ldots, v_N$ be elements of an inner product space
of norm less than $1$. Then
$$ \norm{\frac{1}{N}\sum_{n=1}^N
 v_n}^2\leq 4\Big( \frac{1}{R}\sum_{r=1}^R (1-rR^{-1}) \Re\big(\frac{1}{N}\sum_{ n=1}^N
\langle v_{n+r}, v_n \rangle \big) + R^{-1}+RN^{-1}\Big)
$$
  for every integer  $R$ between $1$ and
$N$ where $\Re(z)$ denotes the real part of $z$.
\end{theorem*}
Using this estimate, in various cases  we  can
  bound the left hand side in \eqref{E:SemiEsti} by
 an expression that involves families of sequences of smaller ``complexity".
The goal is  after a finite number of iterations to get families of sequences
 that are simple  enough to handle directly.
 The details depend on the family of  sequences at hand, but
  typically,  after the first iteration, we get an  upper bound
 by an  average over $r\in \N$ of the $L^2(\mu)$-norm of
   multiple ergodic averages
   with iterates taken from the family of sequences (upon taking integer parts)
  \begin{equation}\label{E:VDC}
  \mathcal{F}_{a,r}\mathrel{\mathop:}=\{\text{sequences of the form } a_i(n)-a(n) , \ a_i(n+r)-a(n), i=1,\ldots, \ell\}
  \end{equation}
 where   $a\in\mathcal{F}$ is  fixed (so in particular independent of $r\in \N$) and chosen appropriately.
 If a choice of $a\in\mathcal{F}$ can be  made  so that the family
 $\mathcal{F}_{a,r}$ has smaller ``complexity'' than $\mathcal{F}$
  except possibly for a finite number of  $r\in\N$,
then,   practically, this means that after applying
  van der Corput's estimate a finite number of times, we will   be able to bound the left hand side in
  \eqref{E:SemiEsti} by  a much simpler expression for which we can prove the desired seminorm estimates directly.

 This strategy has been employed successfully in several instances and produced seminorm estimates of the form \eqref{E:SemiEsti}
for  linear sequences \cite{HK05a}, polynomial sequences
\cite{HK05b, Kou15a, Lei05c},
    variable polynomials of fixed degree \cite{FrW09}, and some sequences arising from smooth functions of polynomial growth \cite{BH09, Fr10, Fr15a}.
A common desirable  feature that these sequences share is that   after taking
successive differences (meaning iterating the operation $a(n)\mapsto
a(n+r)-a(n)$) a finite number of times, we arrive to  sequences that
are either constant or  asymptotically constant.
  This feature is not shared by several other
sequences worth studying, for example, random sequences of integers,  the sequence of primes, and  the
sequences $([n^{\log n}])$ or $([n \sin n])$. In these cases different approaches are needed in order to produce
seminorm estimates for ergodic averages with iterates given by such sequences.

 \subsection{Equidistribution  of polynomial sequences on nilmanifolds}\label{SS:nilequi}
Let  $X=G/\Gamma$ be a nilmanifold. We say that a sequence
$g\colon \N\to X$ is equidistributed in  $X$
if  for every $F\in C(X)$ one has
$$
\lim_{N\to\infty}\frac{1}{N}\sum_{n=1}^N F(g(n))=\int F \ dm_X
$$
 where $m_X$ denotes the Haar measure on $X$. A similar definition applies if $X$ is a
 union of finitely many nilmanifolds.

Let  $b_1,\ldots, b_\ell\in G$,
$x\in X$, and $a_1,\ldots, a_\ell\colon \N\to \Z$ be
sequences. In several of the applications we have in mind, one is at
some point called to deal with equidistribution properties of  sequences of the form
$(g(n)x)$, defined by
$g(n)\mathrel{\mathop:}=b_1^{a_1(n)}\cdots
b_\ell^{a_\ell(n)}$, $n\in\N$.
Such examples cover as special cases
sequences of the form $\big((c_1^{a_1(n)}x_1, \ldots,
c_\ell^{a_\ell(n)}x_\ell)\big)$, defined on the product of the
nilmanifolds $X_1, \ldots, X_\ell$. To see this, let
$X\mathrel{\mathop:}=X_1\times\cdots \times X_\ell$,
$x\mathrel{\mathop:}=(x_1,\ldots, x_\ell)$,
$b_1\mathrel{\mathop:}=(c_1,e_2,\ldots, e_l)$, $\ldots$,
$b_\ell\mathrel{\mathop:}=(e_1,\ldots,e_{\ell-1},c_\ell)$ where
$e_i$ denotes the identity element of the group $G_i$ for $i=1,\ldots, \ell$.

 Such problems are typically much
easier to handle when $X=\T^d$, since in this case one can utilize  the following equidistribution result of Weyl:
If $Y$ is a sub-torus of $\T^d$, then a
sequence
$g\colon \N\to Y$   is equidistributed on $Y$  if and only if for
every  non-trivial  character $\chi \colon Y\to \C$ one has
$$
\lim_{N\to\infty} \frac{1}{N}\sum_{n=1}^N\chi(g(n))=0.
$$
This criterion allows us to
 verify equidistribution properties of arbitrary sequences on $Y$   by  estimating   certain exponential sums. Unfortunately,
 such a convenient reduction is not available for non-Abelian
nilmanifolds, and checking equidistribution in this broader setup
can be very challenging  even
for simple sequences.
The situation is much better understood
when  all the sequences $a_1,\ldots, a_\ell\colon \N\to \Z$ are given by integer polynomials;
in this case we call  $(g(n)x)$ \emph{a polynomial sequence on $X$}.
Next we  state  some key results used to prove equidistribution properties of such polynomial sequences on nilmanifolds.

 Let
 $X = G/\gG$ be a connected nilmanifold.  By  $G_0$ we
 denote the connected component of the identity element
 in $G$, and for technical
  reasons we  assume that $G_0$ is simply connected and that
  $G=G_0\Gamma$. When $X$ is connected, for the problems we are interested in,  we   can always arrange matters so that $G_0$ has these additional properties.
 We let  $Z\mathrel{\mathop:}=G/([G_0,G_0]\Gamma)$
  and   $\pi\colon X\to Z$ be the natural projection.
  It is very desirable to work with the nilmanifold   $Z$ instead of $X$ because $Z$ has much
  simpler structure. Indeed, if $G$ is connected, then $Z$ is a connected
  compact Abelian Lie group, hence, a torus (meaning $\T^d$ for some $d\in \N$), and
  as a consequence every nilrotation in $Z$ is (isomorphic to) a  rotation
  on some torus. In general, the
  nilmanifold $Z$ may be more complicated, but it is  the case that
  every nilrotation in $Z$ is (isomorphic to) a  {\em unipotent  affine transformation}
  on some torus (see Proposition~3.1 in
\cite{FrK05a}), meaning, it has the form  $T(x) = S(x)+b$,
where $S\colon \T^d\to \T^d$  is a  homomorphism
which is unipotent (that is,  there exists $k\in\N$ so that $(S-{\text Id})^{k}=0$) and $b\in\T^d$.
 Iterates of such transformations  can be computed
explicitly
and so one is much more comfortable to be dealing with
equidistribution problems that involve unipotent affine
transformations on some torus than with general nilsystems.

\subsubsection{Qualitative equidistribution}
 The following qualitative equidistribution results were established by A.~Leibman  in
 \cite{Lei05a}:
\begin{itemize}
\item A polynomial
  sequence $(g(n)x)$ is always equidistributed in a finite union of
  sub-nilmanifolds of $X$.

\item  A polynomial
  sequence $(g(n)x)$ is equidistributed in $X$ if and
  only if the sequence  $(g(n)\pi(x))$ is equidistributed in $Z$.
\end{itemize}
The second statement gives   an efficient  way for checking
equidistribution of   polynomial sequences on nilmanifolds. We illustrate this with a simple example.
Suppose that
$b\in G$ is an ergodic nilrotation (meaning the transformation
$x \mapsto bx$ is ergodic) and we want to show that
the polynomial sequence  $(b^{n^2}x)$ is equidistributed in $X$ for every
$x\in X$. In the case where $G$ is connected,  the nilmanifold $Z$ is a torus,
therefore, according to the previous criterion, it suffices to show that
if $\beta$ is an ergodic element of $\T^d$ (this is the case
if the coordinates of $\beta$ are rationally independent), then for every $x\in X$
 the
sequence $(x+n^2\beta)$ is  equidistributed in $\T^d$. This is a
well known fact, and  can be easily verified using Weyl's
equidistribution theorem and van der Corput's  estimate.
If $G$ is not necessarily connected, one needs
to show that if $S\colon \T^d\to \T^d$ is an  ergodic unipotent
affine transformation, then the sequence $S^{n^2}x$ is
equidistributed for every $x\in \T^d$.
 Although this is somewhat
harder to establish, it follows again by    Weyl's
equidistribution theorem modulo some straightforward  computations.

\subsubsection{Quantitative equidistribution}
 Suppose that  one seeks to study  equidistribution properties of the  sequence $(b^{[n^{3/2}]}x)$, or   tries to prove uniform convergence of the sequence $(\frac{1}{N}\sum_{n=1}^NF(b^{n^2}x))$ to the integral of $F$, where $b\in G$ is an ergodic element and $F\in C(X)$.
In  such cases (and several others) one needs to use quantitative
 variants of the previous qualitative equidistribution results.
Such  results were proved by B.~Green and T.~Tao
\cite{GT09c, GT14}. In order to state them we have to introduce some notation.

For simplicity we assume that we work
on a nilmanifold $X=G/\Gamma$ with $G$ connected. As before, we let
$Z=G/([G,G]\Gamma) \cong \T^l$ and $\pi\colon X\to Z$ be the natural
projection map. A \emph{horizontal character} is a character
$\chi\colon Z\to \C$. We have  that $\chi(t)=e^{2\pi i \kappa\cdot t}$, $t\in \T$,
for some $\kappa\in\Z^l$, where $\cdot$ denotes the inner product
operation.
 We let $\norm{\chi}=|\kappa|$.
Suppose that $p\colon\Z\to \R$ has the form $
p(t)=\sum_{j=0}^d n^j \alpha_j $ where $\alpha_j\in\R$ for $j=0,\ldots, d$. We
define
\begin{equation}\label{E:norms}
  \norm{e^{2\pi i p(n)}}_{C^\infty[N]}=\max_{1\leq j\leq d}( N^j \norm{\alpha_j})
\end{equation}
where $\norm{x}=d(x,\Z)$.

Given $N\in\N$, a finite sequence $(g(n)\Gamma)_{1\leq n\leq N}$ is
said to be $\delta$-\emph{equidistributed on $X$}, if
$$
\Big|\frac{1}{N}\sum_{n=1}^N F(g(n)\Gamma)-\int_{X}F \ dm_X\Big|\leq
\delta \norm{F}_{\text{Lip}(X)}
$$
for every Lipschitz function $F\colon X\to \C$, where
$$
\norm{F}_{\text{Lip}(X)}=\norm{F}_\infty+ \sup_{x,y\in X, x\neq
  y}\frac{|F(x)-F(y)|}{d_X(x,y)}
$$
and $d_X$ is the Riemannian metric on $X$.

The quantitative equidistribution result in  \cite[Theorem~2.9]{GT09c}
states that for every nilmanifold $X$,  $d\in \N$, and
 $\delta>0$, there exist a real number $M=M_{X,d,\delta}$
  with the following property:
  For every
  $N\in\N$,
   if $g\colon \Z\to G$ is a polynomial sequence of degree at most $d$
  such that the finite sequence $(g(n)\Gamma)_{1\leq n\leq N}$ is not
  $\delta$-equidistributed,
   then for some non-trivial horizontal character
  $\chi$ with  $\norm{\chi}\leq M$  we have
  \begin{equation}\label{E:badequi}
   \norm{\chi( g(n))}_{C^\infty[N]}\leq  M,
 \end{equation}
   where $\chi$ is thought of as a
     character of the horizontal torus $Z=\T^l$ and $(g(n))$ as a
  polynomial sequence on $\T^l$.

 In the special case where  $X=\T$ (with the standard metric) and the polynomial
  sequence on $\T$ is given by $p(n)=n^d\alpha+q(n)\! \pmod{1}$ where $d\in \N$,
 $\alpha\in \R$, and $q\in\Z[x]$ with $\deg(q)\leq d-1$, the
  previous result implies that: There exists $M=M(d,\delta)>0$ such
  that for every $N\in\N$ and $\delta$ small enough,
  if the finite sequence $\big(n^d\alpha+q(n))_{1\leq n\leq N}$ is not
  $\delta$-equidistributed in $\T$, then $\norm{k\alpha}\leq
 M/N^d$ for some non-zero $k\in \Z$ with $|k|\leq M$.

\subsection{Exploiting randomness.}\label{SS:random}
In some cases we understand the limiting behavior of certain multiple ergodic averages
and we would like to exploit weighted variants of those averages, or averages along subsequences
of the iterates involved. If the weight (or the subsequence) enjoys some randomness features, it is often preferable to study the weighted averages indirectly, by comparing them with the unweighted ones and showing that the difference converges to $0$. To give an example, suppose that $R$ is the set of  positive integers that have an even number of prime factors counted with multiplicity; this set is known to have density $1/2$ and is expected to be distributed randomly.
We are interested in proving mean convergence for the averages
$$
\frac{1}{N}\sum_{n=1}^N w_n \, T^nf\cdot S^ng
$$
where $w_n:={\bf 1}_R(n)$, $n\in\N$. Since mean convergence of the unweighted averages is known, it suffices to show that the difference
\begin{equation}\label{E:diff}
\frac{1}{N}\sum_{n=1}^N (w_n-1/2) \, T^nf\cdot S^ng
\end{equation}
converges to $0$ in $L^2(\mu)$. This can be done using structural decomposition results
for multiple correlation sequences (see \cite{FH15c}), but the most economical way (given
existing knowledge) is to use the van der Corput Lemma twice in order to bound the $L^2(\mu)$ norm
of the averages \eqref{E:diff} by the $U^3(\Z_N)$ Gowers uniformity norm of the sequence $(w_n-1/2)$ (note that $w_n-1/2=\lambda(n)/2$ where $\lambda$ is the Liouville function). Since these norms are known to converge to $0$ as $N\to \infty$ (see \cite{GT08a,GT09d}), the averages \eqref{E:diff} converge to $0$ in $L^2(\mu)$, and we are done.
 In a similar way, we can treat weights of the form $({\bf 1}_R(n))$ where $R=\{n\in\N\colon \norm{n^k\alpha}\leq 1/2\}$  for $k\geq 3$ and $\alpha$  irrational (see \cite{FrLW06}),
 or $R=\{n\in\N\colon \phi(n)=1\}$ where $\phi$ is an aperiodic multiplicative function  with values in $\{\pm 1\}$ (see \cite{FH15c}). Moreover, we can cover  weights of the form  $(\Lambda(n))$ where $\Lambda$ is the von Mangoldt function (see \cite{FrHK07}), or  (unbounded) weights supported on random subsets of the integers
  of zero-density (see \cite{FrLW11}).

\subsection{The approaches of Tao and Walsh}
Mean convergence for averages of the form
$$
\frac{1}{N}\sum_{n=1}^N
T_1^{p_1(n)}f_1\cdots T_\ell^{p_\ell(n)}f_\ell,
$$
 was established
for linear polynomials $p_1,\ldots, p_\ell$  by  T.~Tao \cite{Ta08} and for general polynomials by
M.~Walsh \cite{Wal12}. Although the technical details in these two arguments are  different, the general strategy
is rather similar.
 Both approaches
  proceed by  reformulating the  mean convergence result  as  a finitary quantitative convergence
 problem for averages of products of general bounded sequences. After decomposing one of these sequences  into a sum of a structured component and a component that contributes negligibly in the averaging, one   is left with analyzing the contribution of the structured component. This component is defined using an averaging operation and its precise
structure is not further analyzed;  this contrasts previous ergodic
convergence arguments where the bulk of the proof involves a detailed
analysis of the structured component. What is important, is that
when the original sequence is replaced by its structured
component, the corresponding averages reduce to ``lower
complexity'' ones. The proof then concludes  with an induction on the
complexity of the averages involved (which is similar but not identical  to the PET induction of Section~\ref{SS:PET}), and eventually reduces matters
to the trivial case where one deals with  constant polynomials.

To give an example, suppose that one seeks to prove convergence for
averages of the form $\frac{1}{N}\sum_{n=1}^NT^nf\, S^{n^2}g$, where
the transformation $T$ and $S$ commute. After a few applications of
the previous method, one is reduced to studying averages of the form
$\frac{1}{N}\sum_{n=1}^NT_1^nf_1 \cdots T_r^nf_r$, for some $r\in \N$, where the
transformations $T_1,\ldots, T_r$ commute. After each subsequent
application of the previous method, the number of transformations is reduced by one, hence, after
$r$ additional steps one is reduced to a trivial convergence result.
Let us remark though, that   the
induction hypothesis demands certain uniformity conditions on the
various parameters involved, and  these parameters unfortunately
live strictly  on the finitary universe. This makes it impossible
to carry out this argument entirely on the infinite world
of ergodic theory. On the other hand, see \cite{Au13b} for  an ergodic proof of mean convergence
inspired by the argument of Walsh.

A key advantage of the methods of Tao and Walsh is that they use  ``weak''
decomposition results that do not require detailed
information on the structured component of the sequences involved. This is the main
reason for  the effectiveness and brevity of their arguments and lead to significantly
simpler proofs even for previously known mean convergence results that
relied on complicated  structure theorems. The price to
pay is that these  methods do not give  explicit information for the
 limit function  and it is also not clear if they can be  adjusted
  in order to prove mean convergence
 results  even in the simplest cases where  the
 iterates involved are not given by polynomial sequences.
For such problems, it
seems that in several cases  one still has to rely on  the more elaborate theory of characteristic factors
described in Sections~\ref{SS:2.2}-\ref{SS:nilequi}.

\subsection{Pleasant and magic extensions} Motivated by  work of  T.~Tao \cite{Ta08}, several people, including  H.~Towsner \cite{To09},  T.~Austin \cite{Au09},
and  B.~Host \cite{Ho09}, introduced new tools that help
us handle
 multiple ergodic averages with commuting transformations. In particular, a  key conceptual breakthrough that first appeared in \cite{Au09},
  is that in some instances by working with   suitable extensions of a family of  systems
  (called ``pleasant extensions'' in
 \cite{Au09} and ``magic extensions'' in \cite{Ho09}), characteristic factors of the corresponding multiple ergodic averages may be chosen to have particularly simple structure, a structure that is not visible when one works with the original system (the idea of passing to an extension in order to  simplify some
 convergence problems already appears implicitly in \cite{FuW96}). This is a rather counterintuitive statement
 since characteristic factors of extensions are extensions of characteristic  factors of the original systems. So in order to clarify matters, we explain
 a simple instance where such an approach works.

 Suppose that one wants to prove mean  convergence for the averages
 $$
A_N(T,S,f_i)\mathrel{\mathop:}= \frac{1}{N^2}\sum_{1\leq m,n\leq N}
T^mf_1\cdot S^nf_2 \cdot T^mS^n f_3
 $$
 where $(X,\X,\mu, T,S)$ is a system  and
  $f_1,f_2,f_3\in L^\infty(\mu)$.
 Although an  estimate that relates the $L^2(\mu)$-norm of these averages with the Gowers-Host-Kra seminorms
 of the individual functions  with respect to either $T$ or $S$ is not feasible, the following
 estimate is valid
 $$
\limsup_{N\to\infty} \norm{A_N(T,S,f_i)}_{L^2(\mu)}\leq C\,  \min_{i=1,2,3}{\nnorm{f_i}_{T,S,\mu}},
 $$
for some  constant $C=C_{f_1,f_2,f_3}$,  where
 $$
 \nnorm{f}_{T,S,\mu}^4\mathrel{\mathop:}=\lim_{M\to\infty}\frac{1}{M}\sum_{m=1}^M\lim_{N\to\infty}\frac{1}{N}\sum_{n=1}^N\int f\cdot
 T^m\bar{f} \cdot S^n\bar{f} \cdot T^mS^n f \ d\mu
$$
(it is shown in \cite{Ho09}  that $\nnorm{f}_{T,S, \mu}=\nnorm{f}_{S,T,\mu}$). Now, although factors of the
original systems that control the seminorms  $\nnorm{\cdot }_{T,S,\mu}$ may not admit particularly neat structure, it is shown in \cite{Ho09} that
there exists a new system $(X^*,\mu^*,T^*,S^*)$ that extends the system  $(X,\mu,T,S)$  and
in addition enjoys the following key  property (the term ``magic extension'' from \cite{Ho09}  alludes to this property):
$$
\nnorm{f^*}_{T^*,S^*,\mu^*}=0  \Leftrightarrow f^* \bot\  \cI_{T^*}
\vee \cI_{S^*} $$ where $f^*\in L^\infty(\mu^*)$ and $\cI_{T^*}$ denotes
the $\sigma$-algebra of $T^*$-invariant sets and similarly for $\cI_{S^*}$.
In fact, one
can take $X^*\mathrel{\mathop:}=X^4$,
$T^*\mathrel{\mathop:}=(\text{id},T,\text{id},T)$,
$S^*\mathrel{\mathop:}=(\text{id},\text{id},S,S)$, and define the
measure $\mu^*$  by
$$
\int f_1\otimes f_2 \otimes f_3\otimes f_4 \
d\mu^* \mathrel{\mathop:}=
\lim_{M\to\infty}\frac{1}{M}\sum_{m=1}^M\lim_{N\to\infty}\frac{1}{N}\sum_{n=1}^N\int
f_1\cdot
 T^mf_2 \cdot S^n f_3 \cdot T^mS^n f_4 \ d\mu.
 $$
Note  that mean convergence for the averages $A_N(T,S,f_i)$
follows if we prove mean convergence for the averages $A_N(T^*,S^*,f^*_i)$. Combining all these
observations, we can easily reduce matters to proving mean convergence
for the averages $A_N(T^*,S^*,f^*_i)$ when all functions $f^*_i$ are  $\cI_{T^*} \vee \cI_{S^*}$-measurable. This is a significant simplification of our original problem, and in fact it is now straightforward to  deduce the required convergence property from the mean ergodic theorem. Indeed, this is immediate if $f_3^*$ is a product of a $T^*$-invariant and a $S^*$-invariant function, and the general case follows by linearity and approximation.

This approach has proved  particularly useful for handling mean convergence problems of  multiple ergodic averages of commuting transformations  with linear iterates  (and in some cases for handling multiple recurrence problems) that  previously seemed intractable
 \cite{Au09,   Au10a,  Chu10,   Chu11,Ho09} (see also \cite{Au11c, Au14} for more general group actions).  A drawback to this approach is that it does not give  information
 about the precise form of the limiting function, and  also, up to now, it has not proved to be very  useful
 when some of the iterates are non-linear (for polynomial iterates though there is some progress in this
 direction \cite{Au11a,  Au11b}).

\subsection{Using structural results for  multiple correlation sequences}
Suppose we
want to give sufficient conditions on a bounded sequence of complex numbers $(w_n)$ so that  the weighted averages
\begin{equation}\label{E:weighted}
\frac{1}{N}\sum_{n=1}^N w_n\, T_1^nf_1 \cdots T_\ell^{n}f_\ell
\end{equation}
converge for all choices of systems and functions.
In the absence of a useful description for the characteristic factors of these averages, an alternate way to proceed is to use structural
results for the corresponding multiple correlation sequences.
Using the fact that the multiple correlation sequence
$$
\int f_0\cdot T_1^nf_1 \cdot \ldots \cdot T_\ell^{n}f_\ell\, d\mu, \quad n\in\N,
$$
is equal to a basic $\ell$-step nilsequence modulo a sequence that is small in uniform density
(see the theorem in Section~\ref{SS:526}), we deduce that the averages \eqref{E:weighted} converge weakly for all choices of systems and functions
if and only if the averages
$$
\frac{1}{N}\sum_{n=1}^N w_n\, \psi(n)
$$
converge for all basic $\ell$-step nilsequences $(\psi(n))$ (this condition is also necessary). Similar criterions also apply for weighted averages of multiple correlation sequences with  polynomial iterates (see \cite{BHK05, Chu11} for powers of the same transformation and  \cite{Fr15, FH15c} for general commuting transformations).

Somewhat stronger conditions
imply mean convergence; it suffices  to assume that either the averages
$$
\frac{1}{N-M}\sum_{n=M}^{N-1} w_n\, \psi(n)
$$ converge as $N-M\to \infty$ for every basic $\ell$-step nilsequence $\psi$ (see \cite{FH15c}), or that the averages
$$
\frac{1}{MN}\sum_{m=1}^M\sum_{n=1}^N w_n\, \overline{w_m}\, \psi(m,n)
$$
converge as $M,N\to\infty$ for every basic $(2\ell-1)$-step nilsequence $(\psi(m,n))$ (see \cite{FH15c}).

Although this approach is rather effective for bounded weights (see \cite{Fr15, FH15c, HK09} for applications involving weights of dynamical and arithmetic origin),  it is inadequate
for unbounded weights $(w_n)$, or if one is interested in  ergodic theorems along zero density subsequences of the iterates involved. For such purposes, and in order to fully exploit the potential of this approach, it seems that one needs to prove stronger structural decomposition results like those conjectured in Section~\ref{S:41}.

\subsection{Equivalent problems for sequences}
It is sometimes useful to be aware of the fact   that  mean convergence and multiple
recurrence problems in ergodic theory are intimately related to similar
problems involving bounded sequences of complex numbers.
 We give some explicit examples below.

Given a collection of sequences of integers $\{a_1,\ldots, a_\ell\}$,
it turns out that the following two properties are equivalent:
 \begin{itemize}
\item
 For every  system $(X,\X,\mu,T)$ and  set $A\in \X$ with $\mu(A)>0$,
 there exists $n\in \N$, such that
$$
\mu(A\cap  T^{a_1(n)}A\cap \dots \cap T^{a_\ell(n)}A)>0.
$$

\item For every  bounded sequence $(z(n))$  of non-negative  real numbers  that  satisfies   $\limsup_{M\to\infty}\frac{1}{M}\sum_{m=1}^M z(m)>0$, there exists $n\in \N$, such that
    $$\limsup_{M\to\infty}\frac{1}{M}\sum_{m=1}^M z(m)\cdot z(m+a_1(n))\cdot \ldots \cdot z(m+a_\ell(n))>0.$$
\end{itemize}
Using
the correspondence principle of Furstenberg it is not hard to see
that the first statement implies the second. To see that the second
statement implies the first it suffices to set
$z(m)\mathrel{\mathop:}={\bf 1}_A(T^mx)$, $m\in \N$,  for a suitable point $x\in X$
($\mu$-almost every $x\in X$ works) and use the pointwise ergodic theorem.

 We say that a
sequence of complex numbers $(z(n))$   \emph{ admits correlations along the
sequence of intervals $([1,M_k])$}, where $M_k\to \infty$,  if
 for every $\ell\in \N$ and $n_1,\ldots, n_\ell\in \Z$, the averages
$$\frac{1}{M_k}\sum_{m=1}^{M_k}z_1(m+n_1)\cdots z_\ell(m+n_\ell)$$
converge as $k\to\infty$ for all  sequences $z_1,\ldots, z_\ell\in \{z,\overline{z}\}$.
 Using a diagonal argument it is easy to show  that if $M_k\to\infty$, then any bounded sequence of complex numbers
 $(z(n))$ admits correlations  along some subsequence  of $([1,M_k])$.
Given a collection of sequences of integers $\{a_1,\ldots, a_\ell\}$,
it turns out that the following two properties are equivalent:
 \begin{itemize}
\item
 For every  system $(X,\X,\mu,T)$ and  $f\in L^\infty(\mu)$,
 the averages
$$
\frac{1}{N}\sum_{n=1}^N \int f \cdot T^{a_1(n)}f\cdot \ldots \cdot T^{a_\ell(n)}f \ d\mu
$$
converge  as $N\to\infty$.

\item  For every  bounded sequence $(z(n))$ of complex numbers, which admits correlations
along the sequence of intervals $([1,M_k])$,
the averages
    $$\frac{1}{N}\sum_{n=1}^N\Big(\lim_{k\to\infty}\frac{1}{M_k}\sum_{m=1}^{M_k} z(m)\cdot z(m+a_1(n))\cdot \ldots \cdot  z(m+a_\ell(n))\Big)$$
    converge as $N\to\infty$.
\end{itemize}
 The proof of the previous equivalence  is based on a variation of the correspondence principle of Furstenberg that applies to bounded sequences and can be found, for example, in \cite[Proposition~6.4]{FH15c}.

One can get similar statements for mean convergence, as well as
 for convergence and recurrence properties involving several commuting transformations
(in this case one has to use sequences in several variables).

\section{Some useful notions}
To ease our exposition we collect  some notions that are frequently
used in subsequent sections.
\subsection{Recurrence}
The next  notions are used to describe multiple recurrence properties of a collection of possibly different sequences:
\begin{definition}
The collection of sequences of integers $\{(a_1(n)),\ldots, (a_\ell(n))\}$
   is
   \begin{itemize}
\item      {\em   good for  $\ell$-recurrence of a single transformation} if  for every
system
    $(X,\X,\mu,T)$ and set $A\in \X$  with $\mu(A)>0$, we have
  \begin{equation*}
    \mu(A\cap T^{a_1(n)}A \cap \dots \cap T^{a_\ell(n)}A )>0
  \end{equation*}
  for some  $n\in \N$  with  $a_i(n)\neq 0$ for  $i=1,\ldots,\ell$.

 \item   {\em   good for  $\ell$-recurrence of commuting transformations} if
for every system $(X,\cB, \mu,$ $T_1,\ldots, T_\ell)$ and set $A\in \X$  with $\mu(A)>0$, we have
  \begin{equation*}
    \mu(A\cap T_1^{a_1(n)}A \cap\dots \cap T_\ell^{a_\ell(n)}A )>0
  \end{equation*}
   for some  $n\in \N$  with  $a_i(n)\neq 0$ for  $i=1,\ldots,\ell$.
\end{itemize}
\end{definition}
We remark that in the previous statements the existence of a single
$n\in\N$ for which the multiple intersection has positive measure,
forces the existence of infinitely many $n\in\N$ with the same
property.

 Examples of collections of  sequences
 that are known to be good for  $\ell$-recurrence of commuting transformations
 include collections of: integer polynomials with zero constant term \cite{BL96}
 and integer polynomials with zero constant term
evaluated at the shifted primes  \cite{FrHK11}.
 For
 arbitrary collections of integer polynomials necessary and sufficient conditions for $\ell$-recurrence of a single transformation are  given in \cite{BLL08}.
 For every $\ell\in \N$, the collection $\{([n^{c_1}]), \ldots, ([n^{c_\ell}])\}$, where $c_1,\ldots,c_\ell\in
 \R\setminus\Z$,
is known to be good for  $\ell$-recurrence of a single
transformation \cite{Fr10}, and for $\ell$-recurrence of commuting transformations if the non-integers $c_1,\ldots, c_\ell$ are distinct \cite{Fr15a}.

Next we define notions related to  multiple recurrence properties of collections defined by  a single sequence:
\begin{definition}
  The sequence of integers $(a(n))$
   is
   \begin{itemize}
\item      {\em good for  $\ell$-recurrence of powers} if
for all non-zero $k_1,\ldots, k_\ell\in\Z$  the collection of sequences $\{(k_1a(n)),\ldots, (k_\ell a(n))\}$
is good for  $\ell$-recurrence of a single transformation.

 \item  {\em   good for  multiple recurrence of powers} if it is
      good for  $\ell$-recurrence of powers  for every $\ell\in \N$.

 \item   {\em   good for  $\ell$-recurrence of commuting transformations} if
the collection of sequences $\{(a(n)),\ldots, (a(n))\}$
is good for  $\ell$-recurrence of commuting transformations.

\item   {\em   good for  multiple recurrence of commuting transformations} if it is
 good for  $\ell$-recurrence of commuting transformations for every $\ell\in \N$.
\end{itemize}
\end{definition}

The fact that the sequence $(n)$ is good for  multiple recurrence of powers
corresponds to the multiple recurrence result of H.~Furstenberg~\cite{Fu77}, and the fact that it
is  good for  multiple recurrence of commuting transformations corresponds to the multidimensional extension of this result
 of H.~Furstenberg and Y.~Katznelson~\cite{FuK79}. Further examples of sequences that are good for multiple recurrence of commuting
 transformations include: integer polynomials with zero constant term \cite{BL96}, the shifted primes (see \cite{WZ11} for powers  and
  \cite{BLZ11} in general), integer polynomials with zero constant term evaluated at the shifted primes (see \cite{WZ11} for powers and \cite{FrHK11}
   in general),  arbitrary shifts of integers with an even (or odd) number of prime factors \cite{FH15b}, several generalized polynomial sequences \cite{BeMc10}, and some random sequences of integers of zero density \cite{FrLW11}.
   The sequence $([n^c])$, where $c\in \R$ is positive,
is known to be good for  multiple recurrence of powers \cite{FrW09} (see also \cite{Fr10}), but
if $c\in \R\setminus \Q$ is greater than $1$, it is not known whether it is
good for multiple recurrence of commuting transformations.
Examples of sequences that are good for $\ell$-recurrence of powers but are not good for $(\ell+1)$-recurrence
of powers can be found in \cite{FrLW06, FrLW14, HKM15}.

It is also useful to know some obstructions to recurrence.
One can check that if the sequence $(a(n))$ is good for
$1$-recurrence, then the equation $a(n)\equiv 0 \mod r$  has a
solution for every $r\in\N$; as a consequence, the sequences
$(p_n)$, where $p_n$ is the $n$-th prime,  $(n^2+1)$, $(2^n)$, are not good for $1$-recurrence. More
generally, if the sequence $(a(n))$ is good for $1$-recurrence, then
for every $\alpha\in \R$ and $\varepsilon>0$ the inequality
$\norm{a(n)\alpha}\leq \varepsilon$  has a solution (letting
$\alpha=1/r$ gives the previous congruence condition); as a
consequence, the sequence $([\sqrt{5}n+2])$ is not good for
$1$-recurrence (take $\alpha=1/\sqrt{5}$ and $\varepsilon=1/10$).
This obstruction also implies that
 if a sequence $(a(n))$ is lacunary, meaning, $\liminf_{n\to\infty}{a(n+1)/a(n)}>1$,
 then it is not good for $1$-recurrence, since it is well known that in this case
 there exist an  irrational $\alpha$ and a positive $\varepsilon_0$ such that
 $\norm{a(n)\alpha}\geq \varepsilon_0$ for every $n\in \N$.

 Furthermore, if a sequence $(a(n))$ is good for $2$-recurrence of powers, then for every
$\alpha\in \R$ and $\varepsilon>0$ the inequality
$\norm{a(n)^2\alpha}\leq \varepsilon$  has a solution (this is
implicit in \cite{Fu81a} Section 9.1, and is proved in detail in
\cite{FrLW06, HKM15}). It follows that if one forms a sequence by putting
the elements of the set $\{n\in \N\colon \norm{n^2\sqrt{2}}\in
[1/4,1/2]\}$ in increasing order, then this sequence is not going to
be good for $2$-recurrence of powers (one can show that it is going
to be good for $1$-recurrence \cite{FrLW06, HKM15}). And there are also additional
more restrictive obstructions, if the sequence $(a(n))$ is good for
$2$-recurrence of powers, then for every $\alpha,\beta \in \R$ and
$\varepsilon>0$ the inequality $\norm{[a(n)\alpha]a(n)\beta}\leq
\varepsilon$ has a solution \cite{FrLW09}. More generally, for
$\ell$-recurrence of powers, one can state further obstructions by using
 higher degree polynomials
with zero constant term, or other generalized
polynomials  \cite{HuSY13}.

The correspondence principle of Furstenberg enables us to deduce statements
in density Ramsey theory
from  multiple recurrence statements in ergodic  theory.  Using this principle and some elementary
arguments, one is able to reformulate the previous ergodic notions in
purely combinatorial terms:


\begin{theorem*}
The collection of sequences of integers $\{a_1(n),\ldots, a_\ell(n)\}$
   is
\begin{itemize}
\item
  good for  $\ell$-recurrence of a single transformation
if and only if
 every set $E\subset \Z$ with
  $\bar{d}(E)>0$ contains patterns  of the form
  $$\{m, m +a_1(n),\ldots, m+ a_\ell(n)\}$$ for some $m\in \Z$ and $n \in\N$ with
  $a_i(n)\neq 0$ for $i=1,\ldots, \ell$.

\item    good for  $\ell$-recurrence of commuting transformations
if and only if  for every ${\bf v}_1,\ldots, {\bf v}_\ell\in
\Z^d$,  every set  $E\subset \Z^d$ with
  $\bar{d}(E)>0$ contains patterns of the form
  $$\{{\bf m}, {\bf m} +a_1(n){\bf v}_1,\ldots, {\bf m}+a_\ell(n){\bf v}_\ell\}$$ for some ${\bf m}\in \Z^d$  and  $n \in\N$ with   $a_i(n)\neq 0$
  for $i=1,\ldots, \ell$.
  \end{itemize}
\end{theorem*}
For instance, the fact that the sequence $(n)$ is good for multiple recurrence of powers  corresponds to the theorem of Szemer\'edi
on arithmetic progressions, and the fact that for polynomials $p_1,\ldots, p_\ell\in \Z[t]$ with zero constant term the collection $\{p_1,\ldots, p_\ell\}$ is good for multiple
recurrence of powers corresponds to the polynomial Szemer\'edi theorem.

Let us also remark that the previous notions admit equivalent
uniform versions that are often  useful for applications. For
instance, one can prove the following (see \cite{BHRF00} for an
argument that works for polynomials and \cite{FrLW09} for an argument
that works for  general sequences):
\begin{theorem*}
Let $(a_1(n)), \ldots, (a_\ell(n))$ be sequences of integers.
Then the following statements are equivalent:

\begin{itemize}

\item  The  collection  $\{a_1(n),\ldots,
a_\ell(n)\}$ is good for  $\ell$-recurrence of a single
transformation.

\item For every $\varepsilon>0$ there exist
$\delta>0$ and $N_0$, such that
for every $N\geq N_0$ and   integer set $E\subset[-N,N]$ with
$|E|\geq\varepsilon N$, we have
$$
|E\cap (E-a_1(n))\cap \dots \cap(E-a_\ell(n))|\geq \delta N
$$
for some $n\in [1,N_0]$.

\item For every $\varepsilon>0$ there exist
$\gamma>0$ and $N_1$, such that
for every   system
$(X,\X,\mu,T)$ and $A\in \X$ with $\mu(A)
\geq\varepsilon$, we have that
$$
\mu(A\cap T^{a_1(n)}A\cap \dots\cap T^{a_\ell(n)}A)\geq\gamma
$$
for some $n\in [1,N_1]$.
\end{itemize}
\end{theorem*}
We remark that the constants $\gamma, \delta,$ and $N_0$ in the previous theorem
depend only on $\varepsilon$ and the choice of the sequences $a_1,\ldots, a_\ell$.

\subsection{Convergence}
The next notion is used to describe multiple convergence  properties
of a collection of
 sequences:
\begin{definition}
The collection of sequences of integers $\{(a_1(n)),\ldots, (a_\ell(n))\}$
   is
   \begin{itemize}
\item      \emph{good for $\ell$-convergence of a single transformation} if
for every   system $(X,\mathcal{X},\mu,T)$ and
functions $f_1,\ldots, f_\ell \in L^\infty(\mu)$, the averages
$$
\frac{1}{N}\sum_{n=1}^N T^{a_1(n)}f_1 \cdots  T^{ a_\ell(n)}f_\ell
$$
converge in the mean.

 \item   \emph{good for $\ell$-convergence of commuting transformations} if
for every system $(X,\mathcal{X},\mu,$ $T_1,\ldots, T_\ell)$ and  functions $f_1,\ldots, f_\ell \in L^\infty(\mu)$, the averages
$$
\frac{1}{N}\sum_{n=1}^N T_1^{a_1(n)}f_1 \cdots
T_\ell^{a_\ell(n)}f_\ell
$$
converge in the mean.
\end{itemize}
\end{definition}
Examples  of  collections of sequences (not coming from  multiples of
the same sequence) that are known to be good for $\ell$-convergence
of commuting transformations include those defined by:
integer polynomials \cite{Wal12}, the integer part of polynomials with
real coefficients \cite{Kou15a},   any such sequence evaluated at the prime numbers \cite{FrHK11, Kou15b},
and  pairs of sequences of the form $\{(n), (a_n(\omega))\}$
where $(a_n(\omega))$ is a certain random sequence of integers of zero
density \cite{FrLW11}. Additional examples of collections of sequences that
 are good for $\ell$-convergence of a single transformation include the collection $\{([n^{c_1}]), \ldots, ([n^{c_\ell}])\}$, where $c_1,\ldots,c_\ell\in \R\setminus \Z$ are positive and several other Hardy field sequences of polynomial growth \cite{Fr10},
 these sequences are also known to be  good  for $\ell$-convergence of commuting transformations if the non-integers $c_1,\ldots, c_\ell$ are distinct \cite{Fr15a}.

Next we define notions related  to  multiple convergence properties of collections defined by  a single sequence:
\begin{definition}
The sequence of integers $(a(n))$
   is
   \begin{itemize}
\item      \emph{good for $\ell$-convergence of powers} if
for every $k_1,\ldots, k_\ell\in\Z$  the collection of sequences $\{(k_1a(n)),\ldots, (k_\ell a(n))\}$
is good for  $\ell$-convergence of a single transformation.

\item     \emph{good for multiple convergence of powers}
if it is good for $\ell$-convergence of powers  for every $\ell \in \N$.

 \item   \emph{good for $\ell$-convergence of commuting transformations}
if
the collection of sequences $\{(a(n)),\ldots, (a(n))\}$
is good for  $\ell$-convergence of commuting transformations.

\item      \emph{good for multiple convergence of commuting transformations}
if it is good for $\ell$-convergence of commuting transformations  for every $\ell\in \N$.
\end{itemize}
\end{definition}
Examples  of  sequences that are good for multiple convergence of
commuting transformations include sequences given by
integer polynomials \cite{Wal12}, the integer part of polynomials with
real coefficients \cite{Kou15a},   any such sequence evaluated at the prime numbers \cite{FrHK11, Kou15b},
and some random sequences of integers of  zero density \cite{FrLW11, FrLW14}. Additional examples
of sequences that are known to be good for multiple convergence of
powers include
sequences of the form $([n^a])$ where
$a>0$ and several other  Hardy field sequences of polynomial growth \cite{Fr10}. Examples of sequences that are good for
$\ell$-convergence of powers but are not good for
$(\ell+1)$-convergence of powers  can be found in \cite{FrLW06}.


\section{Problems related to general sequences}
In this section we give a list of problems related to the study of
  multiple ergodic averages and multiple recurrence problems involving iterates given by  general sequences of integers.

\subsection{The structure of multiple correlation sequences}\label{S:41}
Let $(X,\cX,\mu, T_1,\ldots, T_\ell)$ be a system  and  $f_0,f_1,\ldots,f_\ell\in L^\infty(\mu)$
be functions.    We are interested in determining the structure of the  multiple correlation sequences $\mathcal{C}\colon \Z^\ell\to \C$ defined by the formula
\begin{equation}\label{E:MultCor}
\mathcal{C}(n_1,\ldots,n_\ell)\mathrel{\mathop:}=\int f_0\cdot
T_1^{n_1}f_1\cdot \ldots \cdot  T_\ell^{n_\ell}f_\ell \ d\mu, \quad n_1,\ldots, n_\ell\in\Z.
\end{equation}
The next result gives a very satisfactory solution to this problem for $\ell=1$ and
 serves as our model for possible generalizations.
It can be deduced  from Herglotz's theorem on positive definite
sequences (the sequence $n\mapsto \int \bar{f}\cdot T^nf\ d\mu$ is
positive definite) and a polarization identity.
\begin{theorem*}
Let $(X,\mathcal{X},\mu, T)$ be a system and
$f,g\in L^\infty(\mu)$. Then  there exists a complex Borel measure
$\sigma$ on $[0,1)$, with bounded variation, such that for every
$n\in \Z$ we have
\begin{equation}\label{E:SpectralSingle}
\int f \cdot T^ng \ d\mu =\int_0^1 e^{2\pi i nt} \ d\sigma(t).
\end{equation}
\end{theorem*}

Finding a formula analogous to \eqref{E:SpectralSingle},  with the
multiple correlation sequences \eqref{E:MultCor} in place of the
single correlation sequences,  is  a problem of fundamental importance which has been in the
mind of experts for several years.
 A satisfactory solution  is going to give us new insights and significantly improve our ability to deal
with difficult questions involving multiple ergodic averages which at the moment seem out of reach.
 There are indications that  sequences of polynomial nature should  replace the
 ``linear'' sequences $(e^{2\pi i nt})$.
The most reasonable candidates at this point seem to be some
collection of generalized multivariable nilsequences (defined in Section~\ref{SS:nil}).
 For instance, examples of generalized $2$-step nilsequences in $1$-variable are
the sequences
 $\big(e^{i( [n\alpha]n\beta +n\gamma)}\big)$ ($[x]$ denotes the integer part of $x$),
where $\alpha,\beta,\gamma\in \R$, and examples of  generalized  $2$-step nilsequences in
 $2$ variables are  the sequences   $\big(e^{i( [m\alpha]n\beta +m\gamma +n\delta)}\big)$,
where $\alpha,\beta,\gamma,\delta\in \R$.

 Let  $\mathcal{S}$ be  a subset of $\ell^\infty(\N^\ell)$. We say that
a bounded sequence $a\colon \Z^\ell\to \C$ is
\begin{itemize}
\item {\em an integral combination of elements of  $\mathcal{S}$}, if
there exist a complex Borel measure
$\sigma$  of bounded variation on a compact metric space $X$  and
sequences $a_x\in \mathcal{S}$, $x\in X$,  such that for each ${\bf n} \in \Z^\ell$ the map  $x\mapsto a_x({\bf n}) $ is  in
$L^\infty(\sigma)$ and for every ${\bf n}\in \Z^\ell$ we have
$$
a({\bf n})=\int_X a_x({\bf n})\ d\sigma(x).
$$

\item  {\em an approximate integral combination of elements of  $\mathcal{S}$}, if for every $\varepsilon>0$ we have $a=b+c$ where $b$ is an integral combination of elements of $\mathcal{S}$ and $\norm{c}_\infty\leq \varepsilon$.
\end{itemize}
\begin{problem}\label{Conj:StructMultiGeneral}
Determine the structure of the multiple correlation sequences
$(\mathcal{C}(n_1,\ldots,n_\ell))$ defined by \eqref{E:MultCor}.
Is it true that any such sequence is an (approximate) integral combination
of generalized $\ell$-step nilsequences in $\ell$-variables?
\end{problem}
We remark that a similar result may hold even if the transformations $T_1,\ldots, T_\ell$ generate a nilpotent group. But some commutativity assumption on the transformations is needed,
otherwise simple examples show that generalized nilsequences cannot be the only
building blocks. For instance, let $T,S\colon \T\to\T$ be  given by
$Tx\mathrel{\mathop:}=2x$, $Sx\mathrel{\mathop:}=2x+\alpha$, and
$f(x)\mathrel{\mathop:}=e^{- i x}, g(x)\mathrel{\mathop:}=e^{ i x}$.
Then $\int T^nf\cdot S^ng\ dx= e^{i\cdot  (2^n -1)\alpha}$ and one
can show that $(e^{i\cdot 2^n\alpha})$ is not a generalized nilsequence for
$\alpha\in \R\setminus\Q$.



 Even resolving special cases of this problem would be extremely interesting.
 One particular instance is the following:
 \begin{problem1*}
 Let $(X,\cX,\mu,T)$ be a system and $f,g,h\in L^\infty(\mu)$. Is it true that
the sequence  $(\mathcal{C}(n))$ defined by
$$
\mathcal{C}(n)\mathrel{\mathop:}=\int f\cdot T^ng\cdot T^{2n} h\
d\mu, \quad n\in \N,
$$
 is an (approximate) integral combination
of generalized $2$-step nilsequences?
\end{problem1*}
 In \cite{BHK05} it is shown that for ergodic systems one has the decomposition   $$\mathcal{C}(n)=\psi(n)+e(n)$$ where $\psi(n)$ is a (single variable) $2$-step nilsequence and $\lim_{N\to\infty}\frac{1}{N}\sum_{n=1}^{N}|e(n)|=0$. A variant of this result regarding correlation sequences defined by commuting transformations was established in \cite{Fr15a} and an extension regarding correlation sequences in several variables in \cite{FH15c}. Unfortunately, these results do not   provide information on the error term $e(n)$ for zero density subsets of the integers, and as a consequence they are of little use when one studies sparse subsequences of the sequence $\mathcal{C}(n)$.
Lastly,
we remark that an essentially equivalent problem is to characterize the structure of sequences $(\mathcal{C}(n))$ defined by
$$
\mathcal{C}(n)\mathrel{\mathop:}=\lim_{M\to\infty} \frac{1}{M}\sum_{m=1}^M a(m)\cdot a(m+n)\cdot a(m+2n),\quad  n\in\N,
 $$
  where $a\in \ell^\infty(\Z)$ is any bounded sequence that admits correlations along the sequence of intervals $([1,M])$.

It is  natural to inquire whether correlation sequences defined using commuting transformations provide new examples of sequences that cannot be constructed using correlation
sequences that use only powers of the same transformation.    To  formulate precise statements let us define as $\mathcal{C}_{T,S}$ to be the set of all sequences of the form
$$
\Big(\int f\cdot T^ng\cdot S^nh\, d\mu\Big)
$$
where  $(X,\X,\mu,T,S)$ ranges over all systems and  $f,g,h$ over all functions in $ L^\infty(\mu)$. Let also
$\mathcal{C}_{T}$ be the set of sequences of the previous form where we impose the restriction that $T$ and $S$ are powers of the same transformation. It was shown in \cite{Fr15a}, that modulo sequences that are small in uniform density,
the set
$\mathcal{C}_{T,S}$ coincides with the set of  basic $2$-step nilsequences, and also that modulo terms that are small in $\norm{\cdot}_\infty$, every basic $2$-step nilsequence belongs to the set $\mathcal{C}_{T}$. Hence, letting
$
\norm{a}_2:=\limsup_{N-M\to \infty} \frac{1}{N-M}\sum_{n=M}^{N-1}|a(n)|^2$ for $a\in \ell^\infty(\N)$,
one deduces the following:
\begin{theorem*}
Let $a\in \mathcal{C}_{T,S}$. Then for every $\varepsilon>0$  there exists $b\in \mathcal{C}_{T}$ such that $\norm{a-b}_2\leq \varepsilon$.
\end{theorem*}
In fact, we believe that the sets  $\mathcal{C}_{T,S}$ and  $\mathcal{C}_{T}$ coincide:
\begin{problem}\label{C:TS}
Show that $\mathcal{C}_{T,S}=\mathcal{C}_{T}$.
\end{problem}
It can be shown that  a positive answer to Problem~\ref{Conj:StructMultiGeneral} will give  a positive answer
to Problem~\ref{C:TS}.

\subsection{Necessary and sufficient conditions for $\ell$-convergence}


The next result can be deduced from  formula \eqref{E:SpectralSingle} and serves as our model
for giving usable  necessary and sufficient conditions for $\ell$-convergence:
\begin{theorem*}
If $(a(n))$ is a sequence of integers,  then the  following
  statements are equivalent:
  \begin{itemize}
  \item The sequence $(a(n))$ is  good for $1$-convergence.

  \item The sequence $(a(n))$ is  good for $1$-convergence
 for rotations on the circle.

  \item  The sequence $\big(\frac{1}{N}\sum_{n=1}^N
    e^{i a(n)t}\big)$ converges for every $t\in\R$.
  \end{itemize}
\end{theorem*}
Since
 a formula that generalizes \eqref{E:SpectralSingle} to multiple correlation sequences
 is not available, we are unable to prove analogous  necessary and sufficient conditions for $\ell$-convergence.
 Nevertheless, inspired by  Problem~\ref{Conj:StructMultiGeneral} we  make  the following natural guess:
\begin{problem}\label{Conj:CondConvGeneral}
  If $(a_1(n)),\ldots, (a_\ell(n))$ are  sequences of integers, then  show that
    the following 
  statements are equivalent:
  \begin{itemize}
  \item The sequences  $(a_1(n)),\ldots, (a_\ell(n))$ are  good for $\ell$-convergence of commuting transformations.

 \item The sequences  $(a_1(n)),\ldots, (a_\ell(n))$ are good for $\ell$-convergence
 of  $\ell$-step nilsystems.

  \item The sequence $\big( \frac{1}{N}\sum_{n=1}^N
    \psi(a_1(n),\ldots,a_\ell(n))
    \big)$ converges for every  basic generalized $\ell$-step nilsequence  $\psi$ in $\ell$-variables.
  \end{itemize}
\end{problem}
A similar problem was   formulated in \cite{FrJLW10}.
We are mainly interested in knowing if the third (or the second) condition implies the first.
Even the following very special case is open:
\begin{problem2*}
Let $(a(n))$ be a sequence of integers such that the  averages $\frac{1}{N}\sum_{n=1}^N
    \psi(a(n))$ converge for every basic  generalized $2$-step nilsequence  $\psi$.
Show  that for every system $(X,\mathcal{X},\mu,
T, S)$ and functions $f,g\in L^\infty(\mu)$, the averages
$$\frac{1}{N}\sum_{n=1}^N  T^{a(n)}f\cdot S^{a(n)}g $$
converge weakly in $L^2(\mu)$.
\end{problem2*}
  The problem has been solved when  $(a(n))$ is strictly increasing and its range has positive  density;
this follows from
  \cite[Theorem~1.9]{BHK05} for $S=T^2$ and in the general case from
 \cite{Fr15}.
For general sequences the problem is open even when $S=T^2$.

A similar problem that is probably of  equivalent difficulty
 is the following one related to
 weighted averages:
\begin{problem2**}
Let $(w_n)$ be a sequence of complex numbers such that the  averages $\frac{1}{N}\sum_{n=1}^N
    w_n\, \psi(n)$ converge for every  generalized $2$-step nilsequence  $\psi$.
Show  that for every system $(X,\mathcal{X},\mu,
T, S)$ and functions $f,g\in L^\infty(\mu)$, the averages
$$\frac{1}{N}\sum_{n=1}^N w_n\,  T^nf\cdot S^ng $$
converge weakly in $L^2(\mu)$.
\end{problem2**}
This problem is solved  for bounded sequences $(w_n)$
 (see \cite{Fr15a}) and is open for general sequences
 even when $S=T^2$.


\subsection{Sufficient conditions  for $\ell$-recurrence}



Our model result is the next theorem of T.~Kamae and
M.~Mend\`es-France \cite{KaMe78} that gives usable
conditions  for checking that a sequence is good for
$1$-recurrence:
\begin{theorem*}
Let $(a(n))$ be   sequence of
integers that satisfies:
\begin{itemize}
  \item the sequence  $(a(n) \alpha)_{n\in \N}$ is equidistributed
  in $\mathbb{T}$ for every irrational $\alpha$, and

\item the set $\{n\in \N \colon r|a(n)\}$ has positive upper
density for every $r\in \N$.
\end{itemize}

\noindent Then the sequence $(a(n))$ is   good for $1$-recurrence.
\end{theorem*}
Once again,  the proof of this result relies upon the   identity
\eqref{E:SpectralSingle}. Since appropriate generalizations of this identity are not
known  for multiple correlation sequences, we are unable to give a
similar criterion for $\ell$-recurrence when $\ell\geq 2$. To state
a conjectural criterion, we
  extend the notion
   of an irrational
   rotation on the circle  to general connected nilmanifolds:  Given a connected nilmanifold $X=G/\Gamma$,
an {\it irrational nilrotation in $X$} is an
   element $b\in G$ such that the sequence $(b^n\Gamma)_{n\in\N}$
   is equidistributed on $X$.
\begin{problem}\label{Conj:CondRec}
Let   $(a(n))$  be  a sequence that satisfies:
\begin{itemize}
  \item for every connected $\ell$-step  nilmanifold $X$ and
  every irrational nilrotation $b$ in $X$ the sequence
 $(b^{a(n)}\Gamma)_{n \in\N}$ is equidistributed  in $X$, and

\item the set $\{n\in\N \colon r| a(n)\}$ has positive
upper density for every $r\in \N$.
\end{itemize}

\noindent
Show that the sequence $(a(n))$
is  good for $\ell$-recurrence of commuting transformations.
\end{problem}
This  problem was first formulated in \cite{FrLW09} and in the same
article it was solved  for $\ell$-recurrence of powers for sequences that have range a set
of integers with positive density. The stated conditions are
satisfied by any non-constant integer polynomial
sequence with zero constant term (follows from  \cite{Lei05a}),
the sequence  $([n^c])$ for every $c>0$ (follows from   \cite{Fr09}),  the sequence
  of   shifted primes $(p_n-1)$ (follows from  \cite{GT09d}), and random non-lacunary sequences
  of integers (follows from results in \cite{Fr09}).

\subsection{Powers of sequences and recurrence}
It is known that if a sequence is  good for $\ell$-convergence of powers, then
 its first $\ell$  powers are  good for $1$-convergence. More precisely, the following holds (this is implicit in \cite{Fu81a} Section 9.1,
and is proved in detail in  \cite{FrLW06}):

\begin{theorem*}
If $(a(n))$ is  good for $\ell$-convergence of powers, then $(a(n)^k)$ is  good for $1$-convergence for $k=1,\ldots, \ell$.
\end{theorem*}

It is unclear whether a similar property holds for recurrence.
\begin{problem}\label{Conj:PowersRec}
If $(a(n))$ is  good for $\ell$-recurrence of powers, is then $(a(n)^k)$   good for $1$-recurrence for $k=1,\ldots, \ell$?
\end{problem}
This problem was first stated in \cite{FrLW06}
and it  is open even when $\ell=2$.
It is known that if  $(a(n))$ is  good for $2$-recurrence of powers, then the sequence
 $(a(n)^2)$   is good
for Bohr recurrence, meaning it is good for $1$-recurrence for all
rotations on tori (see \cite{Fu81a} Section 9.1, or \cite{FrLW06}).
A well known question of Y.~Katznelson asks whether a set of Bohr
recurrence is necessarily a set of topological recurrence (for
background on this question see \cite{BoGl09, Katz01,We00, HKM15}). Although there exist examples of sets of topological recurrence that are not sets of
 $1$-recurrence \cite{Kri87}, all known examples  are rather complicated. As a consequence,
 an example showing that the answer to  Problem~\ref{Conj:PowersRec} is negative is probably going to be complicated.
\subsection{Commuting vs powers of a single transformation}
If a sequence is  good for $2$-convergence of commuting transformations, then, of course,  it is also
 good for   $2$-convergence of powers.
 Interestingly,
no example that distinguishes the two notions is known and we believe that  there is none:
\begin{problem}\label{Conj:ConvPowersComm}
If     a sequence  is  good for $2$-convergence of powers, then   show that it is   good for $2$-convergence of commuting transformations.
\end{problem}
 The corresponding question for recurrence is also open:
\begin{problem}\label{Conj:RecPowersComm}
Is there    a sequence that is  good for $2$-recurrence of powers but  is not good for $2$-recurrence of commuting transformations?
\end{problem}
This question was first stated  in \cite{Be96} (Question 8) where  V.~Bergelson states that the answer is very likely yes.

\subsection{Fast growing sequences}
Despite the  successes in dealing with multiple recurrence and convergence problems of sequences that
do not grow faster than polynomials, when it comes down to  fast growing sequences our
knowledge is very limited. Say that a sequence $(a(n))$ of positive integers is {\em fast growing} if
$\lim_{n\to\infty} \log(a(n))/\log n=\infty$
(equivalently, if it is  of the form $(n^{b(n)})$ with $b(n)\to\infty$).
\begin{problem}\label{Conj:ConvRecFast}
 Give an explicit example of a fast growing sequence that is  good for  multiple recurrence  and convergence of powers and commuting transformations.
\end{problem}
 Even for $2$-recurrence and convergence of powers, no such example is known. It is known that if a  sequence grows exponentially fast (for example $(2^n)$ or $(n!)$), then it is not good for $1$-recurrence and $1$-convergence (even for irrational circle rotations). On the other hand, several  natural examples of fast growing sequences that do not grow
  exponentially fast   should
   be good for multiple recurrence and convergence,
for instance, the sequences $([n^{(\log n)^a}])$, $(n^{[(\log n)^b]})$, $([e^{n^c}])$  for $a,b>0$ and $c\in (0,1)$.
Unfortunately, it is very hard to work with these sequences, even for issues related to $1$-recurrence and $1$-convergence.
 Probably  a sequence like  $(n^{[\log\log{n}]})$ is  easier to work with. One could also try to construct  (not so explicit) examples using random sequences of integers
(more on this in   Section~\ref{S:random}).

\section{Problems related to polynomial sequences}
In this section we give a list of problems related to the study of
  multiple ergodic averages  involving iterates given by polynomial sequences, and related applications to multiple recurrence.

\subsection{Powers of a single   transformation}\label{SS:3.1}
Let  $\P\mathrel{\mathop:}=\{p_1,\ldots,p_\ell\}$ be a family of
integer polynomials that are \emph{essentially distinct}, meaning,
all polynomials and their differences are non-constant.  First, we
consider averages of the form
\begin{equation}\label{E:PolynomialMulties}
\frac{1}{N}\sum_{n=1}^N T^{p_1(n)}f_1\cdots T^{p_\ell(n)}f_\ell,
 \end{equation}
 where  $(X,\X,\mu,T)$ is a system  and $f_1,\ldots, f_\ell\in L^\infty(\mu)$.

We remark that all the mean convergence results stated in this section work
 equally well for uniform Cesaro averages, meaning averages of the form $\frac{1}{|\Phi_N|}\sum_{n\in \Phi_N}$ in place of the averages $\frac{1}{N}\sum_{n=1}^N$  where
 $(\Phi_N)_{N\in\N}$ is any  F{\o}lner sequence of subsets of $\N$, meaning, $|(\Phi_N+h)\triangle \Phi_N|/|\Phi_N|\to 0$ as $N\to\infty$ for every $h\in\N$.
This property is no longer  true if the iterates  are given by Hardy field sequences or random sequences of integers (see Sections~\ref{S:Hardy} and \ref{S:random}), or if we deal with  pointwise convergence (not even when $\ell=1$ and $p_1(n)=n$).

\subsubsection{Optimal characteristic factors}
Before discussing some problems related to  the  characteristic factors of the averages
\eqref{E:PolynomialMulties}, we state a
result of  B.~Host and B.~Kra \cite{HK05b} and A.~Leibman \cite{Lei05c} that gives useful information about their structure.
\begin{theorem*}
There exists $d=d(\P)$  such that  the factor $\cZ_{d,T}$ (defined in \ref{SS:HK})
is characteristic for mean convergence of the averages \eqref{E:PolynomialMulties}.
\end{theorem*}
We note that in the previous statement the value of $d(\P)$
 can be chosen to depend only on the number  and the maximum degree of the polynomials in $\P$. Given a
family of polynomials $\P$, we denote by  $d_{min}(\P)$ the minimal
value of $d(\P)$  that works in the previous theorem. This value
is in general hard to pin down and depends on the algebraic
relations that the polynomials satisfy. For instance, we know that
$d_{min}(\{n,2n,\ldots,\ell n\}) =\ell-1$ (\cite{HK05a, Zi07}),
and $d_{min}(\P) =1$ when  $\P$ consists of (at least two) rationally
independent polynomials  (\cite{FrK05a, FrK06}). But it is not only linear relations between the polynomials that matter, for
instance, we know that $d_{min}(\{n,2n,n^2\})=2$ while
$d_{min}(\{n,2n,n^3\})=1$ (\cite{Fr08,Lei09}). More
examples of families $\P$ where $d_{min}(\P)$ has been computed can
be found in
 \cite{Fr08,Lei09}. Furthermore,  in \cite{Lei09} a (rather complicated) algorithm is given for computing this value. Despite such progress, the following is still open (the problem is implicit in \cite{BLL07} and was  stated explicitly in \cite{Lei09}):
\begin{problem}\label{Conj:ComplexityBound}
  If $|\P|\geq 2$, show that
$
d_{min}(\P)\leq |\P|-1.
$
\end{problem}
The estimate   is known when  $|\P|=2, 3$ (\cite{Fr08}) and it
is open when $|\P|=4$. The problem is open even when one is
restricted to the class of Weyl systems,  meaning, systems of the
form $(\T^d, \X_{\T^d}, m_{\T^d}, T)$ where $T\colon
\T^d\to \T^d$ is a unipotent affine transformation.
We denote with
$d^{\text{W}}_{min}(\P)$ the minimum value of $d$ such that the
factor $\cZ_{d,T}$ is characteristic for mean convergence  of
the averages \eqref{E:PolynomialMulties} for all Weyl systems
(properties of $d^{\text{W}}_{min}(\P)$ were studied in
\cite{BLL07}).
 \begin{problem10*}  If $|\P|\geq 2$, show that
$ d^{\text{W}}_{min}(\P)\leq |\P|-1. $
\end{problem10*}
 This problem was first stated in \cite{BLL07} (set $W(P)\mathrel{\mathop:}=d^{W}_{min}(\P)+1$ in the remark after \cite[Proposition~5.3]{BLL07}).
  The estimate is known when  $|\P|=2,3,4$ (\cite{Fr08, McCl11}) and it  is open when $|\P|=5$.
  We also remark that no example is known where
 $d_{min}(\P)\neq d^{\text{W}}_{min}(\P)$, so it is natural to suspect that these two values coincide.
 This is known to be the case
when $|\P|=3$ (\cite{Fr08}) and it is open when $|\P|=4$. Obviously one
has $d^{\text{W}}_{min}(\P)\leq d_{min}(\P)$. Some bounds in the
other direction are given in \cite{Lei09}.

\subsubsection{Variable polynomials}
Mean convergence of the averages \eqref{E:PolynomialMulties} was established after
a long series of intermediate results; the papers    \cite{CL84, CL88a, CL88b,Fu77,
FuW96, HK01, HK05a,  Ru95, Zi07}  dealt with the
important case of linear  polynomials, and  using the machinery
 introduced in \cite{HK05a}, convergence for arbitrary polynomials was finally obtained by B.~Host and B.~Kra in~\cite{HK05b}
except for a few cases that were treated by A.~Leibman in~\cite{Lei05c}.
\begin{theorem*}
Let $(X,\mathcal{X},\mu,T)$ be a system,
$f_1,\ldots, f_\ell\in L^\infty(\mu)$ be functions,  and
$p_1,\ldots, p_\ell$ be  integer polynomials. Then the averages
\eqref{E:PolynomialMulties} converge in the mean as $N\to\infty$.
\end{theorem*}
Furthermore, explicit formulas for the limit can be given for
special families of polynomials \cite{Fr08, FrK05a, FrK06,
Lei09,Zi05}, but no explicit formula is known for general families of
polynomials.

We record here a related open problem involving variable polynomials.
We say that
\begin{itemize}
\item the sequence of  polynomials $(p_N)$ where $p_N\in \R[t]$, $N\in\N$, is {\em good} if the polynomials have bounded degree
 and for every non-zero $\alpha\in\R$ we have
$$
\lim_{N\to\infty}\frac{1}{N}\sum_{n=1}^N\, e^{ i p_N(n)\alpha}=0.
$$
\item the sequence of $\ell$-tuples of polynomials $(p_{1,N},\ldots, p_{\ell,N})$ where $p_{i,N}\in \R[t]$, $N\in\N$, $i\in \{1,\ldots, \ell\}$  is {\em good} if
     every non-trivial linear combination of the  sequences of polynomials $(p_{1,N}),\ldots, (p_{\ell,N})$ is good.
\end{itemize}
It can be shown  (see for example \cite[Lemma~4.4]{GT09c}) that if $p_N(t)=\sum_{k=1}^d c_{k,N}t^k$, $c_{k,N}\in \R$, $N\in \N$, then the sequence $(p_N)$ is good if and only if  for every non-zero $\alpha\in \R$ we have $\lim_{N\to\infty} N^j\norm{c_{j,N}\alpha}=\infty$ for at least one $j\in \{1,\ldots, d\}$ where $\norm{\cdot}$ denotes the distance from the closest integer.
Hence, for $\ell=2$, letting $p_{1,N}:=n/N^a$ and $p_{2,N}:=n/N^b$, $n, N\in\N$,  where $0<a<b<1$, produces an example of a sequence of  good pairs of polynomials of degree $1$.
Another example is defined by the polynomials $p_{1,N}:=n/N^a$,  $p_{2,N}:=n^2/N^a$, $\ldots$,
 $p_{\ell,N}:=n^\ell/N^a$, $n, N\in\N$, where $a\in (0,1)$. The next problem is of interest
 even for these collections of variable polynomials.

\begin{problem}
Suppose that the sequence of $\ell$-tuples of polynomials   $(p_{1,N},\ldots, p_{\ell,N})$
is good. Show that for every ergodic system $(X,\X,\mu,T)$ and functions  $f_1,\ldots,f_\ell\in L^\infty(\mu)$ we have
$$
\lim_{N\to\infty}\frac{1}{N}\sum_{n=1}^N\, T^{[p_{1,N}(n)]}f_1\cdots
 T^{[p_{\ell,N}(n)]}f_\ell=\int f_1\, d\mu\, \cdots \int f_\ell\, d\mu
$$
where convergence takes place in $L^2(\mu)$.
\end{problem}
For $\ell=1$ the problem can be solved using the spectral theorem, that is, the identity in  \eqref{E:SpectralSingle}.

\subsubsection{Pointwise convergence}
In most cases, it  is still unknown whether mean convergence of multiple ergodic averages  can be strengthened to pointwise convergence.
We mention two particular cases that are open:
\begin{problem}\label{Conj:PointConvSingPolies}
Let $(X,\X,\mu,T)$ be a system and $f, g, h\in L^\infty(\mu)$ be functions. Show that the averages
$$
\frac{1}{N}\sum_{n=1}^N f(T^nx)\cdot g(T^{2n}x)\cdot h(T^{3n}x)
\quad \text{and the averages} \quad
\frac{1}{N}\sum_{n=1}^N f(T^nx)\cdot g(T^{n^2}x)
$$
converge pointwise almost everywhere.
\end{problem}
Pointwise convergence  of the averages
\eqref{E:PolynomialMulties} is  known when $\ell=1$ \cite{Bou88a} and is also known when
$\ell=2$ and both polynomials are linear \cite{Bou90} (see also
\cite{D07} for an alternative proof). In all other cases the problem
is open even for weak mixing systems. Partial results that deal with
special classes of transformations can be found in \cite{
As98, As05,  Bere85, Bere88, DerL96,HuSY14,  Lei05a, Les93a, LRR03}.

We record also some problems of arithmetic nature regarding pointwise convergence of double averages.
We say that a multiplicative function $\phi\colon \N\to \C$ (meaning,
 $\phi(mn)=\phi(m)\phi(n)$ for all $m,n\in\N$ with $(m,n)=1$)  has {\em convergent means} if
the averages $\frac{1}{N}\sum_{n=1}^N\phi(an+b)$ converge as $N\to \infty$ for every $a,b\in\N$.
It can be shown that every multiplicative function that takes  values in $[-1,1]$  has convergent means, but there are multiplicative functions with values on the complex unit disc that do not have a mean value (for example take $\phi(n)=n^{it}, n\in \N,$ for some non-zero $t\in \R$).
\begin{problem}
Let $(X,\X,\mu,T)$ be a system and $f,g\in L^\infty(\mu)$ be functions.
 If $\Lambda$ is the von Mangoldt function and $\phi$ is a multiplicative function that takes values on the complex unit disc and has convergent means, then  show that the averages
$$
\frac{1}{N}\sum_{n=1}^N \Lambda(n)\,  f(T^nx)\cdot g(T^{2n}x)
\quad \text{and the averages} \quad
\frac{1}{N}\sum_{n=1}^N \phi(n)\, f(T^nx)\cdot g(T^{2n}x)
$$
converge pointwise almost everywhere.
\end{problem}
Pointwise convergence of the first  averages would imply pointwise convergence of the averages
$ \frac{1}{N}\sum_{n=1}^N  f(T^{p_n}x)\cdot g(T^{2p_n}x)$
where $p_n$ denotes the $n$-th prime.
Pointwise convergence of the second  averages when $\phi$ is the Liouville function would imply pointwise convergence of the   averages $ \frac{1}{N}\sum_{n=1}^N  f(T^{s_n}x)\cdot g(T^{2s_n}x)$, where $s_n$ is the $n$-th natural number that has an even number of prime factors counted with multiplicity.
Of course, one also expects similar properties to hold for higher order variants of the previous statements.
Mean convergence for both averages is known (\cite{FrHK07} for the first,  \cite{FH15b} for the second). Also pointwise convergence is known when $g=1$ ((\cite{Wi88} for the first, \cite{FH15b} for the second). The second problem is open even when $\phi$ is the M\"obius or the Liouville function
(in which case the $L^2$-limit is known to be $0$ \cite{FH15b}).

\subsubsection{Subsequences of multiple correlation sequences}
Lastly, we record a problem about correlation sequences involving powers of a single transformations.
We start with  the following
 result of V.~Bergelson,  B.~Host, and B.~Kra \cite{BHK05}:
\begin{theorem*}
Let $(X,\X,\mu,T)$ be an ergodic system and $f_0,\ldots, f_\ell\in L^\infty(\mu)$.
 Then one has a decomposition
$$
\int f_0 \cdot T^n f_1 \cdot \ldots \cdot  T^{\ell n} f_\ell \ d\mu= \psi(n)+e(n)
$$
where $(\psi(n))$ is  an $\ell$-step nilsequence and $\lim_{N\to\infty}\frac{1}{N}\sum_{n=1}^N|e(n)|=0$.
\end{theorem*}
 This result was extended  to non-ergodic
 systems in \cite{L15, L15b} and to polynomial iterates in \cite{Lei11a}; in the latter case the level of nilpotency of $\psi$  depends on $\ell$ and the degree of the polynomials. A key
ingredient in the proofs is the fact that the
factor $\cZ_{\ell,T}$ is characteristic for convergence of the
averages $\frac{1}{N}\sum_{n=1}^N|\int f_0 \cdot T^n f_1 \cdot
\ldots \cdot T^{\ell n} f_\ell \ d\mu|.$  The next problem seeks to explore the extend to which the previous theorem continues to hold for  subsequences of multiple correlation sequences.
\begin{problem}
Let $(a(n))$ be the  sequence of integers  $(p_n)$, where $p_n$ is the $n$-th prime, or $([n^c])$ where $c>0$, or $(2^n)$.
Is it true that for every ergodic system $(X,\X,\mu,T)$  and all  functions $f_0,\ldots, f_\ell\in L^\infty(\mu)$, one has a decomposition
$$
\int f_0 \cdot T^{a(n)} f_1 \cdot \ldots \cdot  T^{\ell a(n)} f_\ell \ d\mu= \psi(a(n))+e(n),
$$
where $(\psi(n))$ is  an ($\ell$-step) nilsequence and $\lim_{N\to\infty}\frac{1}{N}\sum_{n=1}^N|e(n)|=0$?
\end{problem}
We believe that the answer is yes for the first two sequences and no for the third sequence.

One could also ask similar questions for correlation sequences of the form
$$
\mathcal{C}(n)=\int f\cdot T^{a(n)} g \cdot T^{b(n)}h \ d\mu,\quad n\in\N,
 $$
 where $(a(n)), (b(n))$ are particular sequences of  integers. In this case,
  the expected decomposition  has the form
 $$
 \mathcal{C}(n)=\psi(a(n),b(n))+e(n),\quad  n\in\N,
 $$
  where
$(\psi(m,n))$ is a ($2$-step) nilsequence in  two variables and $(e(n))$ is as before.

\subsection{Commuting transformations}
Throughout this section,  $(X,\X,\mu,T_1, \ldots, T_\ell)$ is a system and $f_1,\ldots,f_\ell$ are functions in  $L^\infty(\mu)$.

We start with problems  related to mean convergence of
multiple ergodic averages. After a long series of partial results that dealt with
commuting transformations
 \cite{Au09, Au11a, Au11b,CFH11, CL84,   FrK05a, FuK79, Ho09, Joh11, Les93a, Ta08, To09, Zh96},
M.~Walsh in \cite{Wal12} proved the following result:
\begin{theorem*}\label{T:Walsh}
Let $p_1,\ldots, p_\ell$ be integer valued  polynomials. Then the   averages
\begin{equation}\label{E:Multies}
\frac{1}{N}\sum_{n=1}^{N} T_1^{p_1(n)}f_1\cdots T_\ell^{p_\ell(n)}f_\ell
\end{equation}
converge in the mean as $N\to\infty$.
\end{theorem*}
A similar result also holds when we work with averages of the form
$$
\frac{1}{N}\sum_{n=1}^{N} \big(\prod_{i=1}^\ell T_{i}^{p_{i,1}(n)}\big)f_1\cdots \big(\prod_{i=1}^\ell T_{i}^{p_{i,\ell}(n)}\big)f_\ell
$$
and the transformations $T_1,\ldots, T_\ell$ generate a nilpotent group \cite{Wal12} and when one  averages over F{\o}lner sequences in $\Z$ and uses  polynomials in several variables \cite{Zo15a}.
Moreover, we have mean convergence when the iterates are given  by integer parts of real valued polynomials  \cite{Kou15a}.

\subsubsection{Generalized polynomials}
A \emph{generalized polynomial} is a real valued function that is obtained
from the identity function and real constants by using the
operations of addition, multiplication, and taking  integer parts. An example  is the sequence $([[n\alpha]n\beta+n^2\gamma+ n\delta])$ where $\alpha,\beta, \gamma, \delta\in \R$.
It remains a challenge to extend known mean convergence results with polynomial iterates  to the case where the iterates
are given by generalized polynomials:
\begin{problem}\label{Conj:CommPolConv}
Let $p_1,\ldots, p_\ell$ be integer valued generalized polynomials.  Show that  the  averages
$$
\frac{1}{N}\sum_{n=1}^{N} T_1^{p_1(n)}f_1\cdots T_\ell^{p_\ell(n)}f_\ell
$$
converge in the mean as $N\to\infty$.
\end{problem}
This problem was first stated in \cite{Be06b}.
For $\ell=1$ convergence was proved in \cite{BL07} and follows from the spectral theorem and the fact (proved in  \cite{BL07}) that any sequence of the form $(e^{ip(n)})$, where $p$ is a generalized polynomial, can be represented as a generalized nilsequence. For $\ell=2$ the problem is open even when the transformations are equal and weak mixing.

\subsubsection{Characteristic factors}
As mentioned previously, if all transformations are equal, and the
polynomials are essentially distinct,  then  characteristic factors
of the averages \eqref{E:Multies} can be chosen to have very special
algebraic structure. For general commuting transformations this is
no longer the case; if  one  chooses $p_1=p_2=n$, $T_1=T_2$, and
$f_2=\bar{f}_1$, then the averages \eqref{E:Multies} do not converge
to $0$ unless $f_1=f_2=0$. The same problem persists when two of
the polynomials are \emph{dependent}, meaning, their quotient is
some non-trivial linear combination of the polynomials is constant.
But in all other cases, there is no obvious obstruction
to
 having  ``simple"  characteristic factors with algebraic structure:
  \begin{problem}\label{Conj:CommCharPairIndep}
 Suppose that  the polynomials
  $p_1,\ldots, p_\ell\in \Z[t]$ are pairwise independent. Show
that there exists   $d\in \N$ such that the factors $\cZ_{d,T_1}, \ldots, \cZ_{d,T_\ell}$ are  characteristic  factors for the
averages \eqref{E:Multies}.
\end{problem}
 This is known to be the case when the polynomials have distinct degrees \cite{CFH11}.
  But it is not known  for some simple families of integer polynomials, for instance, for the family
     $\{n^3, n^3+n\}$ and the family  $\{n, n^2, n^2+n\}$. Even for weak mixing transformations the following problem is open:
  \begin{problem15*}
  Suppose that all
  the transformations $T_1, \ldots, T_\ell \colon X\to X$ are weak mixing
  and the polynomials
  $p_1,\ldots, p_\ell\in \Z[t]$ are pairwise independent.
Show that
$$
\lim_{N\to\infty}\frac{1}{N}\sum_{n=1}^{N} T_1^{p_1(n)}f_1\cdots T_\ell^{p_\ell(n)}f_\ell=
\int f_1 \ d\mu\cdots \int f_\ell \ d\mu.
$$
\end{problem15*}
The problem is open even when $\ell=2$ and the family of polynomials is  $\{n^3, n^3+n\}$.

 When all  transformations are equal and the polynomials are in
general position,  characteristic factors for the averages
\eqref{E:Multies} turn out to be  extremely simple; they are given by   the {\em rational Kronecker factor} of the system which is defined as follows: Given a  system
 $(X,\X,\mu,T)$ we let   $\mathcal{K}_{rat}(T)=\bigvee_{d\in \N} \mathcal{I}_{T^d}$ where $\mathcal{I}(T^d):=\{f\in L^2(\mu)\colon T^df=f\}$. The next result is proved in
\cite{FrK05a,
FrK06}:
\begin{theorem*}
  Suppose that the polynomials
  $p_1,\ldots, p_\ell\in \Z[t]$ are rationally independent (meaning, non-trivial linear combinations of them are non-constant).
Then the rational Kronecker factor $\mathcal{K}_{rat}(T)$ is a
characteristic factor for the averages \eqref{E:PolynomialMulties}.
\end{theorem*}
We believe  that this result generalizes to the case of several commuting transformations:
  \begin{problem} \label{Conj:CommCharInd}
  Suppose that the polynomials
  $p_1,\ldots, p_\ell\in \Z[t]$ are rationally independent.
Show that   $\mathcal{K}_{rat}(T_1),\ldots, \mathcal{K}_{rat}(T_\ell)$
are  characteristic factors for the averages \eqref{E:Multies}.
\end{problem}
This  was proved in \cite{CFH11} when  $\ell=2$ for families of the form $\{n,p(n)\}$ where $\text{deg}(p)\geq 2$. In
the same article a somewhat weaker property was proved for all
monomials with distinct degrees. Using techniques from \cite{CFH11} and \cite{Pan15}
it is very likely that the problem can be solved for all polynomial families  with distinct degrees.

\subsubsection{Optimal lower bounds for multiple recurrence} We mention also a closely related  multiple recurrence problem:
\begin{problem}\label{Conj:LowerBounds}
Suppose that the polynomials $p_1,\ldots,p_\ell\in \Z[t]$ are
rationally independent and have zero constant term. Show that for
every $A\in \X$ and every $\varepsilon>0$, there exists $n\in \N$
such that
\begin{equation}\label{E:lowerbounds}
\mu(A\cap T_1^{p_1(n)}A\cap \cdots\cap T_\ell^{p_\ell(n)}A)\geq \mu(A)^{\ell+1}-\varepsilon.
\end{equation}
\end{problem}
In fact, the set of integers  $n$ for which \eqref{E:lowerbounds}
holds is expected to have bounded gaps. The lower bounds are known  when all
transformations are equal \cite{FrK05b} and they are also known for general commuting transformations
when the polynomials are
monomials with distinct degrees \cite{CFH11}.
Combining techniques from \cite{CFH11} and  \cite{Pan15}
it is very likely that the problem can be solved for arbitrary polynomials with distinct degrees.   The result fails if two of
the polynomials are distinct and  dependent; in this case no
fixed power of $\mu(A)$ works as a lower bound in \eqref{E:lowerbounds} for arbitrary systems and sets \cite{BHK05}.  On the
other hand, the assumption that the polynomials are rationally
independent is not necessary, for instance, the result is expected to
work  for the family of polynomials $\{n,n^2,n^2+n\}$
 (this is known to be the case when all  transformations are equal \cite{Fr08}).
We remark that if Problem~\ref{Conj:CommCharInd} is solved,  then  the conjectured
lower bounds of Problem~\ref{Conj:LowerBounds} will follow rather
easily. On the other hand, Problem~\ref{Conj:LowerBounds} is open even when $\ell=2$.

\subsubsection{Multiple recurrence for intersective polynomials}
A family of integer valued polynomials $\{p_1,\ldots, p_\ell\}$ is called {\em intersective}
if for every $r\in \N$ there exists $n\in\N$ such that $p_i(n)\equiv 0\! \! \pmod{r}$
for $i=1,\ldots, \ell$.
When $\ell=1$ an example of an
intersective polynomial with no  linear factors  is $p(n)=(n^2-13)(n^2-17)(n^2-221)$. Examples of periodic systems show that if a family of polynomials is good for multiple recurrence of commuting transformations, then this
family has to be intersective.
The next problem was first stated in \cite{BLL08} and seeks to show that the condition of intersectivity is also sufficient for multiple recurrence of commuting transformations.
\begin{problem}\label{Conj:Intersective}
Let $(X,\X,\mu,T_1,\ldots, T_\ell)$ be a system and   $\{p_1,\ldots, p_\ell\}$
  be a family of intersective integer polynomials. Show that for every set $A\in \X$ with $\mu(A)>0$ one has
$$
\mu(A\cap T_1^{p_1(n)}A\cap \cdots\cap T_\ell^{p_\ell(n)}A)>0
$$
for some $n\in\N$.
\end{problem}
The problem is open even when $\ell=2$ and $p_1=p_2$.
The case where all the polynomials have zero constant term is covered by the Polynomial Szemer\'edi Theorem \cite{BL96}. The case where all the transformations are equal was
dealt in \cite{BLL08}. The argument used in  \cite{BL96} depends crucially on the fact that the polynomials have zero constant term. The argument
used in \cite{BLL08} depends crucially on the fact that characteristic factors for the corresponding multiple ergodic averages are (inverse limits) of nilsystems;  a
substitute for  this result that is useful for the problem at hand is not currently available.

Note that a solution to Problem~\ref{Conj:Intersective}  would imply that for every collection of intersective polynomials $\{p_1,\ldots, p_\ell\}$,  for every $d,\ell\in\N$, set $E\subset \N^d$ with $\bar{d}(E)>0$, and vectors  ${\bf v}_1,\ldots,  {\bf v}_\ell\in \N^d$,   there exist ${\bf m}\in \N^d$ and $n\in \N$ such that
$${\bf m}, {\bf m}+p_1(n) {\bf v}_1,\ldots, {\bf m}+p_\ell(n){\bf v}_\ell\in E.
$$

\subsubsection{Pointwise convergence for commuting transformations}
Regarding pointwise convergence of
 multiple ergodic averages of commuting transformations, progress has been scarce,
 even in some seemingly simple cases.
 The following is a well known open  problem:
\begin{problem}\label{Conj:PointwiseCommuting}
Let $(X,\X,\mu, T,S)$ be a system
 and $f,g\in L^\infty(\mu)$ be functions. Show that the averages
$$
\frac1N \sum_{n=1}^N f(T^nx)\cdot  g(S^nx)
$$
converge pointwise almost everywhere.
\end{problem}
For a list of partial results that apply to special classes of transformations see the list after Problem~\ref{Conj:PointConvSingPolies}. Moreover, pointwise convergence in the
case where the system $(X,\X,\mu, T,S)$ is distal was recently established in \cite{DoS15}, and prior to this, in the case where $S=T^k$ for some $k\in \N$ in \cite{HuSY14}.
See also \cite{DoS16} for more general statements involving an arbitrary number of commuting transformations.

\subsubsection{Optimal error term for decomposition results}\label{SS:526}
We end this subsection with a problem related to the structure of multiple correlation sequences defined by commuting transformations. The following result was proved in
\cite{Fr15}:
\begin{theorem*}
Let $(X,\X,\mu,T_1,\ldots, T_\ell)$ be a system. Then for every $\varepsilon>0$ one has a decomposition
$$
\int f_0 \cdot T_1^n f_1  \cdot \ldots \cdot  T_\ell^n f_\ell \ d\mu= \psi(n)+e(n)
$$
where $(\psi(n))$ is  a basic $\ell$-step nilsequence and $\limsup_{N\to\infty}\frac{1}{N}\sum_{n=1}^N|e(n)|\leq \varepsilon$.
\end{theorem*}
A similar result also holds when the iterates are given by integer polynomials \cite{Fr15}, or integer parts of real valued polynomials \cite{Kou15a},
and when one uses correlation sequences in several variables  \cite{FH15c}.
When $T_1,\ldots, T_\ell$ are powers of the same transformation, the main results in
\cite{BHK05, L15, L15b} give that one can have an analogous  decomposition  with $\varepsilon=0$
and $\psi$ an $\ell$-step nilsequence.
It is thus natural to ask whether a similar result holds for commuting transformations.
\begin{problem}\label{Conj:CommDecNil}
Is it true that one always has a decomposition
$$
 \int f_0 \cdot T_1^n f_1 \cdot \ldots \cdot   T_\ell^n f_\ell \ d\mu= \psi(n)+e(n)
$$
where $(\psi(n))$ is  an $\ell$-step nilsequence and $\lim_{N\to\infty}\frac{1}{N}\sum_{n=1}^N|e(n)|=0$?
\end{problem}
 Of course one could ask a similar question when the iterates are given by
polynomial sequences.

\subsection{Transformations that do not commute}
All problems in  the previous sections were stated for families of   transformations
  that commute.
  When one works with arbitrary families of invertible
measure preserving transformations    the next result shows that   one
cannot expect to have similar convergence and recurrence results:
\begin{theorem*}
Let $a,b\colon \N\to\Z\setminus\{ 0\}$ be $1-1$ sequences. Then
there exist  invertible Bernoulli measure preserving transformations
$T$ and $S$ acting on the same probability space $(X,\X,\mu)$ such
that
\begin{itemize}
  \item for some   $f,g\in L^\infty(\mu)$
the averages $\frac{1}{N}\sum_{n=1}^N\int T^{a(n)} f\cdot S^{b(n)}g \ d\mu$ diverge;

  \item for some $A\in \mathcal{X}$ with $\mu(A)>0$ we have $T^{a(n)}A\cap S^{b(n)}A=\emptyset$ for every $n\in \N$.
\end{itemize}
\end{theorem*}
To construct such examples it suffices to modify examples of
D.~Berend (see Ex $7.1$ in \cite{Bere85}) and H.~Furstenberg (page
$40$ in \cite{Fu81a}) that cover the case $a(n)=b(n)=n$ (the
details  appear in \cite{FrLW11}). When $a(n)=b(n)$,  it is also
known that given any finitely generated solvable group $G$ of
exponential growth, there exist invertible measure preserving
transformations $T,S$, with $\langle T,S\rangle \subset G$, and such that  the
conclusion of the previous theorem holds for those $T$ and $S$.
 It is interesting that despite such negative news,  once one introduces an extra variable,
 several convergence
(and very likely recurrence) results can be   extended to arbitrary
families of measure preserving transformations. We mention a positive result
from \cite{CF11}:
\begin{theorem*}
Let $(X,\X,\mu)$ be a probability space, $T_1,\ldots, T_\ell \colon
X\to X$ be invertible measure preserving transformations,
$f_1,\ldots,f_\ell \in L^\infty(\mu)$ be functions,
$p_1,\ldots,p_\ell$ be essentially distinct  polynomials of degree at most  $k\in \N$, and $a\in
(0,1/k)$. Then  the averages
\begin{equation}\label{E:Polynomial}
\frac{1}{N^{1+a}}\sum_{1\leq m\leq N, 1\leq n\leq N^a} f_1(T_1^{m+p_1(n)}x)\cdots f_\ell(T_\ell^{m+p_\ell(n)}x)
\end{equation}
converge pointwise almost everywhere as $N\to \infty$.\footnote{For $k=1$ one can actually take $a=1$, this was established in \cite{HuSY14} but also follows from the method of proof in \cite{CF11}.}
\end{theorem*}
The assumption that the polynomials are
essentially distinct is necessary.
 It was also shown in \cite{CF11} that there exists $d\in \N$ such that
the
 factors $\cZ_{d,T_1},\ldots,\cZ_{d,T_\ell}$ are characteristic  for pointwise convergence of the averages \eqref{E:Polynomial}.
Despite these facts, the corresponding multiple recurrence result
 (that would generalize the polynomial Szemer\'edi theorem) remains open:
 \begin{problem}\label{Conj:RecNonCommEasy}
Let $(X,\X,\mu)$ be a probability space, $T_1,\ldots, T_\ell \colon X\to X$ be invertible measure preserving transformations,
and $p_1,\ldots,p_\ell$ be distinct polynomials with zero constant term.
 Show that for every $A\in \X$ with $\mu(A)>0$ we have
$$
\mu(A\cap T_1^{m+p_1(n)}A\cap\cdots\cap T_\ell^{m+p_\ell(n)}A)>0
$$
for some  $m,n\in \N$.
\end{problem}
A solution to this problem would  imply that given a
 countable amenable group $G$ and  arbitrary elements $g_1,\ldots,g_\ell\in G$, for every  $E\subset G$
 that has  positive upper density with respect to some  F{\o}lner sequence in $G$, there exist
   $g\in G$ and $m,n\in\N$
such that
$$
g,\  g_1^{m+p_1(n)}g,\ldots,   g_\ell^{m+p_\ell(n)}g \in  E.
$$

The assumption in Problem~\ref{Conj:RecNonCommEasy} that the polynomials are distinct is necessary since,
as mentioned before, there exist  (non-commuting) transformations
$T,S$, acting on the same probability space $(X,\X,\mu)$, and a set
$A\in \X$ with $\mu(A)>0$
 such that
$\mu(T^nA\cap S^nA)=0$ for every  $n\in\N$. The multiple recurrence
property is known  when the polynomials are rationally independent \cite{FZ15}
and when all the transformations are weak
mixing   (because in this case the characteristic factors are trivial \cite{CF11}).  For general measure preserving transformations and linear polynomials
some of the simplest cases are open:
 \begin{problem21*}
Let $(X,\X,\mu)$ be a probability space  and $T, S, R \colon X\to X$ be  invertible measure preserving transformations.
 Show that for every $A\in \X$ with $\mu(A)>0$ there exist $m,n\in \N$ such that
$$
\mu(A\cap T^{m}A\cap S^{m+n}A\cap R^{m+2n}A)>0.
$$
\end{problem21*}

\section{Problems related to sequences arising from smooth functions}\label{S:Hardy}
In this section we give a list of problems related to the study of
multiple ergodic averages involving iterates given by sequences arising from smooth
functions, and related applications  to  multiple recurrence.

We are going to restrict ourselves, almost entirely,  to a class of
non-oscillatory functions that is rich enough to contain several
interesting examples. Its formal definition is the following: Let
$B$ be the collection of equivalence classes of real valued
  functions  defined on some half-line $(c,\infty)$, where we
  identify two functions if they agree eventually. A
  \emph{Hardy field} is a subfield of the ring $(B,+,\cdot)$ that is
  closed under differentiation.
  With $\H$ we denote the \emph{union of all
  Hardy fields}. It is easy to check that if a function belongs in $\H$,  then it is eventually monotonic
  and the same holds for its derivatives, so if $a\in \H$, all limits $\lim_{t\to\infty}a^{(k)}(t)$ exist (they may be infinite). We call  a \emph{Hardy sequence} any  sequence of the form $([a(n)])$ where $a\in \H$.

An explicit  example of a Hardy field to keep in mind  is the set
$\LE$ that consists of all  \emph{logarithmic-exponential functions}
(introduced by Hardy in \cite{Ha10}), meaning  all
 functions defined on some  half-line $(c,\infty)$
using  a finite combination of the symbols $+,-,\times, :,$ $\log, \exp$, operating on the real variable $t$
and on real constants.
For example, all rational functions and the functions
 $t^{\sqrt{2}}$,  $t\log{t}$,
$t^{\sqrt{\log \log t}}/\log(t^2+1)$ belong to $\LE$.
Also all {\em polynomials with fractional powers},
meaning, functions of the form  $\sum_{k=1}^d \alpha_k \, t^{c_k}$, where $\alpha_k, c_k, \in \R$, $k=1,\ldots, d$, belong to $\LE$.


 The main advantage we get by working with elements of $\H$ is
 that it is possible to relate  their growth rates with the growth rates  of their
 derivatives.
 As a consequence,  a single growth condition encodes a lot of useful information and this enables us to give
 more transparent and aesthetically pleasing  statements.
To give an example note that if $a\in \H$ and $b\in \mathcal{LE}$, then there exists a Hardy field that contains both $a$ and $b$. As a consequence, the
 limit $\lim_{t\to\infty} a'(t)/b'(t)$ exists (it may be infinite), and so assuming that $a(t),b(t)\to\infty$, we get (using L'Hospital's rule) that the quotients $ a(t)/b(t)$ and $a'(t)/b'(t)$
 have the same limit as $t\to \infty$. We deduce, for instance, that if $a\in \H$ satisfies  $a(t)/t^2\to \infty$, then $a'(t)/t\to \infty$ and $a''(t)\to \infty$.

 Background material on Hardy fields can be found in \cite{Bos82, Bos83, Bos94, Ha10, Ha12, R83a}.

\subsection{Powers of a single  transformation}
 To avoid repetition, we remark that in this subsection  we always   work with a family  $\mathcal{F}\mathrel{\mathop:}=\{a_1(t),\ldots,a_\ell(t)\}$
  of functions of \emph{polynomial growth} (meaning, for some $k\in \N$ we have
  $a_i(t)/t^k\to 0$ for $i=1,\ldots,\ell$) that belong to the same Hardy field. With
  $\text{span}^*(\mathcal{F})$ we denote the set of all \emph{non-trivial} linear combinations of elements of $\mathcal{F}$.
\subsubsection{Necessary and sufficient conditions for $\ell$-convergence}
We first state two problems from  \cite{Fr10}  related to the mean
convergence of multiple ergodic averages involving iterates given by Hardy
sequences. The following result was proved in \cite{Fr10} (the case $\ell=1$ was
first handled in \cite{BKQW05}):
\begin{theorem*}
Let $a\in \H$ have polynomial growth. Then the sequence $([a(n)])$ is good for multiple
convergence of powers
 if and only if  one of  the following conditions is satisfied:
\begin{itemize}
 \item   $|a(t)-cp(t)|/ \log t\to \infty$ for every $c\in\R$ and  $p\in \Z[t]$; \text{ or }

 \item $a(t)-cp(t)\to d$ for some $c,d\in\R$; \text{ or }

 \item $|a(t)-t/m|\leq C \log{t}$ for some $m\in\Z$ and $C>0$.
\end{itemize}
\end{theorem*}
For instance, the sequences $(n^2)$, $([n^{3/2}])$, $([n \log n])$,
$([n^2+(\log n)^2])$, $([n^2+n\sqrt{2}+\log\log n])$ are all good
for multiple convergence of powers, but the sequences $([n^2+\log
n])$, $([n^2\sqrt{2}+\log \log n])$ are not good for
$1$-convergence. Unlike the
case of polynomial sequences, if $a\in \H$ satisfies
$a(t)/t^{k-1}\to \infty$ and  $a(t)/t^k\to 0$ for some $k\in \N$, then
it can be shown that the sequence $([a(n)])$ takes odd (respectively even) values in
arbitrarily large intervals. As a consequence, when $T$ is the
rotation by $1/2$ on the circle and $f={\bf
1}_{[0,1/2]}$, the $L^2(\mu)$-limit
 $\lim_{N\to\infty}\frac{1}{|\Phi_N|}\sum_{n\in \Phi_N} T^{[a(n)]}f$
 does not exist for some appropriately chosen  F{\o}lner sequence
 $(\Phi_N)$  of subsets of $\N$.

The next problem seeks to give similar necessary and sufficient conditions for $\ell$-convergence of arbitrary collections of sequences arising from functions of polynomial growth that belong to the same Hardy field. We remind the reader that in such circumstances one is seeking to prove mean  convergence  for averages of the form
\begin{equation}\label{E:fes}
   \frac1N \sum_{n=1}^N T^{[a_1(n)]}f_1\cdots
   T^{ [a_\ell(n)]}f_\ell
 \end{equation}
 for all systems $(X,\X,\mu,T)$ and functions $f_1,\ldots, f_\ell\in L^\infty(\mu)$.

\begin{problem}\label{C:ConjConv}
Let    $\mathcal{F}$ be as above. Show that   the family of
sequences $\{([a_1(n)]),\ldots,$ $([a_\ell(n)])\}$ is good for
$\ell$-convergence of a single transformation
 if and only if  every  function  $a\in\text{span}^*(\mathcal{F})$ satisfies one of  the following conditions:
\begin{itemize}
 \item   $|a(t)-cp(t)|/ \log t\to \infty$ for every $c\in\R$ and  $p\in \Z[t]$; \text{ or }

 \item $a(t)-cp(t)\to d$ for some $c,d\in\R$; \text{ or }

 \item $|a(t)-t/m|\leq C \log{t}$ for some non-zero $m\in\Z$ and $C>0$.
\end{itemize}
\end{problem}
Mean convergence of the averages \eqref{E:fes} was
proved in \cite{Fr10} but  under much more restrictive conditions than
those advertised here.  The problem is open even when one assumes that the functions $a_1,\ldots, a_\ell$ are polynomials with fractional powers.
Also the collection of sequences
  $\{([n\log n]),([n^2\log n]),\ldots, ([n^\ell \log n])\}$ is another explicit example that is expected
   to be good for $\ell$-convergence of a single transformation but this is not known yet (not even for all weak mixing systems, or all nilsystems).

\subsubsection{Convergence to the product of the integrals}
When the multiple ergodic averages of a collection of Hardy sequences of polynomial growth
converge in the mean, one would like to have  an explicit formula
for their limit. In general,
such a limit formula can be extremely complicated, but when the sequences are in ``general position''
the limit is expected to be very simple:
\begin{problem}\label{C:ConjProduct}
Let $\mathcal{F}$ be as above and
 suppose that  for every function  $a\in \text{span}^*(\mathcal{F})$ we have $|a(t)-cp(t)|/ \log{t}\to\infty$ for every $c\in \R$ and
 $p\in \Z[t]$.
Show that  for every ergodic system
$(X,\mathcal{B},\mu,T)$
 and functions
   $f_1,\dots,f_\ell\in L^\infty(\mu)$ we have
  \begin{equation}\label{E:product'}
\lim_{N\to\infty}    \frac1N \sum_{n=1}^N T^{[a_1(n)]}f_1\cdots
   T^{ [a_\ell(n)]}f_\ell=\int f_1\, d\mu\, \cdots \, \int f_\ell\, d\mu
 \end{equation}
where the convergence takes place in $L^2(\mu)$.
\end{problem}
 The limit formula is known when  $a_i(t)=t^{c_i}$,
$i=1,\ldots,\ell$,  where $c_1,\ldots, c_\ell\in \R\setminus \Z$ are
distinct and positive \cite{Fr10} (this was established first in
\cite{BH09} when $c_i\in (0,1)$, or when the system is weak mixing). More generally,  identity
 \eqref{E:product'} is established in \cite{Fr10}  when the functions $a_1,\ldots,a_\ell$ and their pairwise differences belong to the  set $\mathcal{LE}\cap \{a\colon a(t)/t^{k+\varepsilon}\to \infty, a(t)/t^{k+1} \to 0, \text{ for some } k\geq 0 \text{ and } \varepsilon>0\}$ and $a_i(t)/a_j(t)\to 0$ or $\infty$  as $t\to \infty$ for all $i\neq j$.
Problem~\ref{C:ConjProduct} is open even when one assumes that the functions $a_1,\ldots, a_\ell$ are polynomials with fractional powers.
 If some function $a\in
\text{span}(\mathcal{F})$ satisfies the estimate $|a(t)-cp(t)|\leq C \log{t}$ for
some $c\in \R$, $C>0$,  and $p\in \Z[t]$ with $\deg(p)\geq 2$, then  one
easily sees that \eqref{E:product'} fails for  $T$  given by an
appropriate  rotation on $\T^\ell$.

An intermediate step that would help solve the previous two problems
is to find suitable characteristic factors for the relevant multiple
ergodic averages:
If $\mathcal{F}$ is as above,
$a_i(t)/ \log{t}\to\infty$,  and $(a_i(t)-a_j(t))/\log{t}\to\infty$
whenever $i\neq j$, then
for large enough $d\in \N$ the factor $\cZ_{d,T}$ is expected to be characteristic for mean convergence of the averages
\eqref{E:fes}.
This   is known when  for some $\varepsilon>0$ we have $a_i(t)/t^\varepsilon\to\infty$ and $(a_i(t)-a_j(t))/ t^\varepsilon\to\infty$
whenever
 $i\neq j$ \cite{Fr10}, and the methods of \cite{Fr10} (see the proof of Theorem~2.4 there) can be used to show that it also holds when
 $a_i(t)=ia(t)$
for $i=1,\ldots, \ell$ and $a(t)/\log{t}\to \infty$.
\subsubsection{Pointwise convergence}
Regarding variants of the identity \eqref{E:product'} that deal with pointwise convergence,  progress has been very scarce. The case
$\ell=1$ was treated in  \cite{BKQW05}, but other than this, even
the simplest cases remain open.
\begin{problem}\label{C:Conjpointwise}
Let $a,b$ be distinct positive  non-integers.  Show that for every
ergodic system $(X,\mathcal{X},\mu,T)$
 and functions
   $f, g \in L^\infty(\mu)$, we have
 $$
 \lim_{N\to\infty} \frac1N \sum_{n=1}^N f(T^{[n^a]}x) \cdot
g(T^{[n^b]}x)=\int f \ d\mu \cdot \int g\ d\mu
$$
for almost every $x\in X$.
\end{problem}
As mentioned before, mean convergence is known \cite{Fr10}. For pointwise convergence all  cases where both  $a$ and $b$ are greater than $1$ are open.

\subsubsection{Necessary conditions for $\ell$-recurrence}
Next, we state  some problems related to multiple recurrence. The following result was proved in
\cite{Fr10} (see also \cite{FrW09} for a special case):
\begin{theorem*}
Let $a\in \H$ have polynomial growth and suppose that
 $|a(t)-cp(t)|\to \infty$
for every  $c\in \R$ and  $p\in \Z[t]$. Then the sequence $([a(n)])$ is good for multiple recurrence of powers.
\end{theorem*}
 It follows that the sequences $([n^{\sqrt{2}}])$, $([n \log n])$,
$([n^2+(\log n)^2])$, $([n^2+\log n])$, $([n^2\sqrt{2}+\log \log
n])$ $([n^2+n\sqrt{2}])$ are all  good for multiple recurrence  of
powers.

Next, we  seek to give  necessary
conditions for $\ell$-recurrence of collections of
sequences arising from functions of polynomial growth that belong to
the same Hardy field. We remind the reader that in such
circumstances we seek  to prove that for every system $(X,\X,\mu,T)$ and $A\in X$ with  $\mu(A)>0$,   we
have $ \mu(A\cap T^{-[a_1(n)]}A\cap\cdots \cap T^{-[a_\ell(n)]}A)>0
$
 for some $n\in \N$.
\begin{problem}\label{C:ConjRec1}
Let $\mathcal{F}$ be as above and
 suppose that  for every function  $a\in \text{span}^*(\mathcal{F})$ we have  $|a(t)-cp(t)|\to \infty$
for every  $c\in \R$ and  $p\in \Z[t]$. Show that the collection of
sequences $\{([a_1(n)]),\ldots,([a_\ell(n)])\}$ is good for
$\ell$-recurrence of a single transformation.
\end{problem}
This is known for $\ell=1$ \cite{FrW09}.
For $\ell=2$ the problem is open even when one assumes that the sequences  are polynomials with fractional powers.


\subsubsection{Hardy sequences of super-polynomial growth}
Despite the fact that multiple recurrence and convergence properties of Hardy sequences of polynomial growth are relatively well  understood, when it comes down to sequences that grow faster than polynomials, even the most basic problems are open.
\begin{problem}
Find an example of a function $a\in \H$ that grows faster than polynomials, meaning, $a(t)/t^k \to\infty$
for every $k\in \N$,  such that the sequence
$[a(n)]$ is good for multiple recurrence and convergence of powers.
\end{problem}
The sequences  $([n^{(\log n)^{a}}])$, $([e^{n^b}])$, where $a>0$
and $b\in (0,1)$, seem to be  natural candidates; unfortunately they
are extremely hard to work with.
 Even when $\ell=1$, the relevant exponential sum estimates needed to prove
 convergence appear to be out of reach in most cases; for the first sequence such estimates are available only
 when $a\in (0,1/2)$ \cite{Kar71}, and no estimates are available for the second sequence. On the other hand,  a slower growing sequence, like the sequence
  $([n^{\log\log n}])$ may be easier to handle. But even for this sequence, $2$-recurrence and $2$-convergence is not known for all weak mixing systems or all nilsystems.

 \subsubsection{Hardy sequences evaluated at the primes}\label{SSS:HardyPrimes}
With $p_n$ we denote the $n$-th prime.
\begin{problem}
 Let  $c$ be a   positive  non-integer.
 Show that  the sequence $([p_n^c])$ is   good for  multiple recurrence
and convergence  of powers.
\end{problem}
Proving multiple recurrence   is known when $c<1$ since in this case the range of the sequence $([p_n^c])$ misses
at most finitely many positive integer values.
It is known
 that if $c$  is a positive non-integer, then  the sequence of fractional parts $(\{p_n^c\})$ is equidistributed in $[0,1]$
(see \cite{St74} or \cite{Wo75} for $c<1$ and \cite{Leit76} for
$c>1$). Probably the techniques used to prove these equidistribution
results  also gives $1$-recurrence and $1$-convergence
(it suffices to show that the sequence
$(\{p_n^c\alpha\})$ is equidistributed in the unit interval for
every non-zero $\alpha\in \R$), but
  the problem is  open for $\ell$-recurrence and $\ell$-convergence when  $\ell\geq 2$.

\subsubsection{Oscillatory sequences}
All the previous problems deal with sequences that do not oscillate. Multiple recurrence and convergence properties of oscillatory sequences are not well studied  and even
 analyzing  some simple looking sequences
leads to  very challenging   problems:
\begin{problem}\label{Conj:ConvRecOscillatory}
 Show that  the sequence $([n \sin n])$ is   good for multiple recurrence and convergence
 of powers.
\end{problem}
Quite likely one can say more; the averages $\frac{1}{N}\sum_{n=1}^N
T^{a(n)}f_1\cdot T^{2a(n)}f_2\cdot \ldots \cdot T^{\ell a(n)}f_\ell$
have the same limit in $L^2(\mu)$  when $a(n)=[n\sin n]$ and when $a(n)=n$.
This is known for $\ell=1$;  and follows from equidistribution results
in \cite{BK90} (see also related results in \cite{BBK95, BBK02}). The problem has
not been studied for $\ell\geq 2$,
even for particular classes of measure preserving systems, like
nilsystems or  weakly mixing systems.

\subsection{Commuting transformations}
 \label{SS:ConCommuting}
 As mentioned before, if a Hardy sequence has polynomial growth and stays away from constant multiples of integer polynomials, then it is going to be good for multiple recurrence and convergence of powers. The next  problem seeks to extend these results to the case of commuting transformations.

\begin{problem}\label{C:ConjCommuting2}
 Show that  if $c>1$ is not an integer, then  the sequence $([n^c])$ is good for multiple recurrence and convergence of commuting transformations.
 Moreover,  show that if $(X,\X,\mu, T_1,\ldots,T_\ell)$ is
 a system  and $f_1,\ldots,f_\ell\in L^\infty(\mu)$, then
  the $L^2(\mu)$-limit
 \begin{equation}\label{E:abra}
\lim_{N\to\infty}  \frac1N \sum_{n=1}^N T^{[n^c]}_1f_1 \cdots T^{[n^c]}_\ell f_\ell
 \end{equation}
  is equal to the $L^2(\mu)$-limit
 $\lim_{N\to\infty}  \frac1N \sum_{n=1}^N T^{n}_1 f_1\cdots T^{n}_\ell f_\ell $.
\end{problem}
Both the recurrence and the convergence problem
 is open even when $\ell=2$ and all transformations are assumed to be weak mixing.
 The result is known when the $T_1,\ldots, T_\ell$ are powers of the same transformation \cite{Fr10}, but the proof  relies crucially on the precise algebraic structure of suitable
 characteristic factors for  the corresponding multiple ergodic averages;  an advantage that we do not have when we work with general commuting transformations.
Moreover,
the methods used in \cite{Wal12} to prove mean convergence when the iterates are polynomial,
 seem hard to adjust in order to deal with  fractional powers; in any case, they give no information
 about the limit and so do not seem suitable for proving multiple recurrence results or obtaining limit formulas.

\subsection{Configurations in the primes}\label{SSS:ConPrimes}
The theorem of Szemer\'edi on
arithmetic progressions \cite{Sz75}, and its polynomial extension
\cite{BL96}, have been instrumental in proving that the primes
contain arbitrarily long arithmetic progressions \cite{GT08} and
polynomial progressions  \cite{TZ08}. Thus, it is   natural to expect
that the various known Hardy field extensions of the theorem of
Szemer\'edi \cite{Fr10, FrW09} can be used to prove  that the primes
contain the corresponding Hardy field patterns. We mention a relevant problem:
 \begin{problem}\label{C:HardyInPrimes}
 Let $\ell \in \N$ and  $c,c_1,\ldots, c_\ell$ be  positive real numbers. Show that    the  prime numbers  contain patterns of the form
 $$
 \{m,m+[n^{c}],m+2[n^{c}], \ldots, m+\ell[n^c]\} \quad \text{ and } \quad
\{m,m+[n^{c_1}],\ldots,m+[n^{c_\ell}]\}
$$
for infinitely many $n\in\N$.
\end{problem}
When all exponents are rational the existence of such patterns
follows immediately from \cite{TZ08}.

\section{Problems related to random sequences}\label{S:random}
In this section we give a list of problems related to the study of multiple ergodic averages
involving iterates given by random sequences of integers.

The  random sequences  that we work with  are
 constructed by selecting a positive  integer $n$ to be a member of our sequence
with probability $\sigma_n\in [0,1]$.
More precisely, let $(\Omega,\mathcal{F}, \mathbb{P})$ be a probability space, and  let $(X_n)_{n\in\N}$ be a sequence of independent $0-1$ valued random variables with
$$
\mathbb{P}(X_n=1)\mathrel{\mathop:}=\sigma_n  \ \text{ and } \
\mathbb{P}(X_n=0)\mathrel{\mathop:}=1-\sigma_n
 $$
 where $\sigma_n$ is a decreasing sequence of positive real numbers that satisfies $\sum_{n=1}^\infty\sigma_n=\infty$ (in which  case $\sum_{n=1}^\infty X_n(\omega)=+\infty$ almost surely).
The random sequence $(a_n(\omega))_{n\in\N}$ is constructed by
taking the positive integers $n$ for which $X_n(\omega)=1$ in
increasing order. Equivalently, $a_n(\omega)$ is the smallest $k\in
\N$ such that $X_1(\omega)+\cdots+X_k(\omega)=n$. If
$\sigma_n=n^{-a}$ for some $a\in (0,1)$, then one can show that almost surely
  $a_n(\omega)/n^{1/(1-a)}$  converges to a non-zero constant. On the other hand, if $\sigma_n=1/n$, then
  almost surely  there exists a subsequence $(n_k)$ of the integers, of density arbitrarily close to
  one, such that  the sequence $(a_{n_k}(\omega))$ is lacunary
  \cite{JLW99} (this is no longer the case if  $n\sigma_n\to\infty$). So  it makes  sense to call
  \emph{random non-lacunary sequences}
   the random sequences of integers one gets  when  $\sigma_n$
  satisfies
  $n\sigma_n\to \infty$.

We say that a property holds almost surely for the sequences $(a_n(\omega))$, if there exists a universal set $\Omega_0\in \mathcal{F}$, such that  $\mathbb{P}(\Omega_0)=1$, and for every $\omega\in \Omega_0$ the sequence $(a_n(\omega))$ satisfies the given property.

The next result was proved by  M.~Boshernitzan \cite{Bos83} for mean
convergence and by J.~Bourgain \cite{Bou88a} for pointwise
convergence (see also \cite{RW95} for a nice exposition of these
results).
\begin{theorem*}
If $n\sigma_n\to\infty$,  then almost surely the following holds:
For every   system  $(X,\X,\mu,T)$ and function
$f\in L^\infty(\mu)$, the averages
\begin{equation}\label{E:formulaAP}
 \frac{1}{N}\sum_{n=1}^N T^{a_n(\omega)}f
 \end{equation}
converge in the mean
 and their limit  equals the $L^2(\mu)$-limit of the averages  $\frac{1}{N}\sum_{n=1}^N T^{n}f$.
 Furthermore, if  $n\sigma_n/(\log\log n)^{1+\delta}\to\infty$ for some $\delta>0$, then the conclusion also holds pointwise.
\end{theorem*}
It is known that the mean convergence result fails if $\sigma_n=1/n$
 and the pointwise convergence result fails if
$\sigma_n=(\log\log n)^{1/3}/n$ (for both results see \cite{JLW99}).
It is not known whether the pointwise convergence result fails when
say $\sigma_n=\log\log n/n$.

One would naturally like to extend the previous convergence result to
  multiple ergodic averages:

\begin{problem}\label{Conj:ConvRandomComm1}
Suppose that  $n\sigma_n\to\infty$. Show that  almost surely
the sequence $(a_n(\omega))$ is good for multiple recurrence and convergence of commuting transformations.
Moreover, show that almost surely the following holds:
For every system  $(X,\X,\mu, T_1,\ldots,T_\ell)$  and functions $f_1,\ldots, f_\ell \in L^\infty(\mu)$, the averages  
 $$
 \frac{1}{N}\sum_{n=1}^N T_1^{a_n(\omega)}f_1\cdots
 T_\ell^{ a_n(\omega)}f_\ell
 $$
converge in $L^2(\mu)$
 and their limit  equals the limit of the averages
  $\frac{1}{N}\sum_{n=1}^N T_1^{n}f_1\cdots T_\ell^nf_\ell$.
\end{problem}
This problem was first mentioned  in 2004 by M.~Wierdl in \cite{Katz04}.
For $\ell\geq 2$  the result
 is known when  $\sigma_n=n^{-a}$ where
$a\in (0,1/2^{\ell-1})$ \cite{FrLW11} (the argument in \cite{Chr11} also gives this).
 When $T_1,\ldots, T_\ell$ are powers of the same transformation, a different argument
was given in \cite{FrLW14} that works for all $a\in (0,1/(\ell+1))$ (this result is superior
 to the previous one for $\ell> 3$). Under the assumption that
$n\sigma_n\to\infty$, pointwise convergence is
known when the transformations are
nilrotations  and the functions are continuous \cite{Fr09}. On the other hand, when $\sigma_n=n^{-a}$
for some $a>1/2$, $\ell=2$, and $T_2=T_1^2$, mean convergence is not
 known even when we restrict ourselves to the class of  weak mixing
systems.

A solution to Problem~\ref{Conj:ConvRandomComm1} would  imply  that   almost surely the following holds:  For every $d,\ell\in\N$,  every ${\bf v_1,\ldots,v_\ell}\in \Z^d$, every set  $\Lambda\subset\mathbb{Z}^d$ with
  $\bar{d}(\Lambda)>0$ contains patterns of the form
$$
\{{\bf m}, {\bf m}+a_n(\omega){\bf v_1},\ldots, {\bf m}+a_n(\omega){\bf v_\ell}\}
$$ for some ${\bf m}\in \Z^d$ and $n\in\N$.

Next we mention a problem that combines deterministic and random iterates.
\begin{problem}\label{Conj:ConvRandomComm2}
Suppose that  $n\sigma_n\to\infty$. Show that  almost surely the following holds:
For every system $(X,\X,\mu, T,S)$
and functions $f, g \in L^\infty(\mu)$, we have
 $$
 \lim_{N\to\infty} \frac{1}{N}\sum_{n=1}^N T^nf\cdot S^{a_n(\omega)}g=
 \E(f|\mathcal{I}_T)\cdot \E(g|\mathcal{I}_S)
 $$
where the limit is taken in $L^2(\mu)$.
Furthermore, if $\sigma_n=n^{-a}$ for some $a\in (0,1)$, show that
the convergence also holds pointwise almost everywhere.
\end{problem}
The limit formula (in $L^2(\mu)$ and pointwise)  is only known
when $a\in (0,1/14)$ \cite{FrLW11}. When $T=S$, using a different method, the eligible range was improved to $a\in (0,1/2)$ in  \cite{FrLW14}.

In the next problem we are going to work with two random sequences with different
 growth rates that are chosen independently of each other. More precisely,
 let $X_1,Y_1,X_2,Y_2,\ldots$ be a sequence of independent
 $0$-$1$ valued random variables with
$\mathbb{P}(X_n=1)\mathrel{\mathop:}=n^{-a}$ and
$\mathbb{P}(Y_n=1)\mathrel{\mathop:}=n^{-b}$ for some $a,b\in
(0,1)$. We construct the random sequence $(a_n(\omega))$
 by taking the positive integers $n$ for which
$X_n(\omega)=1$ in increasing order, and the random sequence
$(b_n(\omega))$
 by taking the positive integers $n$ for which
$Y_n(\omega)=1$ in increasing order.
\begin{problem}\label{Conj:ConvRandomComm3}
Suppose that $a,b\in(0,1)$ and $a\neq b$. Show that almost surely the following holds:
For every system $(X,\X,\mu,T,S)$ and functions $f, g
\in L^\infty(\mu)$ we have
 $$
 \lim_{N\to\infty} \frac{1}{N}\sum_{n=1}^N T^{a_n(\omega)}f\cdot S^{b_n(\omega)}g=
 \E(f|\mathcal{I}_T)\cdot \E(g|\mathcal{I}_S)
 $$
where the limit is taken in $L^2(\mu)$ or pointwise.
\end{problem}
The problem seems non-trivial even when $T=S$ is a weak mixing
transformation. Moreover, no values of $a,b\in (0,1)$ are known for which the conclusion holds.

\section{Problems related to systems with multiplicative structure}
 Let $T_n, n\in\Z,$ be  invertible measure preserving
transformations acting on a probability space $(X,\X,\mu)$ that satisfy
 $T_0=T_1=\text{id}$ and $T_m\circ T_n=T_{mn}$ for every $m,n\in \N$. We say that
 the quadruple $(X,\X,\mu,T_n)$ is a {\em measure preserving system with multiplicative structure}.
 Motivated by partition regularity problems of quadratic equations in three variables, in \cite{FH15a}
 multiple recurrence properties of systems with multiplicative structure were studied. A sample result is the following:
 \begin{theorem*}
Let $(X,\X,\mu, T_n)$ be a measure preserving system with multiplicative structure and
$A\in \X$ with  $\mu(A)>0$. Then there exist $m,n\in \N$ such that
\begin{equation}\label{E:part}
\mu\bigl(T_{m(m+n)} A\cap
T_{(m+2n)(m+3n)} A
\bigr)>0.
\end{equation}
\end{theorem*}
Using variants of the previous result one can prove that for every partition of $\N$ into
finitely many cells, there exist distinct $x$ and $y$ belonging to the same cell such that
$9x^2+16y^2=\lambda^2$ for some $\lambda\in \N$. In fact, our result proves density regularity, namely,  that such solutions can be found on every subset of $\N$ that has positive  density with respect
to any dilation invariant density on the integers. In order to prove a similar statement for the equation $x^2+y^2=\lambda^2$ (i.e., prove  partition regularity of Pythagorean pairs)
it suffices to get a positive answer for the following variant of the previous result:
\begin{problem}\label{P:Pairs}
Let $(X,\X,\mu, T_n)$ be a measure preserving system with multiplicative structure and
$A\in \X$ with  $\mu(A)>0$. Is it true that there exist $m,n\in \N$, $m> n$,  such that
$$
\mu\bigl(T_{2mn} A\cap
T_{(m-n)(m+n)} A
\bigr)>0\, ?
$$
\end{problem}
In order to prove \eqref{E:part} it turns out to be useful to analyse the limiting behavior of averages of the form
$$
\frac{1}{N^2}\sum_{1\leq m,n\leq N} \phi\big(m(m+n)\big)\, \overline{\phi}\big((m+2n)(m+3n)\big)
$$
where $\phi\colon \N\to\C$ is an arbitrary completely multiplicative function (i.e. satisfies $\phi(mn)=\phi(m)\phi(n)$ for all $m,n\in\N$) with modulus $1$. A key technical point in the proof  is that the expressions $\phi\big(m(m+n)\big)\, \overline{\phi}\big((m+2n)(m+3n)\big)$ are real and non-negative when $n=0$; this is an advantage that we do not have when we deal with similar expressions related to  Problem~\ref{P:Pairs}.

It would also be interesting to prove higher order multiple recurrence results for systems of multiplicative structure. The following   problem is a typical one:
\begin{problem}\label{P:higher}
Let $(X,\X,\mu, T_n)$ be a measure preserving system with multiplicative structure and
$A\in \X$ with  $\mu(A)>0$. Is it true that there exist $m,n\in \N$ such that
$$
\mu\bigl(T_{m(m+n)} A\cap
T_{(m+2n)(m+3n)}\cap
T_{(m+4n)(m+5n)} A
\bigr)>0\, ?
$$
\end{problem}
If the answer is positive, then this will imply partition regularity (in fact, density regularity with respect to any dilation invariant  density on the integers) for some non-trivial
quadratic equation in three variables with all three variables belonging to the same partition cell. A fundamental new difficulty  is that although
the single correlation sequences
$\mu(T_rA\cap T_s A)$ can be expressed as integral combinations of sequences of the form
$\phi(r)\, \overline{\phi}(s)$ where $\phi$ is a completely multiplicative function, this is no longer the case for the higher order correlation sequences
$\mu(T_rA\cap T_s A\cap T_tA)$. So one has to take an alternate approach.

\section{Extended bibliography organized by topic}
We give a rather extensive bibliography with material  that is  directly related to the problems
discussed before, organized by topic.
 We caution the reader that this is \emph{not} a comprehensive list of articles in ergodic Ramsey theory,
and in fact articles in several  important topics in this area
are missing from this list.  For instance, the reader will find very few articles  related to
actions  of measure preserving transformations for groups other than $\Z^d$, the richness of return times in various multiple recurrence results,   topological dynamics, and applications in partition Ramsey theory.
  There are several  excellent places to look for such topics, for instance,
   the   survey articles of V.~Bergelson \cite{Be96, Be06a, Be06b}  cover a substantial part of related  material and
   contain an extensive bibliography up to 2006.

\bigskip



\bigskip
Linear sequences: \cite{AN00,AKLR14,AbY14,  As98, As05,   As09,
As10,  As15, ADM15, AM15, AM15a, AM15b, AP12, AP13, Au09,  Au10a, Au13, Au13b,  AET11, Bere85,
Bere88,  BB84, BB86,  Be00a, BeFW06, BHK05, BJLR14, BL02,
BL04,Bert12, BFW16,BoGl09,  BGL99,    Chu10, Chu11, ChuZ11, CL84, CL88a, CL88b, CoP12,
D07, DTT08,  DT10,  DoS14, DoS15, EiKo13, EiKo16, EiT12,  Fi11a, Fo91,  Fo95, Fr04, FrK05b,
 Fu77, FuK79, FuK85, FuKO82,  FuKW90, Gr09, Gr15,
Ho09, HK01,  HK02, HK04, HK05a, HK09, HK12, HKM15, HuSY14, HuSY16,  Jen97, KiV14, Lei02b,
 Les87, Les93a, LRR03, Mc93, Mei91,  MoR16, Ro16, Ru95, Ru98, Rue09,  Ta06, Ta08,  T15a, T15b,  Th79,
 To09, Zh96,Zi99,Zi06,Zi07,Zo13,  Zo15a, Zo15b, Zo15c}.

\bigskip
Polynomial sequences: \cite{Au11c, Au11a, Au11b,  Be87a,  BBB94, BH96,
BHRF00,BL96, BL15,  BLL07, BLL08, BLM11, BeMc96, BeMc00, BeMc10, BeMo16,  BeRo15, BjB15,   Bou88a,
Bou88b, Bou88c, Bou89,Bou90, BouC16, Bu16,  Chu09,
 CF11, CFH11, DerL96,   Fi11b, FiB16, Fr08, Fr10, Fr15,   FH15a, FH15b, FH15c, FrJLW10, FrK05a, FrK06,  FrLW09, FZ15,   FuW96,
 Gr11, HafK15, HK05b, HKM14,  Joh11, Ki10, Ko16b,  Kou15a, Lei98,  Lei05c, L15, L15b,
 Les84, LyMa11, McCl11, Pan15, Po11, Th90, Wal12}.

\bigskip
Other sequences (rational, smooth, random, prime numbers, generalized
polynomials):   \cite{BH96,  BH09, BKMST14, BKS15, BL07, BLA15,  BLZ11,
BeMo13, BeMo15,  Bos83, BKQW05,  BW96, Chr11, Ei16a, EiKr14, FS10,  Fr10, Fr15a,
 FrHK07, FrHK11, FrLW11,  FrLW14,   FrW09, GT08, GT09b,
 JLW99,  Kou15b, LPWR94,  Le10,  Lei11b, Mc05, McQ09, Pan16, Su15,  TZ08, Wi88, WZ11, ZoK14}.

\bigskip
Equidistribution on nilmanifolds and other nil-stuff:  \cite{AGH63,
BL16,
CaSz11, Can16a, Can16b, Do14, DoDMSY13, EiZo13, Fr09, GT09c, GT14, GT09d,  GuMV16a, GuMV16b, GuMV16c,  HK08a, HK08b, HK11a, HKM10, HM07, HuSY13,
Ko16a,
Lei02a, Lei05a, Lei05b, Lei06, Lei07, Lei08, Lei09, Lei11a, Lei11b,
Les89, Les91, Ma51, Pa69, Pa70, Pa71, Pa73, Sh94, Sh98, ShY12, TZ15, Tu13,  Zi05}.


\bigskip
Books, survey articles, and thesis on related topics: \cite{Aa97, As03b,
AGH63, Au10b,  Be87b, Be96, Be00b, Be06a, Be06b, Chu10b, CFS82, Co82, Do15, EiWa11, EiFHN15, FrM09,
Fu81a, Fu81b, Fu90, Fu91, Fu96, Fu10, Gl03, Ha10,
Ju09,  KM10,  Katz04,  Kra06a,  Kra06b,  Kra07, Kra11,  Kre85, KN74, Le14, Ma51,  Mc99, Pe89,Qu10, RW95,
Ru90, Ta07a,  Ta07b, Ta09,  Wa82, War09, We00, Zo13b}.







\begin{thebibliography}{99}

\bibitem{Aa97} J.~Aaronson.  An introduction to infinite ergodic
theory. {\em Mathematical Surveys and Monographs} \textbf{50},
American Mathematical Society, Providence, RI, 1997.

\bibitem{AN00} J.~Aaronson, H.~Nakada.  Multiple
recurrence of Markov shifts and other infinite measure preserving
transformations. {\em   Israel J. Math.} \textbf{117} (2000),
285--310.

\bibitem{AKLR14} E.~el Abdalaoui, J.~Ku{\l}aga-Przymus, M. Lema\'nczyk, T.~de la Rue. The Chowla and the Sarnak conjectures from ergodic theory point of view.
    To appear in {\em Discrete Contin. Dyn. Syst.}, \texttt{arXiv:1410.1673}.

\bibitem{AbY14} E.~el Abdalaoui, X.~Ye.
A cubic nonconventional ergodic average with M\"obius and Liouville weight.
Preprint (2015), \texttt{arXiv:1504.00950}.

\bibitem{As98} I.~Assani.
Multiple recurrence and almost sure convergence for weakly mixing dynamical systems.
{\em Israel J. Math.} \textbf{103} (1998), 111ñ-124.

\bibitem{As03b} I.~Assani. Wiener Wintner ergodic theorems. {\em World Scientific Publishing Co., Inc.},
River Edge, NJ, 2003.

\bibitem{As05} I.~Assani.  Pointwise convergence of nonconventional averages. {\em Colloq. Math.} \textbf{102} (2005), no. 2, 245-ñ262.


\bibitem{As09} I.~Assani. Averages along cubes for not necessarily commuting
m.p.t. Ergodic theory and related fields. {\em Contemp. Math.} \textbf{430},
Amer. Math. Soc., Providence, RI, (2007), 1--19.




\bibitem{As10}
I.~Assani.
Pointwise convergence of ergodic averages along
cubes.
{\em J. Analyse Math.} {\bf 110} (2010), 241--269.


\bibitem{As15}
I.~Assani.
Pointwise double recurrence and nilsequences. Preprint (2015),
\texttt{arXiv:1504.05732}.


\bibitem{ADM15}
 I.~Assani, D.~Duncan, R.~Moore. Pointwise characteristic factors for Wiener-Wintner double recurrence theorem.  {\em Ergodic Theory Dynam. Systems} {\bf 36}  (2016),  1037--1066.

 \bibitem{AM15}
 I. Assani,  R.~Moore. Extension of Wiener-Wintner double recurrence theorem to polynomials. To appear in {\em J. Analyse Math.}, \texttt{arXiv:1409.0463}.

\bibitem{AM15a}
 I.~Assani,  R.~Moore.
A good universal weight for nonconventional ergodic averages in norm.
To appear in {\em Ergodic Theory Dynam. Systems}, \texttt{arXiv:1503.08863}.

\bibitem{AM15b}
 I.~Assani,  R.~Moore.
A good universal weight for multiple recurrence averages with commuting transformations in norm.
Preprint (2015), \texttt{arXiv:1506.06730}.

\bibitem{AP12}
 I.~Assani,  K.~Presser.
Pointwise characteristic factors for the multiterm return times theorem.
     {\em Ergodic Theory Dynam. Systems}  {\bf 32} (2012),  341--360.

\bibitem{AP13}
 I.~Assani,  K.~Presser.
A survey of the return times theorem.
Proceedings of the 2011-2012 UNC-Chapel Hill  workshops,
{\em Walter De gruyter}, (2013).


\bibitem{AGH63} L.~Auslander, L.~Green, F.~Hahn. Flows on
  homogeneous spaces.  With the assistance of L.~Markus and W.~Massey,
  and an appendix by L.~Greenberg, {\it Annals of Mathematics
    Studies}, \textbf{53}, Princeton University Press, Princeton, N.J.
  (1963).

  \bibitem{Au09} T.~Austin. On the norm convergence of nonconventional ergodic averages.
\emph{Ergodic Theory Dynam. Systems} {\bf 30} (2010), 321--338.

\bibitem{Au10a}
T.~Austin. Deducing the multidimensional Szemer\'edi theorem from an infinitary removal lemma.
{\em J. Analyse Math.} {\bf 111} (2010), 131--150.

\bibitem{Au10b}
T.~Austin. Multiple recurrence and the structure of probability-preserving systems.
{\em Phd Thesis}, UCLA (2010), \texttt{arXiv:1006.0491}.





\bibitem{Au11c} T.~Austin. Norm convergence of continuous-time polynomial multiple ergodic averages.
{\em Ergodic Theory Dynam. Systems} {\bf 32} (2012), no. 2, 361--382.

\bibitem{Au13} T.~Austin. Ergodic-theoretic implementations of the Roth density-increment argument,
{\em  Online J. Anal. Comb.} {\bf 8} (2013), 33pp.


\bibitem{Au11a}
 T.~Austin.
Pleasant extensions retaining algebraic structure, I.
  {\em J. Analyse Math.} {\bf 125} (2015), 1--36.


\bibitem{Au11b}
 T.~Austin.
Pleasant extensions retaining algebraic structure, II.
   {\em J. Analyse Math.} {\bf 126} (2015), 1--111.



\bibitem{Au13b} T.~Austin. A proof of Walsh's convergence theorem using couplings. {\em Int. Math. Res. Not.} {\bf 15} (2015),  6661--6674.

\bibitem{Au14} T.~Austin. Non-conventional ergodic averages for several commuting actions of an amenable group. To appear in  {\em J. Anal. Math.}, \texttt{arXiv:1309.4315}.

\bibitem{AET11}
T.~Austin, T.~Eisner, T.~Tao.
Nonconventional ergodic averages and multiple recurrence for von Neumann dynamical systems.
 {\em Pacific J. Math.}   {\bf 250} (2011), 1--60



\bibitem{Bere85} B.~Berend. Joint ergodicity and mixing. {\em J. Analyse Math.}
\textbf{45} (1985), 255--284.

\bibitem{Bere88} B.~Berend. Multiple ergodic theorems. {\em J. Analyse Math.}
\textbf{50} (1988), 123--142.

\bibitem{BB84}
D.~Berend, V.~Bergelson.
Jointly ergodic measure-preserving transformations.
{\em Israel J. Math.} {\bf  49} (1984), no. 4, 307-ñ314.

\bibitem{BB86}
D.~Berend, V.~Bergelson.
Characterization of joint ergodicity for noncommuting transformations.
{\em Israel J. Math.} \textbf{56} (1986), no. 1, 123--128.




\bibitem{BBK95} D.~Berend, M.~Boshernitzan,  G.~Kolesnik.
Distribution modulo $1$ of some oscillating sequences II.
{\em Israel J. Math.} \textbf{92} (1995), no. 1-3, 125--147.

\bibitem{BBK02} D.~Berend, M.~Boshernitzan,  G.~Kolesnik.
Distribution modulo 1 of some oscillating sequences III.
{\em Acta Math. Hungar.} \textbf{95} (2002), no. 1-2, 1--20.


\bibitem{BK90} D.~Berend,  G.~Kolesnik.
Distribution modulo $1$ of some oscillating sequences. {\em Israel
J. Math.} \textbf{71} (1990), no. 2, 161--179.


\bibitem{Be87a} V.~Bergelson. Weakly mixing PET. {\em Ergodic Theory
Dynam. Systems} \textbf{7} (1987), no. 3, 337--349.

\bibitem{Be87b} V.~Bergelson. Ergodic Ramsey Theory. {\em Logic and
Combinatorics} (editted by S. Simpson). {\em Contemporary Mathematics}
\textbf{ 65}, (1987), 63--87.



\bibitem{Be96} V.~Bergelson. Ergodic Ramsey Theory -- an update,
Ergodic Theory of $\Z^d$-actions (edited by M. Pollicott and K.
Schmidt). {\em London Math. Soc. Lecture Note Series}
\textbf{228} (1996), 1--61.




\bibitem{Be00a} V.~Bergelson. The multifarious Poincare recurrence theorem.
 Descriptive set theory and dynamical systems.
{\em London Math. Soc. Lecture Note Ser.} {\bf 277}, Cambridge Univ. press, Cambridge,
  (2000),
31--57.

\bibitem{Be00b} V.~Bergelson.
Ergodic theory and Diophantine problems. Topics in symbolic dynamics and applications.
{\em London Math. Soc. Lecture Note Ser.} {\bf 279}, Cambridge Univ. Press, Cambridge, (2000), 167--205.


\bibitem{Be06a} V.~Bergelson.  Combinatorial and Diophantine Applications of Ergodic Theory (with appendices by A. Leibman and by A. Quas and M. Wierdl). {\em  Handbook of Dynamical Systems}, Vol. 1B, B. Hasselblatt and A. Katok, eds., Elsevier, (2006),  745--841.

\bibitem{Be06b} V.~Bergelson. Ergodic Ramsey Theory: a dynamical approach to static theorems.
 {\em Proceedings of the International Congress of Mathematicians}, Madrid 2006,  Vol. II, 1655--1678.



\bibitem{BBB94} V.~Bergelson, M.~Boshernitzan, J.~Bourgain.
Some results on nonlinear recurrence. {\em J. Analyse Math.}
{\bf 62} (1994), 29--46.

\bibitem{BeFW06} V.~Bergelson, H.~Furstenberg, B.~Weiss.
Piecewise-Bohr sets of integers and combinatorial number theory. {\em Topics in Discrete Mathematics} (Dedicated to Jarik Nesetril on the occasion of his 60th birthday) Springer, 2006, 13--37.

\bibitem{BH96} V.~Bergelson, I.~H\aa{}land-Knutson. Sets of recurrence and
  generalized polynomials.  Convergence in ergodic theory and
  probability (Columbus, OH, 1993), {\it Ohio State Univ. Math. Res.
    Inst. Publ.}, \textbf{5}, de Gruyter, Berlin, (1996), 91--110.

\bibitem{BH09} V.~Bergelson, I.~H{\aa}land-Knutson. Weak mixing implies mixing of higher orders along tempered functions. {\em Ergodic
    Theory Dynam. Systems} {\bf 29} (2009), no. 5,  1375--1416.


\bibitem{BHM06} V.~Bergelson, I.~H\aa land-Knutson, R.~McCutcheon.
IP Systems, generalized polynomials and recurrence. {\em Ergodic
Theory Dynam. Systems} {\bf 26} (2006), no. 4, 999--1019.

\bibitem{BHK05} V.~Bergelson, B.~Host, B.~Kra, with an appendix by I. Ruzsa.
Multiple recurrence and nilsequences.
{\em Inventiones Math.} {\bf 160} (2005), no. 2, 261--303.



\bibitem{BHRF00}
V.~Bergelson, B.~Host, R.~McCutcheon, F.~Parreau. Aspects of
uniformity in recurrence. {\em Colloq. Math.}  \textbf{84/85}
(2000), no. 2, 549--576.

\bibitem{BJLR14} V.~Bergelson, A.~del Junco, M. Lema\'nczyk, J.~Rosenblatt.
Rigidity and non-recurrence along sequences.
{\em Ergodic Theory Dynam. Systems}  {\bf 34}  (2014),  1464--1502.

\bibitem{BKMST14} V.~Bergelson, G.~Kolesnik,
M.~Madritsch, Y.~Son, R.~Tichy.
Uniform distribution of prime powers and sets of recurrence and van der Corput sets in $\Z^k$.
  {\em Israel J. Math.} {\bf 201} (2014), 729--760.

\bibitem{BKS15} V.~Bergelson, G.~Kolesnik, Y.~Son.
Uniform distribution of subpolynomial functions along primes and applications.
Preprint (2015), \texttt{arXiv:1503.04960}.

\bibitem{BL96}
V.~Bergelson, A.~Leibman. Polynomial extensions of van der
Waerden's and Szemer\'edi's theorems.  {\em J. Amer. Math. Soc.}
\textbf{9} (1996), 725--753.

\bibitem{BL02}
V.~Bergelson, A.~Leibman. A nilpotent Roth theorem. {\em Inventiones
Mathematicae} {\bf 147} (2002), 429--470.


\bibitem{BL04}
V.~Bergelson, A.~Leibman. Failure of Roth theorem for solvable groups of exponential growth.
{\em Ergodic Theory Dynam. Systems} \textbf{24} (2004), no. 1, 45--53.

\bibitem{BL07}
V.~Bergelson, A.~Leibman.
Distribution of values of bounded generalized polynomials. {\em Acta Math. } {\bf 198}  (2007), 155--230.

\bibitem{BL15}
 V.~Bergelson,  A.~Leibman.
Sets of large values of correlation function for polynomial cubic configurations.
To appear in {\em Ergodic Theory Dynam. Systems}.

\bibitem{BL16}
 V.~Bergelson,  A.~Leibman.
IPr$^*$ recurrence and nilsystems. Preprint (2016), \texttt{1604.02489}.



\bibitem{BLL07} V.~Bergelson, A.~Leibman, E.~Lesigne.
Complexities of finite families of polynomials, Weyl systems, and constructions in combinatorial number theory.
{\em J. Analyse Math.} \textbf{103}  (2007), 47--92.

\bibitem{BLL08} V.~Bergelson, A.~Leibman, E.~Lesigne.
  Intersective polynomials and the polynomial Szemer\'edi theorem. {\it Adv. Math.} \textbf{219}
  (2008), no. 1, 369--388.

\bibitem{BLM11} V.~Bergelson, A.~Leibman, C.~Moreira.  From discrete to
continuous-time ergodic theorems. {\em Ergodic Theory Dynam. Systems} {\bf 32} (2012), no. 2,  383--426.

\bibitem{BLA15}
 V.~Bergelson,  A.~Leibman,  Y.~Son.
Joint ergodicity along generalized linear functions.
To appear in {\em Ergodic Theory Dynam. Systems},
\texttt{arXiv:1409.7151}.

\bibitem{BLZ11} V.~Bergelson, A.~Leibman,  T.~Ziegler.  The shifted
primes and the multidimensional Szemer\'edi and polynomial van
der Waerden Theorems. {\em Comptes Rendus Mathematique} {\bf 349} (2011), no. 3-4, 123--125.



\bibitem{BeMc96} V.~Bergelson, R.~McCutcheon.
Uniformity in polynomial Szemer\'edi theorem,  Ergodic Theory of $\Z^d$-actions (edited by M. Pollicott and K. Schmidt). {\em London Math. Soc. Lecture Note Series}
 \textbf{228} (1996),  273--296.


\bibitem{BeMc00}  V.~Bergelson, R.~McCutcheon.
 An ergodic $Ip$ polynomial Szemer\'edi theorem for commuting transformations.
{\em Mem. Amer. Math. Soc.} \textbf{146} (2000), no. 695.

\bibitem{BeMc10}  V.~Bergelson, R.~McCutcheon. Idempotent ultrafilters, multiple weak mixing and Szemer\'edi's theorem for generalized polynomials.  {\em J. Analyse Math.} \textbf{111} (2010), 77--130.

\bibitem{BeMo13}  V.~Bergelson, J.~Moreira.
Ergodic theorem involving additive and multiplicative groups of a field and \{x+y,xy\} patterns.
To appear in {\em Ergodic Theory Dynam. Systems},
\texttt{arXiv:1307.6242}.

\bibitem{BeMo15}  V.~Bergelson, J.~Moreira.
Measure preserving actions of affine semigroups and \{x+y,xy\} patterns.  To appear in {\em Ergodic Theory Dynam. Systems},
\texttt{arXiv:1509.07574}.


\bibitem{BeMo16}  V.~Bergelson, J.~Moreira.
Van der Corput's difference theorem: some moders developments. To appear in {\em Indag. Math.},
\texttt{arXiv:1510.07332}.

\bibitem{BeRo15}  V.~Bergelson, B.~Robertson.
 Polynomial multiple recurrence over rings of integers. {\em Ergodic Theory Dynam. Systems}
{\bf 36} (2016),  1354--1378.


\bibitem{Bert12} J. F.~Bertazzon. Note sur la continuit\'e de la projection dans les facteurs de Host-Kra.
 {\em C. R. Acad. Sci. Paris Sr. I Math}
 {\bf 350} (2012),  699--702.

\bibitem{BjB15} M.~Bj\"orklund, K.~Bulinski.
Twisted patterns in large subsets of $\Z^N$.
Preprint (2015), \texttt{arXiv:1512.01719}.



\bibitem{Bo81} M.~Boshernitzan. An extension of Hardy's class L
  of ``Orders of Infinity".  {\em J. Analyse Math.} \textbf{39}
  (1981), 235--255.

\bibitem{Bos82} M.~Boshernitzan. New ``Orders of Infinity".  {\em
    J. Analyse Math.} \textbf{41} (1982), 130--167.

\bibitem{Bos83} M.~Boshernitzan.
Homogeneously distributed sequences and Poincar\'e sequences of integers of sublacunary growth.
{\it Monatsh. Math.} \textbf{96} (1983), no. 3, 173--181.








\bibitem{Bos94} M.~Boshernitzan. Uniform distribution and Hardy
  fields.  {\em J. Analyse Math.} \textbf{62} (1994), 225--240.

\bibitem{BFW16} M.~Boshernitzan, N.~Frantzikinakis, M.~Wierdl.
Under recurrence in the Khintchine recurrence theorem.
Preprint (2016), \texttt{arXiv:1603.07720}.

\bibitem{BoGl09}
M.~Boshernitzan, E.~Glasner. On two recurrence problems. {\em Fund.
Math.} \textbf{206} (2009), 113--138.

\bibitem{BKQW05} M.~Boshernitzan, G.~Kolesnik, A.~Quas,
  M.~Wierdl.  Ergodic averaging sequences.  {\it J. Analyse Math.}
  \textbf{95} (2005), 63--103.

\bibitem{BW96} M.~Boshernitzan,
  M.~Wierdl.
Ergodic theorems along sequences and Hardy fields.
{\em Proc. Nat. Acad. Sci. U.S.A.}  {\bf 93} (1996), no. 16, 8205--8207.


\bibitem{Bou88a} J.~Bourgain.
 On the maximal ergodic theorem for certain subsets of the
positive integers. {\em  Israel J. Math.} {\bf 61} (1988),
39--72.


\bibitem{Bou88b} J.~Bourgain.
An approach to pointwise ergodic theorems.
{\em Lecture Notes in Math.} {\bf 1317}, Springer, Berlin, (1988), 204--223.

\bibitem{Bou88c} J.~Bourgain.
 A nonlinear version of Roth's theorem for sets of positive density in the real line. {\em  J. Analyse Math.} {\bf 50} (1988),
169--181.

\bibitem{Bou89} J.~Bourgain.
Pointwise ergodic theorems for arithmetic sets.
With an appendix by the author, H.~Furstenberg, Y.~Katznelson and D.~Ornstein.
{\em Inst. Hautes …tudes Sci. Publ. Math.} \textbf{69} (1989), 5--45.

\bibitem{Bou90} J.~Bourgain.
Double recurrence and almost sure convergence.
{\em J. Reine Angew. Math.} \textbf{404} (1990), 140--161.

\bibitem{BouC16} J.~Bourgain,  M.~Chang.
Nonlinear Roth type theorems in finite fields.
Preprint (2016),  \texttt{arXiv:1608.05448}.

\bibitem{BGL99} T.~Brown, R.~Graham, B.~Landman. On the set of common differences in van der Waerden's theorem on arithmetic progressions. {\em Canad. Math. Bull.} \textbf{42} (1999), 25--36.

\bibitem{Bu16} K.~Bulinski.
Spherical Recurrence and locally isometric embeddings of trees into positive density subsets of $\Z^d$.
Preprint (2016), \texttt{arXiv:1606.07596}.


\bibitem{CaSz11}
O.~Camarena, B.~Szegedy.
Nilspaces, nilmanifolds and their morphisms. Preprint (2010), \texttt{arXiv:1009.3825}.

\bibitem{Can16a} P.~Candela.
Notes on nilspaces: algebraic aspects. Preprint (2016), \texttt{arXiv:1601.03693}.

\bibitem{Can16b} P.~Candela.
Notes on compact nilspaces. Preprint (2016), \texttt{arXiv:1605.08940}.

\bibitem{Chr11} M.~Christ. On random multilinear operator inequalities. Preprint (2011), \texttt{arXiv:1108.5655}.


\bibitem{Chu09} Q.~Chu.
Convergence of weighted polynomial multiple ergodic averages.
{\em Proc. Amer. Math. Soc.} \textbf{137} (2009), no. 4, 1363--1369.

\bibitem{Chu10} Q.~Chu.
Convergence of multiple ergodic averages along cubes for several commuting transformations.
{\em Studia Math.} \textbf{196} (2010), no. 1, 13--22.

\bibitem{Chu10b} Q.~Chu.
Quelques problËmes de convergence et de r\'ecurrence multiple en th\'eorie ergodique.
{\em Phd Thesis},  Universit\'e Paris-Est (2010), \texttt{tel-00587631}.


\bibitem{Chu11} Q.~Chu.
Multiple recurrence for two commuting transformations.  {\em Ergodic
Theory Dynam. Systems} \textbf{31} (2011), no.3, 771--792.


\bibitem{CF11} Q.~Chu, N.~Franzikinakis. Pointwise convergence for cubic and polynomial ergodic averages of non-commuting transformations.
 {\em Ergodic Theory Dynam. Systems}  {\bf 32}  (2012), 877--897.

\bibitem{CFH11} Q.~Chu, N.~Franzikinakis, B.~Host. Ergodic averages of commuting transformations with distinct degree polynomial iterates.
{\em Proc. Lond. Math. Soc.} {\bf 102} (2011), 801--842.

\bibitem{ChuZ11} Q.~Chu, P.~Zorin-Kranich. Lower bound in the Roth theorem for amenable groups. {\em Ergodic  Theory Dynam. Systems} {\bf 35} (2015), 1746--1766.

\bibitem{CL84} J-P.~Conze, E.~Lesigne.
 Th\'eor\`emes ergodiques pour des mesures diagonales.
 {\em Bull.
Soc. Math. France} {\bf 112}  (1984), no. 2,  143--175.




\bibitem{CL88a} J-P.~Conze, E.~Lesigne.
Sur un th\'eor\`eme ergodique pour des mesures diagonales. {\em
Probabilit\'es, Publ. Inst. Rech. Math. Rennes}, 1987-1, Univ.
Rennes I, Rennes, (1988), 1--31.

\bibitem{CL88b}
J-P.~Conze, E.~Lesigne.
\newblock Sur un th\'{e}or\`{e}me ergodique pour des mesures diagonales.
\newblock {\em C. R. Acad. Sci. Paris, S\'{e}rie {I}}
\textbf{306} (1988), 491--493.

\bibitem{CFS82}
I.~Cornfeld, S.~Fomin,  Y.~Sinai.
Ergodic theory.
Translated from the Russian by A. B. Sosinskii Grundlehren der Mathematischen Wissenschaften {\em Fundamental Principles of Mathematical Sciences} \textbf{245} Springer-Verlag, New York, 1982.


\bibitem{Co82} L.~Corwin,  F.~Greenleaf.
Representations of nilpotent Lie groups and their applications. Part I.
Basic theory and examples. {\em Cambridge Studies in Advanced Mathematics} \textbf{18}, Cambridge University Press, Cambridge, 1990.

\bibitem{CoP12}
G.~Costakis, I.~Parissis. Szemer\'edi's theorem, frequent hypercyclicity and multiple recurrence.
{\em Math. Scand.} {\bf 110} (2012),  251--272.

\bibitem{D07} C.~Demeter.
Pointwise convergence of the ergodic bilinear Hilbert transform.
{\em Illinois J. Math.}  \textbf{51} (2007), no. 4, 1123--1158.


\bibitem{DTT08} C.~Demeter, T.~Tao, C.~Thiele. Maximal multi-linear operators. {\em Trans. Amer. Math. Soc.} {\bf 360}
(2008), no. 9, 4989--5042.

\bibitem{DT10} C.~Demeter,  C.~Thiele.
On the two-dimensional bilinear Hilbert transform.
{\em Amer. J. Math.} \textbf{132} (2010), no. 1, 201--256.






\bibitem{DerL96} J-M.~Derrien, E.~Lesigne. Un th\'eor\`eme ergodique polynomial ponctuel pour les endomorphismes exacts et les
K-syst\`emes. {\em Ann. Inst. H. Poincar\'e Probab. Statist.} \textbf{32} (1996), no. 6, 765--778.


\bibitem{DoDMSY13}  P.~Dong, S.~Donoso, A.~Maass, S.~Shao, X. Ye.
 Infinite-step nilsystems, independence
and complexity. {\em  Ergodic Theory Dynam. Systems}, {\bf 33} (2013), 118--143.

\bibitem{Do14} S.~Donoso.
Enveloping semigroups of systems of order $d$. {\em Discrete Contin. Dyn. Syst.} {\bf 14}
(2014),  2729--2740.

\bibitem{Do15} S.~Donoso.
Contributions to ergodic theory and topological dynamics: cube structures and automorphisms.
{\em Phd Thesis}, \texttt{tel-01271433}.



\bibitem{DoS14} S.~Donoso, W.~Sun.
A pointwise cubic average for two commuting transformations.
To appear in {\em Israel Journal of Mathematics}, \texttt{arXiv:1410.4887}.

\bibitem{DoS15} S.~Donoso, W.~Sun.
Pointwise multiple averages for systems with two commuting transformations.
Preprint (2015), \texttt{arXiv:1509.09310}.

\bibitem{DoS16} S.~Donoso, W.~Sun.
Pointwise convergence of some multiple ergodic averages.
Preprint (2016), \texttt{arXiv:1609.02529}.


\bibitem{EiWa11}
M.~Einsiedler,  T.~Ward.
Ergodic theory with a view towards number theory.
{\em Graduate Texts in Mathematics} \textbf{259}, Springer-Verlag London, Ltd., London, 2011.

\bibitem{Ei16a}
T.~Eisner.
 Nilsystems and ergodic averages along primes. Preprint (2016), \texttt{arXiv:1601.00562}.




\bibitem{EiFHN15}  T.~Eisner, B.~Farkas, M.~Haase, R.~Nagel.
Operator Theoretic Aspects of Ergodic Theory.
{\em Graduate Texts in Mathematics}, Springer, (2015).

 \bibitem{EiKo13}
T.~Eisner, D.~Kov\'acs.
On the entangled ergodic theorem. {\em Ann. Scuola Norm. Sup. Pisa Cl. Sci.} {\bf 12}  (2013), 141--156.

\bibitem{EiKo16}
T.~Eisner, D.~Kov\'acs.
On the pointwise entangled ergodic theorem.
Preprint (2015),  \texttt{arXiv:1509.05554}.



 \bibitem{EiKr14}
T.~Eisner,  B.~Krause.
 (Uniform) convergence of twisted ergodic averages. To appear in {\em  Ergodic Theory Dynam. Systems},
 \texttt{arXiv:1407.4736}.

\bibitem{EiT12}
  T.~Eisner, T.~Tao.
Large values of the Gowers-Host-Kra seminorms. {\em J. Anal. Math.} {\bf 117} (2012), 133--186.

 \bibitem{EiZo13}
  T.~Eisner, P.~Zorin-Kranich.
   Uniformity in the Wiener-Wintner theorem for nilsequences.
   {\em  Discrete Contin. Dyn. Syst.} {\bf 33} (2013), 3497--3516.





\bibitem{FS10} A.~Fan,  D.~Schneider.
Recurrence properties of sequences of integers.
{\em Sci. China Math.} {\bf 53} (2010), no. 3, 641--656.





\bibitem{Fi11b} A.~Fish. Polynomial largeness of sumsets and totally ergodic sets.
{\em Online J. Anal. Comb.} {\bf 5} (2010).

\bibitem{Fi11a} A.~Fish. Solvability of linear equations within weak mixing sets.
{\em Isr. J. Math.}  {\bf 184} (2011), 477--504.
 	

\bibitem{FiB16} A.~Fish, M. Bj\"orklund. Characteristic polynomial patterns in difference sets of matrices.
{\em Bull. London Math. Soc.} {\bf 48} (2016),  300--308.

\bibitem{Fo91} A.~Forrest. The construction of a set of recurrence which is not a set of
strong recurrence. {\em Israel J. Math.} {\bf 76}  (1991), no. 1-2, 215--228.

\bibitem{Fo95} A.~Forrest.
Two techniques in multiple recurrence.
Ergodic theory and its connections with harmonic analysis (Alexandria, 1993), 273ñ290,
{\em London Math. Soc. Lecture Note Ser.} \textbf{205}, Cambridge Univ. Press, Cambridge, 1995.

\bibitem{Fr04} N.~Frantzikinakis.
The structure of strongly stationary systems.
{\em J. Analyse Math.} \textbf{93} (2004), 359--388.



\bibitem{Fr08} N.~Frantzikinakis.  Multiple ergodic averages for
  three polynomials and applications.   {\it Trans. Amer.
    Math. Soc.}  \textbf{360} (2008), no. 10, 5435--5475.


\bibitem{Fr09} N.~Frantzikinakis. Equidistribution of sparse sequences on nilmanifolds.
 {\em J. Analyse Math.} {\bf 109} (2009), 353--395.

\bibitem{Fr10} N.~Frantzikinakis. Multiple recurrence and convergence for Hardy sequences of polynomial
growth.  {\em J. Analyse Math.}  \textbf{112} (2010), 79--135.

 \bibitem{Fr15a}
 N.~Frantzikinakis. A multidimensional Szemer\'edi theorem for Hardy sequences of polynomial growth
   {\it Trans. Amer.
    Math. Soc.}   {\bf  367}  (2015), 5653--5692.

\bibitem{Fr15}
 N.~Frantzikinakis.
Multiple correlation sequences and nilsequences. {\em Invent. Math.} {\bf 202} (2015), no. 2, 875--892.






\bibitem{FH15a}
 N.~Frantzikinakis,   B.~Host.
Higher order Fourier analysis of multiplicative functions and applications. 
 {\em J. Amer. Math. Soc.} {\bf 30},  (2017), 67--157.




\bibitem{FH15b}
 N.~Frantzikinakis,   B.~Host.
Multiple ergodic theorems for arithmetic sets. To appear in {\em  Trans. Amer.
    Math. Soc.},
\texttt{arXiv:1503.07154}.

\bibitem{FH15c}
 N.~Frantzikinakis,   B.~Host. Weighted multiple ergodic averages and correlation sequences.
To appear in
 {\em Ergodic
Theory Dynam. Systems},  \texttt{arXiv:1511.05945}.


\bibitem{FrHK07} N.~Frantzikinakis, B.~Host, B.~Kra.  Multiple
  recurrence and convergence for sets related to the primes. {\em J.
    Reine Angew. Math.} \textbf{611} (2007), 131--144.

\bibitem{FrHK11} N.~Frantzikinakis, B.~Host, B.~Kra. The polynomial multidimensional Szemer\'edi Theorem along shifted primes.
    {\em Isr. J. Math.} {\bf 194} (2013), 331--348.





\bibitem{FrJLW10} N.~Frantzikinakis, M.~Johnson, E.~Lesigne, M.~Wierdl. Powers of sequences
and convergence. {\em Ergodic Theory Dynam. Systems} \textbf{30}  (2010),
no. 5,  1431--1456.

\bibitem{FrK05a} N.~Frantzikinakis, B.~Kra.  Polynomial averages
  converge to the product of integrals. {\em Isr. J. Math.}
  \textbf{148} (2005), 267--276.




\bibitem{FrK05b} N.~Frantzikinakis, B.~Kra. Convergence of multiple ergodic averages
for some commuting transformations. {\em Ergodic Theory Dynam. Systems} \textbf{25}  (2005),
no. 3,  799--809.

\bibitem{FrK06} N.~Frantzikinakis,  B. Kra. Ergodic averages for
independent polynomials and applications. {\em J. London Math.
Soc.} \textbf{74} (2006), no. 1,  131--142.


\bibitem{FrLW06} N.~Frantzikinakis, E.~Lesigne, M.~Wierdl.
 Sets of $k$-recurrence but not $(k+1)$-recurrence.
 {\em Ann. Inst. Fourier (Grenoble)} \textbf{56} (2006), no. 4,
 839--849.

\bibitem{FrLW09} N.~Frantzikinakis, E.~Lesigne, M.~Wierdl.
  Powers of sequences and recurrence.  {\it Proc. Lond. Math. Soc.} (3)
  \textbf{98} (2009), no. 2,  504--530.



\bibitem{FrLW11} N.~Frantzikinakis, E.~Lesigne, M.~Wierdl.
Random sequences and pointwise convergence of multiple ergodic
averages. {\em Indiana Univ. Math. J.}  {\bf 61} (2012), 585--617.


\bibitem{FrLW14} N.~Frantzikinakis, E.~Lesigne, M.~Wierdl.
Random differences in Szemer\'edi's theorem and related results.
To appear in {\em J. Analyse Math.}, \texttt{arXiv:1307.1922}.

\bibitem{FrM09}  N.~Frantzikinakis, R.~McCutcheon. Ergodic Theory: Recurrence.
Survey.  {\it Encyclopedia of Complexity and Systems Science.}
Springer, (2009), Part 5, 3083--3095.



\bibitem{FrW09} N.~Frantzikinakis, M.~Wierdl. A Hardy field extension
of Szemer\'edi's theorem.  {\em Adv. Math.} \textbf{222}, (2009), 1--43.

\bibitem{FZ15} N.~Frantzikinakis, P.~Zorin-Kranich.

Multiple recurrence for non-commuting transformations along rationally independent polynomials.
{\em Ergodic Theory Dynam. Systems} {\bf 35}   (2015), 403--411.

\bibitem{Fu77} H.~Furstenberg.
Ergodic behavior of diagonal measures and a theorem of Szemer\'edi
on arithmetic progressions. {\em J. Analyse Math.} \textbf{31}
(1977), 204--256.



\bibitem{Fu81a} H.~Furstenberg.
Recurrence in ergodic theory and combinatorial number theory. {\em
Princeton University Press}, Princeton, 1981.

\bibitem{Fu81b} H.~Furstenberg.
Poincar\'e recurrence and number theory.
{\em Bull. Amer. Math. Soc. (N.S.)} \textbf{5} (1981), no. 3, 211--234.


\bibitem{Fu90} H.~Furstenberg.
Nonconventional ergodic averages. The legacy of John von Neumann (Hempstead, NY, 1988),
{\em Proc. Sympos. Pure Math.}, \textbf{50}, Amer. Math. Soc., Providence, RI, (1990), 43--56.

\bibitem{Fu91} H.~Furstenberg.
Recurrent ergodic structures and Ramsey theory. {\em
Proceedings of the International Congress of Mathematicians}, Vol. I, II (Kyoto, 1990),  Math. Soc. Japan, Tokyo, (1991), 1057--1069.


\bibitem{Fu96} H.~Furstenberg.
A polynomial Szemer\'edi theorem.  Combinatorics, Paul Erd$\ddot{\text{o}}$s is eighty, Vol. 2 (Keszthely, 1993)
{\em Bolyai Soc. Math. Stud.}, \textbf{2}, J\`anos Bolyai Math. Soc., Budapest, (1996),  253--269.


\bibitem{Fu10} H.~Furstenberg.
 Ergodic structures and non-conventional ergodic theorems.
 {\em Proceedings of the  International Congress of Mathematicians}, Hyderabad 2010, 286--298.

\bibitem{FuK79} H.~Furstenberg, Y.~Katznelson.
 An ergodic Szemer\'edi theorem for commuting transformations.
{\em J. Analyse Math.} \textbf{34} (1979), 275--291.

\bibitem{FuK85} H.~Furstenberg, Y.~Katznelson.
 An ergodic Szemer\'edi theorem for $IP$-systems
 and combinatorial theory.
{\em J. Analyse Math.} \textbf{45} (1985), 117--168.

\bibitem{FuK91} H.~Furstenberg, Y.~Katznelson.
A density version of the Hales-Jewett theorem. {\em  J. Analyse
Math.} \textbf{57} (1991), 64--119.

\bibitem{FuKO82} H.~Furstenberg, Y.~Katznelson, D.~Ornstein.
 The ergodic theoretical proof of Szemer\'edi's theorem.
 {\em  Bull.
Amer. Math. Soc. (N.S.)} {\bf 7} (1982), no. 3, 527--552.

\bibitem{FuKW90} H.~Furstenberg, Y.~Katznelson, B.~Weiss.
 Ergodic theory and configurations in sets of positive
density. Mathematics of Ramsey theory,  {\em Algorithms Combin.},
\textbf{5}, Springer, Berlin, (1990), 184--198.



\bibitem{FuW96}  H.~Furstenberg, B.~Weiss.
 A mean ergodic theorem for $(1/N)\sum\sp N\sb {n=1}$
$f(T\sp nx)$ $\ g(T\sp {n\sp 2}x)$. Convergence in ergodic theory
and probability (Columbus, OH, 1993),
 Ohio State Univ. Math. Res. Inst. Publ., {\bf 5}, de Gruyter, Berlin,
(1996), 193--227.

\bibitem{Gl03}
E.~Glasner.
Ergodic theory via joinings.
{\em Mathematical Surveys and Monographs} \textbf{101}. American Mathematical Society, Providence, RI, 2003.




\bibitem{Gow01} W.~Gowers. A new proof of Szemer\'edi's theorem.
{\em Geom. Funct. Anal.} {\bf 11} (2001), 465--588.


\bibitem{GT08a} B.~Green, T.~Tao. An inverse theorem for the Gowers  $U^3(G)$-norm.
  {\em Proc. Edinb. Math. Soc. (2)} {\bf 51} (2008), no. 1, 73--153.

\bibitem{GT08} B.~Green, T.~Tao.
The primes contain arbitrarily long arithmetic progressions.  {\it Annals Math.}
\textbf{167} (2008), 481--547.


\bibitem{GT09b} B.~Green, T.~Tao. Linear equations in primes.
 {\em  Annals  Math.} \textbf{171} (2010), 1753--1850.



\bibitem{GT09c}
B.~Green, T.~Tao.
The quantitative behaviour of polynomial orbits on nilmanifolds.
{\em Ann. of Math.}  {\bf 175}  (2012),  465--540.


\bibitem{GT14}
 B.~Green,  T.~Tao.
On the quantitative distribution of polynomial nilsequences-erratum. {\em Ann. of Math.}   {\bf 179} (2014),  1175--1183, \texttt{arXiv:1311.6170v3}.



\bibitem{GT09d} B.~Green, T.~Tao. The M\"obius function is strongly orthogonal to
nilsequences.
{\em
Ann. of Math.} {\bf 175}
 (2012),  541--566.



\bibitem{Gr09}
J.T.~Griesmer.
Ergodic averages, correlation sequences, and sumsets. {\em Phd Thesis}, Ohio State University (2009).
Available at \texttt{http://rave.ohiolink.edu/etdc/view?acc\_num=osu1243973834}.

\bibitem{Gr11} J.T.~Griesmer.
Sumsets of dense sets and sparse sets. {\em
    Israel J. Math.} {\bf 190} (2012), 229--252.

\bibitem{Gr15}
Recurrence, rigidity, and popular differences.
Preprint (2015), \texttt{arXiv:1509.03901}.

\bibitem{GuMV16a}
Y.~Gutman, F.~Manners, P.~Varj\'u.
The structure theory of Nilspaces I.
Preprint (2016), \texttt{arXiv:1605.08945}.

\bibitem{GuMV16b}
Y.~Gutman, F.~Manners, P.~Varj\'u.
The structure theory of Nilspaces II: Representation as nilmanifolds.
Preprint (2016), \texttt{arXiv:1605.08948}.

\bibitem{GuMV16c}
Y.~Gutman, F.~Manners, P.~Varj\'u.
The structure theory of Nilspaces III: Inverse limit representations and topological dynamics.
Preprint (2016), \texttt{arXiv:1605.08950}.

\bibitem{HafK15} Y.~Hafouta, Y.~Kifer.
Nonconventional Polynomial CLT. Preprint (2015), \texttt{arXiv:1504.00655}.


\bibitem{Ha10} G.~Hardy.  Orders of Infinity. The ëInfinitarcalcuulí of Paul du Bois-Reymond.
 Reprint of the 1910 edition. {\em Cambridge
    Tracts in Math. and Math. Phys.}, 12,  Hafner Publishing Co., New York, 1971.


\bibitem{Ha12} G.~Hardy. Properties of logarithmic-exponential
  functions. {\em Proc. London Math. Soc.} (2) \textbf{10} (1912),
  54--90.


\bibitem{Ho09} B.~Host. Ergodic seminorms for commuting transformations and applications.
{\emph Studia Math.} {\bf 195} (1) (2009), 31--49.

\bibitem{HK01} B.~Host,  B.~Kra.
Convergence of Conze-Lesigne averages.
{\em Ergodic Theory Dynam. Systems} \textbf{21} (2001), no. 2, 493--509.


\bibitem{HK02} B.~Host,  B.~Kra. An Odd Furstenberg-Szemer\'edi Theorem and quasi-affine systems.
 {\em J. Analyse Math.}, \textbf{86} (2002), 183--220.

\bibitem{HK04} B.~Host,  B.~Kra. Averaging along cubes. Modern dynamical systems
and applications. {\em Cambridge Univ. Press}, Cambridge, (2004), 123--144.

\bibitem{HK05a} B.~Host,  B.~Kra.
Non-conventional ergodic averages and nilmanifolds. {\em Annals
Math.}  \textbf{161}  (2005), 397--488.

\bibitem{HK05b} B.~Host, B.~Kra.  Convergence of polynomial
  ergodic averages. {\it Isr. J. Math.} \textbf{149} (2005), 1--19.

\bibitem{HK08a} B.~Host, B.~Kra.
Analysis of two step nilsequences. {\em Ann. Inst. Fourier (Grenoble)} \textbf{58} (2008), 1407--1453.

\bibitem{HK08b}  B.~Host, B.~Kra.
Parallelepipeds, nilpotent groups and Gowers norms.
{\em Bull. Soc. Math. France} {\bf  136} (2008), no. 3, 405--437.

  \bibitem{HK09} B.~Host, B.~Kra. Uniformity seminorms on $l^{\infty}$ and applications.
{\em J. Analyse Math.} \textbf{108} (2009), 219--276.

\bibitem{HK11a} B.~Host, B.~Kra. Nil-Bohr sets of integers.
{\em Ergodic Theory Dynam. Systems} {\bf 31} (2011), 113--142.

\bibitem{HK12}
B.~Host, B.~Kra.
A point of view on Gowers uniformity norms.  {\em New York J.  Math.}
{\bf 18} (2012), 213--248.





\bibitem{HKM10}
 B.~Host, B.~Kra, A.~Maass.
Nilsequences and a topological structure theorem.
 \emph{Adv. in Math.} \textbf{224} (2010),  103--129.

\bibitem{HKM14}
 B.~Host, B.~Kra, A.~Maass. Complexity of Nilsystems and systems lacking nilfactors.
         {\em J. Anal. Math.} {\bf 124} (2014), 261--295.

\bibitem{HKM15}
 B.~Host, B.~Kra, A.~Maass.
 Variations on topological recurrence. {\em  Monatsh. Math.}
        {\bf  179} (2015), Issue 1,  57--89.

\bibitem{HM07}
 B.~Host, A.~Maass.
Two step nilsystems and parallelepipeds.
{\em Bull. Soc. Math. France}  {\bf 135} (2007), no. 3, 367--405.

\bibitem{HuSY13} W.~Huang, S.~Shao,  X.~Ye.
Nil Bohr-sets and almost automorphy of higher order. {\em Memoirs of the AMS} {\bf 241}
(2016), no. 1143.

\bibitem{HuSY14} W.~Huang, S.~Shao,  X.~Ye. Pointwise convergence of multiple ergodic averages and strictly ergodic models. Preprint (2014),  \texttt{arXiv:1406.5930}.

\bibitem{HuSY16} W.~Huang, S.~Shao,  X.~Ye. Strictly ergodic models under face and parallelepiped group actions. Preprint (2016).



\bibitem{Jen97}  E.~Jenvey. Strong stationarity and De Finetti's theorem.
{\em  J.
d'Analyse Math.} \textbf{73} (1997), 1--18.

\bibitem{Joh11} M.~Johnson.
Convergence of polynomial ergodic averages for some commuting transformations.
  \emph{Illinois J. Math.} \textbf{53} (2009), no. 3, 865--882.


\bibitem{JLW99} R.~Jones, M.~Lacey, M.~Wierdl. Integer sequences with big gaps and the pointwise
ergodic theorem. {\em Ergodic Theory Dynam. Systems} \textbf{19} (1999), no.5, 1295--1308.



\bibitem{Ju09}  A.~del Junco. Ergodic theorems.
{\it Encyclopedia of Complexity and Systems Science.}
Springer, (2009), Part 5, 241--263.




\bibitem{KM10} S.~Kalikow,  R.~McCutcheon.
An outline of ergodic theory.
{\em Cambridge Studies in Advanced Mathematics}, \textbf{122}. Cambridge University Press, Cambridge, 2010.


\bibitem{KaMe78}
T.~Kamae, M.~Mend\`es-France. Van der Corput's difference theorem.
{\it Isr. J. Math.} \textbf{31} (1978), 335--342.


\bibitem{Kar71} A.~Karatsuba. Estimates for trigonometric sums by Vinogradov's method, and some
applications. {\em Trudy Mat. Inst. Steklov} \textbf{112} (1971), 241--255; English Transl.,
{\em Proc. Steklov Inst. Math.} \textbf{112} (1973), 251--265.

\bibitem{Katz01} Y.~Katznelson. Chromatic numbers of Cayley graphs on $\Z$ and
recurrence. {\em Combinatorica} \textbf{21} (2001) 211--219.



\bibitem{Katz04} Open problems. Conference in honor of Y.~Katznelson, Stanford University,  (2004).
Available at \texttt{http://mii.stanford.edu/seminars/pastevents/yk70/ik70-op.pdf}.

\bibitem{Ki10}
Y.~Kifer. Nonconventional limit theorems. {\em Probab. Theory Relat. Fields} {\bf 148} (2010),   71--106.



\bibitem{KiV14}  Y.~Kifer,  S.~Varadhan. Nonconventional limit theorems in discrete and continuous time via martingales.    {\em Ann. Probab.}
    {\bf 42}  (2014), 649--688.





\bibitem{Ko16a} J.~Konieczny.
Combinatorial properties of Nil-Bohr sets.  Preprint (2015),
\texttt{arXiv:1507.07370}.

\bibitem{Ko16b} J.~Konieczny.
Weakly mixing sets and polynomial equations. To appear in {\em Quaterly J. Math.},
\texttt{arXiv:1602.00343}.

\bibitem{Kou15a} A.~Koutsogiannis.
Integer part polynomial correlation sequences. To appear in
{\em Ergodic Theory Dynam.
Systems}, \texttt{arXiv:1512.01313}.

\bibitem{Kou15b} A.~Koutsogiannis.
Closest integer polynomial multiple recurrence along shifted primes. To appear in
{\em Ergodic Theory Dynam.
Systems}, \texttt{arXiv:1512.02264}.

\bibitem{Kra06a}  B.~Kra. The Green-Tao Theorem on arithmetic progressions in the primes: an ergodic point of view. {\em Bull. Amer. Math. Soc.} \textbf{43} (2006), 3--23.

\bibitem{Kra06b} B.~Kra.
From combinatorics to ergodic theory and back again.
{\em Proceedings of International Congress of Mathematicians}, Madrid 2006, Vol. III, 57--76.

\bibitem{Kra07}  B.~Kra.
Ergodic methods in additive combinatorics.
Additive combinatorics, 103-143, {\em CRM Proc. Lecture Notes}, \textbf{43}, Amer. Math. Soc., Providence, RI, 2007.

 \bibitem{Kra11}  B.~Kra. Poincare recurrence and number theory: thirty years later.
          {\em Bull. Amer. Math. Soc.} {\bf 48}  (2011), 497--501.

\bibitem{Kre85} U.~Krengel.
Ergodic theorems.
With a supplement by Antoine Brunel. {\em de Gruyter Studies in Mathematics}, \textbf{6}, Walter de Gruyter \& Co., Berlin, (1985).

\bibitem{Kri87} I.~Kriz.
Large independent sets in shift-invariant graphs: solution of Bergelson's problem.
{\em Graphs Combin.} \textbf{3} (1987), no. 2, 145--158.


\bibitem{KN74} L.~Kuipers, H.~Niederreiter.  Uniform distribution
  of sequences.  Pure and Applied Mathematics. {\em
    Wiley-Interscience}, New York-London-Sydney, 1974.


\bibitem{LPWR94} M.~Lacey,  K.~Petersen,  M.~Wierdl,  D.~Rudolph.
Random ergodic theorems with universally representative sequences.
{\em Ann. Inst. H. Poincare Probab. Statist.} \textbf{30} (1994), no. 3, 353--395.

\bibitem{Le10}
T.H.~L\^{e}. Intersective polynomials and the primes. {\em J. Number Theory} {\bf 130} (2010), 1705--1717.

\bibitem{Le14}
T.H.~L\^{e}.  Problems and results on intersective sets.
 Combinatorial and additive number theory,
{\em Springer Proceedings in Mathematics \& Statistics}
{\bf  101}  (2014), 115--128.







\bibitem{Lei98} A.~Leibman. Multiple recurrence theorem for measure preserving actions of a nilpotent group.
{\em Geom. Funct. Anal.}  \textbf{8} (1998), 853--931.



\bibitem{Lei02a} A.~Leibman. Polynomial mappings of groups.  {\em
    Israel J. Math.} \textbf{129} (2002), 29--60.

\bibitem{Lei02b} A.~Leibman. Lower bounds for ergodic averages.
{\em Ergodic Theory Dynam.
Systems} {\bf 22} (2002), no. 3, 863--872.

\bibitem{Lei05a} A.~Leibman.
Pointwise convergence of ergodic averages for polynomial
sequences of rotations of a nilmanifold. {\em Ergodic Theory Dynam.
Systems} {\bf 25}  (2005), no. 1,   201--213.

\bibitem{Lei05b} A.~Leibman.
Pointwise convergence of ergodic averages for polynomial actions
of $\mathbb{Z}^d$ by translations on a nilmanifold. {\em Ergodic
Theory Dynam. Systems} \textbf{25} (2005), no. 1,  215--225.

\bibitem{Lei05c} A.~Leibman.  Convergence of multiple ergodic
  averages along polynomials of several variables.  {\it Isr. J.
    Math.} \textbf{146} (2005), 303--316.

\bibitem{Lei06} A.~Leibman.
Rational sub-nilmanifolds of a compact nilmanifold.
{\em Ergodic Theory Dynam. Systems} \textbf{26} (2006), no. 3, 787--798.


\bibitem{Lei07} A.~Leibman. Orbits on a nilmanifold under the action of a polynomial sequence of translations. {\em Ergodic
Theory Dynam. Systems} \textbf{27}  (2007), no. 4, 1239--1252.

\bibitem{Lei08} A.~Leibman.
Ergodic components of an extension by a nilmanifold.
{\em Illinois J. Math.} \textbf{52} (2008), no. 3, 957--965.


\bibitem{Lei09} A.~Leibman. Orbit of the diagonal in the power of a nilmanifold. {\em Trans. Amer. Math. Soc.}
{\bf 362} (2010), 1619--1658.

\bibitem{Lei11a} A.~Leibman. Multiple polynomial sequences and nilsequences.  {\em Ergodic Theory Dynam. Systems}  {\bf 30} (2010), no. 3, 841--854.

\bibitem{Lei11b} A.~Leibman. A canonical form and the distribution of values of generalized polynomials.
 {\em Isr. J.
    Math.} {\bf 188} (2012), 131--176.

    \bibitem{L15}
 A.~Leibman.
Nilsequences, null-sequences, and multiple correlation sequences. {\em Ergodic Theory Dynam. Systems} {\bf 35} (2015), no. 1, 176--191.

\bibitem{L15b}
 A.~Leibman.
Correction to the paper ``Nilsequences, null-sequences, and multiple correlation sequences'',  \texttt{arXiv:1205.4004}.



\bibitem{Leit76} D.~Leitmann.
On the uniform distribution of some sequences.
{\em J. London Math. Soc.} (2) \textbf{14} (1976), no. 3, 430--432.

 \bibitem{Les84} E.~Lesigne. Sur la convergence ponctuelle de certaines moyennes ergodiques.
  {\em C. R. Acad. Sci. Paris S\'er. I Math.} \textbf{298} (1984), no. 17, 425--428.

\bibitem{Les87} E.~Lesigne.
Th\'eor\`emes ergodiques ponctuels pour des mesures diagonales. Cas des syst\`emes distaux.   [Pointwise ergodic theorems for diagonal measures. The case of distal systems]
{\em Ann. Inst. H. Poincar\'e Probab. Statist.} \textbf{23} (1987), no. 4, 593--612.

\bibitem{Les89} E.~Lesigne.  Th\'eor\`emes ergodiques pour une
  translation sur un nilvari\'et\'e. {\it Ergodic Theory Dynam.
    Systems} \textbf{9} (1989), no. 1, 115--126.

\bibitem{Les91} E.~Lesigne. Sur une nil-vari\'et\'e, les parties
  minimales associ\'ees \`a une translation sont uniquement
  ergodiques.  {\it Ergodic Theory Dynam. Systems} \textbf{11}
  (1991), no. 2,   379--391.

\bibitem{Les93a} E.~Lesigne.
\'Equations fonctionnelles, couplages de produits gauches et th\'eor\`emes ergodiques pour mesures diagonales.
{\it   Bull. Soc. Math. France} \textbf{121}  (1993),  no. 3, 315--351.


\bibitem{LRR03} E.~Lesigne, B.~Rittaud, T.~de la Rue. Weak disjointness of measure preserving dynamical systems. {\em Ergodic Theory Dynam. Systems} \textbf{23} (2003), 1173--1198.

\bibitem{LyMa11} N.~Lyall, A.~Magyar.
Simultaneous polynomial recurrence.
{\em Bull. Lond. Math. Soc.}  {\bf 43} (2011), Issue 4,  765--785.

\bibitem{Ma51} A.~Malcev.
On a class of homogeneous spaces.
{\em Amer. Math. Soc. Translation} \textbf{39} (1951).

\bibitem{McCl11} D.~McClendon. On the maximal Weyl complexity of families of four polynomials.
\emph{Preprint}, Available at
\texttt{http://www.swarthmore.edu/NatSci/dmcclen1/research.html}.


\bibitem{Mc93} R.~McCutcheon.
Three results in recurrence.  Ergodic theory and its connections
with harmonic analysis (Alexandria, 1993). {\em London Math. Soc.
Lecture Note Ser.}, \textbf{205}, Cambridge Univ. Press,
Cambridge, (1995), 349--358.

\bibitem{Mc99} R.~McCutcheon.
Elemental methods in ergodic Ramsey theory.
{\em Lecture Notes in Mathematics} \textbf{1722}. Springer-Verlag, Berlin, 1999.


\bibitem{Mc05} R.~McCutcheon. FVIP systems and multiple recurrence.
{\em Israel J. Math.} {\bf 146} (2005), 157--188.


\bibitem{McQ09} R.~McCutcheon, A.~Quas.
Generalized polynomials and mild mixing.
{\em Canad. J. Math.} \textbf{61} (2009), no. 3, 656--673.

\bibitem{Mei91} D.~Meiri.  Generalized correlation sequences. {\em M.Sc. Thesis},  Tel Aviv University (1991).
Available at \texttt{http://taalul.com/David/Math/ma.pdf}.

\bibitem{MoR16} J.~Moreira, F. K. Richter
A spectral refinement of the Bergelson-Host-Kra decomposition and new multiple ergodic theorems.
Preprint (2016),  \texttt{arXiv:1609.03631}.


\bibitem{Pan15} H.~Pan. Note on polynomial recurrence. Preprint (2015),  \texttt{arXiv:1502.07203}.

\bibitem{Pan16} H.~Pan.
Ergodic recurrence and bounded gaps between primes. Preprint (2016),  \texttt{arXiv:1608.04111}.


\bibitem{Pa69} W.~Parry. Ergodic properties of affine
  transformations and flows on nilmanifolds. {\it Amer. J. Math.}
  \textbf{91} (1969), 757--771.

\bibitem{Pa70} W.~Parry.
Dynamical systems on nilmanifolds.
{\em Bull. London Math. Soc.} {\bf 2} (1970), 37--40.

\bibitem{Pa71} W.~Parry.
Metric classification of ergodic nilflows and unipotent affines.
{\em Amer. J. Math.} \textbf{93} (1971), 819--828.




\bibitem{Pa73} W.~Parry.
Dynamical representations in nilmanifolds.
{\em Compositio Math.} \textbf{26} (1973), 159--174.

\bibitem{Pe89} K.~Petersen. Ergodic theory.
{\em Cambridge Studies in Advanced Mathematics}, \textbf{2}.
Cambridge University Press, Cambridge, 1989.

\bibitem{Po11} A.~Potts.
Multiple ergodic averages for flows and an application.     {\em Illinois J. Math.}
{\bf 55} (2011), 589--621.


\bibitem{Qu10} M.~Queff\'elec.
Substitution dynamical systemsóspectral analysis.
Second edition. {\em Lecture Notes in Mathematics} \textbf{1294},  Springer-Verlag, Berlin, 2010.


\bibitem{Ro15}
D.~Robertson.
Characteristic factors for multiple recurrence and combinatorial applications. {\em Phd Thesis}, Ohio State University (2015).
Available at \texttt{http://rave.ohiolink.edu/etdc/view?acc\_num=osu1429633268}.


\bibitem{Ro16}
D.~Robertson.  Characteristic factors for commuting actions of amenable groups.
To appear in {\em J. Analyse Math.},  \texttt{arXiv:1402.3843}.




\bibitem{RW95} J. Rosenblatt, M. Wierdl. Pointwise ergodic
theorems via harmonic analysis. Ergodic theory and its connections
with harmonic analysis (Alexandria, 1993). {\em London Math. Soc.
Lecture Note Ser.}, \textbf{205}, Cambridge Univ. Press,
Cambridge, (1995), 3--151.

\bibitem{R83a} M.~Rosenlicht. Hardy fields. {\em J. Math. Anal.
    Appl.} \textbf{91}   (1983), no. 2,   297--311.



\bibitem{Ru90} D.~Rudolph.
Fundamentals of measurable dynamics.
Ergodic theory on Lebesgue spaces. {\em Oxford Science Publications.} The Clarendon Press, Oxford University Press, New York, 1990

\bibitem{Ru95} D.~Rudolph.
Eigenfunctions of $T\times S$ and the Conze-Lesigne algebra. Ergodic theory and its connections with harmonic analysis (Alexandria, 1993).
{\em London Math. Soc. Lecture Note Ser.}, \textbf{205}, Cambridge Univ. Press, Cambridge, (1995),
369--432.


\bibitem{Ru98} D.~Rudolph. Fully generic sequences and a multiple-term return-times theorem. {\em Invent. Math.} {\bf  131}, 199--228.

\bibitem{Rue09}T.~de la Rue.
Notes on Austin's multiple ergodic theorem. (2009), \texttt{arXiv:0907.0538}.


\bibitem{Sh94} N.~Shah.
Limit distributions of polynomial trajectories on homogeneous spaces.
{\em Duke Math. J.}  \textbf{75} (1994), no. 3, 711--732.

\bibitem{Sh98} N.~Shah.
Invariant measures and orbit closures on homogeneous spaces for actions of subgroups generated by unipotent elements.  Lie groups and ergodic theory (Mumbai, 1996),
{\em Tata Inst. Fund. Res. Stud. Math.}, \textbf{14}, Tata Inst. Fund. Res., Bombay, (1998),  229--271.


\bibitem{ShY12} S.~Shao, X.~Ye.
Regionally proximal relation of order dd is an equivalence one for minimal systems and a combinatorial consequence.
{\em Adv. in Math.} {\bf 231} (2012), 1786--1817.



\bibitem{St74}
I.~Stux.
On the uniform distribution of prime powers.
{\em Comm. Pure Appl. Math.} \textbf{27} (1974), 729--740.


\bibitem{Su15}
W.~Sun.
Multiple recurrence and convergence for certain averages along shifted primes.
{\em Ergodic Theory Dynam. Systems} {\bf 35} (2015), no. 5, 1592--1609.

\bibitem{Sz75} E.~Szemer\'edi.
 On sets of integers containing no $k$ elements in arithmetic
progression. {\em Acta Arith.} {\bf 27} (1975), 299--345.

\bibitem{Ta06} T.~Tao. A quantitative ergodic theory proof of Szemer\'edi's theorem.
{\em Electron. J. Combin.} {\textbf 13} (2006), 1--49.

\bibitem{Ta07a} T.~Tao.
The ergodic and combinatorial approaches to Szemer\'edi's theorem.
 Additive combinatorics,
{\em CRM Proc. Lecture Notes} \textbf{43}, Amer. Math. Soc., Providence, RI, (2007), 145--193.

\bibitem{Ta07b} T.~Tao.
The dichotomy between structure and randomness, arithmetic progressions, and the primes.
{\em International Congress of Mathematicians. Vol. I}, Eur. Math. Soc., Z¸rich, (2007),  581--608.

\bibitem{Ta08} T.~Tao.
Norm convergence of multiple ergodic averages for commuting transformations.
{\em Ergodic Theory Dynam. Systems} \textbf{28} (2008), no. 2,  657--688.

\bibitem{Ta09} T.~Tao.
Poincar\'e's legacies, pages from year two of a mathematical blog. Part I.
{\em American Mathematical Society}, Providence, RI, 2009.



\bibitem{T15a}
 T.~Tao.
Deducing a weak ergodic inverse theorem from a combinatorial inverse theorem. \emph{Blog entry},
\texttt{https://terrytao.wordpress.com/2015/07/23}.


\bibitem{T15b}
 T.~Tao.
Deducing the inverse theorem for the multidimensional Gowers norms from the one-dimensional version. \emph{Blog entry},
\texttt{https://terrytao.wordpress.com/2015/07/24}.





\bibitem{TZ08} T.~Tao, T.~Ziegler.
 The primes contain arbitrarily long polynomial progressions.
 {\em Acta Math.} \textbf{201} (2008), 213--305.

\bibitem{TZ15} T.~Tao, T.~Ziegler.
Concatenation theorems for anti-Gowers-uniform functions and Host-Kra characteristic factors.
{\em Discrete Analysis} (2015), \texttt{arXiv:1603.07815}.

 \bibitem{Th79} J-P.~Thouvenot.
La demonstration de Furstenberg du th\'eor\`eme de Szemer\'edi sur les progressions arithm\'etiques.
 Seminaire Bourbaki, Vol. 1977/78, Expose No.518, {\em Lect. Notes Math.} \textbf{710} (1979), 221--232.

\bibitem{Th90} J-P.~Thouvenot.
On the almost-sure convergence of ergodic means following certain subsequences of the integers (after Jean Bourgain). Seminaire Bourbaki, Vol. 1989/90, 42\`eme ann\'ee, {\em Ast\'erisque} \textbf{189-190}, Exp. No.719, (1990) 133--153.

\bibitem{To09} H.~Towsner. Convergence of diagonal ergodic averages.
{\em Ergodic Theory Dynam. Systems} \textbf{29} (2009),  1309--1326.

\bibitem{Tu13} S.~Tu.
Nil Bohr$_0$-sets and polynomial recurrence.
{\em J. Math. Anal. Appl.} {\bf 409} (2014), 890--898.


\bibitem{Wal12}
M.~Walsh. Norm convergence of nilpotent ergodic averages. {\em Annals of Mathematics} {\bf 175} (2012), no. 3, 1667--1688.

\bibitem{Wa82} P.~Walters. An introduction to ergodic theory. {\em Graduate
Texts in Mathematics}, \textbf{79}, Springer-Verlag, New
York-Berlin, (1982).



\bibitem{War09} T.~Ward. Ergodic theory: Interactions with combinatorics and number theory.
{\it Encyclopedia of Complexity and Systems Science.}
Springer, (2009), Part 5, 3040--3053.

\bibitem{We00} B.~Weiss. Single orbit dynamics.  CBMS
Regional Conference Series in Mathematics, \textbf{95}, {\em
American Mathematical Society}, Providence, RI, (2000).


\bibitem{Wi88} M.~Wierdl.
Pointwise ergodic theorem along the prime numbers. {\em Israel J.
Math.} \textbf{64}  (1988), 315--336.


\bibitem{Wo75} D.~Wolke.
Zur Gleichverteilung einiger Zahlenfolgen.
{\em Math. Z.} \textbf{142} (1975), 181--184.

\bibitem{WZ11} T.~Wooley, T.~Ziegler.  Multiple recurrence and
convergence along the primes.  {\em Amer. J. of Math.}
{\bf 134} (2012), 1705--1732.

\bibitem{Zh96}
Q.~Zhang. On convergence of the averages
$\frac{1}{N}\sum_{n=1}^{N}f_1(R^nx) f_2(S^nx)f_3(T^nx)$. {\em
Monatsh. Math.} {\bf 122} (1996), 275--300.

\bibitem{Zi99} T.~Ziegler.
An application of ergodic theory to a problem in geometric Ramsey theory.
{\em Israel J.  Math.} \textbf{114} (1999), 271--288.

\bibitem{Zi05} T.~Ziegler.
A nonconventional ergodic theorem for a nilsystem. {\em Ergodic
Theory Dynam. Systems} \textbf{25} (2005), no. 4, 1357--1370.


\bibitem{Zi06} T.~Ziegler. Nilfactors of $\R^m$-actions and configurations
in sets of positive upper density in $\R^m$. {\em J. Analyse Math.}
\textbf{99}  (2006), 249--266.

\bibitem{Zi07} T.~Ziegler.
Universal characteristic factors and Furstenberg averages. {\em J.
Amer. Math. Soc.} \textbf{20} (2007), 53--97.


\bibitem{Zo13}
 P.~Zorin-Kranich.
Cube spaces and the multiple term return times theorem. {\em Ergodic Theory Dynam. Systems} {\bf 34} (2013),  1747--1760.

\bibitem{Zo13b} P.~Zorin-Kranich.
Ergodic theorems for polynomials in nilpotent groups. {\em Phd Thesis},
University of Amsterdam (2013), \texttt{arXiv:1309.0345}.

\bibitem{Zo15a}
 P.~Zorin-Kranich. Norm convergence of multiple ergodic averages on amenable groups.
To appear in {\em J. Analyse Math.}, \texttt{arXiv:1111.7292}.

\bibitem{Zo15b}
P.~Zorin-Kranich. A uniform nilsequence Wiener-Wintner theorem for bilinear ergodic averages.
Preprint (2015),   \texttt{arXiv:1504.04647}.

\bibitem{Zo15c}
P.~Zorin-Kranich. A double return times theorem. Preprint (2015), \texttt{arXiv:1506.05748}.

\bibitem{ZoK14}
P.~Zorin-Kranich, B.~Krause.
A random pointwise ergodic theorem with Hardy field weights. To appear in
 {\em Illinois J. Math.}, \texttt{arXiv:1410.0806}.


\end{thebibliography}
\end{document}